\documentclass[journal,10pt]{elsarticle}
\usepackage[OT1]{fontenc}  	
\usepackage[utf8]{inputenc}
\usepackage[greek,swedish,spanish,french,english,USenglish]{babel} 		
\usepackage{amssymb}  			
\usepackage{amsmath}	
\usepackage{amsthm}	
\usepackage{mathrsfs}	
\usepackage{esint}  				
\usepackage{bbm}   				
\usepackage{color,graphicx} 			
\usepackage[small]{caption}
\usepackage{cases}							
\usepackage{mathpazo}  

\newcommand{\be}{\begin{equation}}
\newcommand{\ee}{\end{equation}} 

\renewcommand{\Im}{\mathop{\rm Im}\nolimits}
\renewcommand{\Re}{\mathop{\rm Re}\nolimits}
\newcommand{\sh}{\mathop{\rm sh}\nolimits}
\newcommand{\ch}{\mathop{\rm ch}\nolimits}
\renewcommand{\th}{\mathop{\rm th}\nolimits}
\newcommand{\cth}{\mathop{\rm cth}\nolimits}
\newcommand{\tg}{\mathop{\rm tg}\nolimits}
\newcommand{\ctg}{\mathop{\rm ctg}\nolimits}
\newcommand{\arctg}{\mathop{\rm arctg}\nolimits}

\newcommand{\arcth}{\mathop{\rm arcth}\nolimits}
\newcommand{\li}{\mathop{\rm li}\nolimits}
\newcommand{\Ei}{\mathop{\rm Ei}\nolimits}

\makeatletter
\newcommand{\specialnumber}[1]{
\def\tagform@##1{\maketag@@@{(\ignorespaces##1\unskip\@@italiccorr#1)}}}
\makeatother

\makeatletter
\def\ps@pprintTitle{%
\let\@oddhead\@empty
\let\@evenhead\@empty
\let\@oddfoot\@empty
\let\@evenfoot\@oddfoot
}\makeatother

\begin{document}

\begin{frontmatter}

\title{Two series expansions for the logarithm of the gamma function involving Stirling numbers and containing only rational\\ coefficients for
certain arguments related to $\pi^{-1}$\vspace{-2mm}\footnote{\phantom{-} \\[1mm]
\texttt{\underline{Note to the readers of the 9th arXiv version:}
this is the final arXiv version of this article (first arXiv version was released on 18 August 2014) and it will not be updated in the future. 
This version differs from the journal version of the paper, which has been published by the Journal of Mathematical Analysis and Applications (Elsevier), 
vol.~442, no.~2, pp.~404-434, 2016,
DOI 10.1016/J.JMAA.2016.04.032 http://www.sciencedirect.com/science/article/pii/S0022247X16300701 In particular, the journal version
does not contain appendices, but the results contained in both versions are the same (results obtained in appendices of the arXiv version are
given in a sketch form in the journal version at the end of Section~2.2, footnotes 30-32 explain how these results were obtained).}}}

\author{Iaroslav V.~Blagouchine\corref{cor1}} 
\ead{iaroslav.blagouchine@univ-tln.fr}
\cortext[cor1]{Corresponding author. Phones: +33--970--46--28--33, +7--953--358--87--23, +7--812--744--76--79.}
\address{University of Toulon, France.\vspace{-12mm}}

\begin{abstract}
In this paper, two new series for the logarithm of the
$\Gamma$-function are presented and studied. Their polygamma analogs
are also obtained and discussed. These series involve the Stirling
numbers of the first kind and have the property to contain only
rational coefficients for certain arguments related to $\pi^{-1}$. In
particular, for any value of the form
$\ln\Gamma(\frac{1}{2}n\pm\alpha\pi^{-1})$ and
$\Psi_k(\frac{1}{2}n\pm\alpha\pi^{-1})$, where $\Psi_k$ stands for the
$k$th polygamma function, $\alpha$ is positive rational greater than
$\frac{1}{6}\pi$, $n$ is integer and $k$ is non-negative integer, these
series have rational terms only. In the specified zones of convergence,
derived series converge uniformly at the same rate as $\sum(n\ln^m\!
n)^{-2}$, where $m=1, 2, 3,\ldots$\,, depending on the order of the
polygamma function. Explicit expansions into the series with rational
coefficients are given for the most attracting values, such as
$\ln\Gamma(\pi^{-1})$, $\ln\Gamma(2\pi^{-1})$,
$\ln\Gamma(\tfrac{1}{2}+\pi^{-1})$, $\Psi(\pi^{-1})$,
$\Psi(\tfrac{1}{2}+\pi^{-1})$ and $\Psi_k(\pi^{-1})$. Besides, in this
article, the reader will also find a number of other series involving
Stirling numbers, Gregory's coefficients (logarithmic numbers, also
known as Bernoulli numbers of the second kind), Cauchy numbers and
generalized Bernoulli numbers. Finally, several estimations and full
asymptotics for Gregory's coefficients, for Cauchy
numbers, for certain generalized Bernoulli numbers and
for certain sums with the Stirling numbers are obtained. In particular, these include sharp bounds for
Gregory's coefficients and for the Cauchy numbers of the second kind.
\end{abstract}


\end{frontmatter}

\section{Introduction}
\subsection{Motivation of the study}
Numerous are expansions of the logarithm of the $\Gamma$--function and of polygamma functions
into various series. For instance
\begin{eqnarray}
&&\displaystyle\label{jxc2mxp2mxdx}
\ln\Gamma(z) \,=\, \left(z-\frac{1}{\,2\,}\right)\!\ln z -z +\frac{1}{\,2\,}\ln2\pi
+\sum_{n=1}^N \frac{{B}_{2n}}{2n(2n-1)z^{2n-1}} + \,O\big(z^{-2N-1}\big)\,,\qquad \quad
\begin{array}{l}
\displaystyle
N<\infty \\
\displaystyle
|z|\to\infty \\
\displaystyle
|\arg z|<\frac{\pi}{2}
\end{array}\qquad 
\end{eqnarray}

\begin{eqnarray}
&&\displaystyle
\ln \Gamma(z)  \,=\, -\gamma z -\ln z + \sum_{n=1}^\infty \left[\frac{z}{n}-\ln\!\left(1+\frac{z}{n}\right)\right] \,,
\qquad z\in\mathbbm{C}, \quad z\neq0, -1,-2, \ldots  \\[1mm]
&&\displaystyle
\ln\Gamma(z) \,=\, \left(z-\frac{1}{\,2\,}\right)\!\ln z -z +\frac{1}{\,2\,}\ln2\pi
+\sum_{n=0}^\infty \left[\left(z+n+\frac{1}{2} \right)\ln\frac{z+n+1}{z+n} -1\right]\,,
\qquad \quad z\neq0, -1,-2, \ldots \\[1mm]
&&\displaystyle
\ln\Gamma(z)  \,=\,-z\,(\gamma+\ln2\pi) \,-\, \frac{1}{2}\ln\frac{\sin\pi z}{\pi}\,  +\,\frac{1}{2}(\gamma+\ln2\pi)
\,+\, \frac{1}{\pi}\sum_{n=1}^{\infty}\frac{\sin 2\pi n z \cdot\ln{n}}{n}\,, \qquad 0<z<1  \\[1mm]
&&\displaystyle\label{098u43hfver}
\ln\Gamma(z)  \,=\,-\gamma z -\ln z + \sum_{n=2}^{\infty}\frac{(-1)^nz^n}{n}\zeta(n)\,, \qquad |z|<1  \\[1mm]
&&\displaystyle
\ln\Gamma(z) \,=\, \left(z-\frac{1}{\,2\,}\right)\!\ln z -z +\frac{1}{\,2\,}\ln2\pi
+\frac{1}{2}\sum_{n=1}^\infty \frac{n\cdot\zeta(n+1, z+1)}{\,(n+1)(n+2)\,}\,,\qquad\quad  \Re{z}>0\\[1mm]
&&\displaystyle\label{fckj304jdf}
\ln\Gamma(z) \,=\, \left(z-\frac{1}{\,2\,}\right)\!\ln\!\left(z-\frac{1}{\,2\,}\right) -z +\frac{1}{\,2\,}+\frac{1}{\,2\,}\ln2\pi
-\sum_{n=1}^\infty \frac{\zeta(2n, z)}{\,2^{2n+1}n(2n+1)\,}\,,\qquad\quad \Re{z}\geqslant\frac{1}{2}
\end{eqnarray}
which are respectively known as Stirling's series\footnote{This series expansion 
is one of the oldest and was known to Gauss \cite[p.~33]{gauss_02}, Euler \cite[part II, Chapter VI, \S~159, p.~466]{euler_02}, 
Stirling \cite[p.~135]{stirling_01} and De Moivre \cite{demoivre_01}, who originated it. 
Note that this series should be used with care because for $N\to\infty$, as was remarked by De Moivre himself, 
it diverges. For more information, see \cite[p.~329]{bromwich_01}, \cite[p.~111]{norlund_02}, \cite[\S~12-33]{watson_01}, 
\cite[\S~15-05]{jeffreys_02}, \cite[p.~530]{knopp_01}, \cite[p.~1]{copson_01},
\cite{kratzer_01}, \cite[\S~4.1, pp.~293--294]{olver_01}, \cite{evgrafov_03_eng}, \cite[pp.~286--288]{gelfond_01}, 
\cite[\no 6.1.40--6.1.41]{abramowitz_01}, \cite{murray_01}, 
\cite[pp.~43--50]{campbell_01}, \cite{forsyth_01}, \cite{hayman_01}, \cite{bromwich_02}, 
\cite{lanczos_01}, \cite{spouge_01}, \cite[p.~267]{comtet_01}, \cite{nemes_02}, \cite{mortici_01}, \cite{batir_01}, \cite{paris_01}.}, 
Weierstrass' series\footnote{This series follows straightforwardly from the well--known 
Weierstrass' infinite product for the $\Gamma$--function \cite[p.~10]{godefroy_01}, 
\cite[p.~12]{nielsen_01}, \cite[p.~14]{artin_01}, \cite[p.~236]{whittaker_01}, 
\cite[p.~20]{campbell_01}, \cite{kratzer_01}, \cite[pp.~21--22]{kuznecov_01}.},
Guderman's series\footnote{See e.g.~\cite[p.~111]{norlund_02}, \cite[p.~76, \no 661]{gunter_03_eng}.},
Malmsten--Kummer's series\footnote{This series is usually referred to as
\emph{Kummer's series} for the logarithm of the $\Gamma$--function, see, e.g., \cite[vol.~I, \S~1.9.1]{bateman_01}, 
\cite{whittaker_01}, \cite{srivastava_03}. However, it was
comparatively recently that we discovered that it was first obtained by Carl Malmsten and not by Ernst Kummer, see 
\cite[Sect.~2.2]{iaroslav_06}.}, Legendre's
series\footnote{See e.g.~\cite[vol.~I, eq.~1.17(2)]{bateman_01}, \cite[eq.~(21)]{wrench_01}.},
Binet's series\footnote{See e.g.\cite[vol.~I, p.~48, Eq.~(10)]{bateman_01}. To Binet are also due several 
other series whcih we discuss on p.~\pageref{ckj3w0emnd}.}
and Burnside's formula\footnote{See e.g. \cite{burnside_01}, \cite{wilton_02}, 
\cite{murray_01}, \cite[vol.~I, p.~48, Eq.~(11)]{bateman_01}. Note that in the latter reference, 
there are two errors related to Burnside's formula, i.e.~to our formula \eqref{fckj304jdf}.}  
for the logarithm of the $\Gamma$--function.\footnote{Some other series expansions for $\ln\Gamma(z)$ 
may be also found in \cite[pp.~335--339]{binet_01}, \cite[p.~1076, Eq.~(6)]{ser_01}, \cite[vol.~I, \S~1.17]{bateman_01},
\cite[pp.~240--251]{norlund_02}, see also a remark on p.~\pageref{ckj3w0emnd}. 
For further information on the $\Gamma$--function, see \cite{godefroy_01}, \cite{nielsen_01}. 
\cite[Chapt.~XII]{whittaker_01}, \cite{artin_01}, \cite{campbell_01},  \cite[vol.~I, Chapt.~I]{bateman_01}, \cite{davis_01}.}
Usually, coefficients of such expansions are either highly transcendental, or seriously suspected to be so.
Expansions into the series with rational coefficients are much less investigated, and 
especially for the gamma and polygamma functions of ``exotic'' arguments, such as, for example, $\pi^{-1}$,
logarithms or complex values.

In one of our preceding works, in exercises \no 39--49 \cite[Sect.~4]{iaroslav_06}, we have evaluated several curious integrals 
containing inverse circular and hyperbolic functions,
which led to the gamma and polygamma functions at rational multiple of $\pi^{-1}$. 
It appears that some of these integrals are particularly suitable for power series expansions.
In this paper, we derive two series expansions for the logarithm of the $\Gamma$--function,
as well as their respective analogs for the polygamma functions, by making use of such a kind of integrals.
These expansions are not simple and cannot be explicitly written in powers of $z$ up to a given order, 
but they contain rational coefficients for any argument of the form 
$z=\frac{1}{2}n\pm\alpha\pi^{-1}$, where $\alpha$ is positive rational greater than $\frac{1}{6}\pi$ and
$n$ is integer, and therefore, may be of interest 
in certain situations. As examples, we provide explicit expansions into the series with rational coefficients 
for $\ln\Gamma(\pi^{-1})$, $\ln\Gamma(2\pi^{-1})$, $\ln\Gamma(\tfrac{1}{2}+\pi^{-1})$, 
$\Psi(\pi^{-1})$, $\Psi(\tfrac{1}{2}+\pi^{-1})$ and $\Psi_k(\pi^{-1})$.
Coefficients of discovered expansions involve the Stirling numbers of the first kind, which often appear in combinatorics, as well as in various ``exotic'' 
series expansions, such as, for example
\be\label{fmse}
\begin{array}{ll}
\displaystyle
\ln\ln2 \; \, & \displaystyle= -\frac{1}{2}+\frac{5}{24}-\frac{1}{8}+\frac{251}{2880}-\frac{19}{288}
+\frac{19\,087}{362\,880}-\frac{751}{17\,280} +\frac{107\,001}{29\,030\,400} - \ldots\\[5mm]
\displaystyle &\displaystyle
= \sum_{n=1}^{\infty}\!\frac{\, (-1)^n\,}{n}\cdot
\frac{1}{\,n!\,}\!\sum_{l = 1 }^{n} \!\frac{\,\big|S_1(n, l)\big|\,}{l+1}
\end{array}
\ee
which is probably due to Arthur Cayley who gave it in 1859,\footnote{Cayley \cite{cayley_01} did not present the formula in the same form as
we did, he only
gave first six coefficients for formula \eqref{ki1d039dn321} from p.~\pageref{ki1d039dn321}. 
He noticed that the law for the formation of these coefficients ``is a very complicated one'' and that they are related in some way 
to Stirling numbers. 
The exact relationship between series \eqref{fmse} and Stirling's polynomials (and, hence, Stirling numbers) was established later by Niels Nielsen \cite[p.~76]{nielsen_01},
\cite[p.~36. eq.~(8)]{nielsen_03}, see also \cite[vol.~III, p.~257, eqs.~19.7(58)--19.7(63)]{bateman_01}. 
By the way, coefficients of this particular series are also strongly correlated with \emph{Cauchy numbers of the second kind},
see \eqref{i2039dj23r}, \eqref{j38ndbr893r}, footnote \ref{c23409jc} and \eqref{xe21fd4}.} 
or a very similar expansion converging to Euler's constant
\be\label{fms}
\begin{array}{ll}
\displaystyle
\gamma \; \, & \displaystyle= \frac{1}{2}+\frac{1}{24}+\frac{1}{72}+\frac{19}{2880}+\frac{3}{800}
+\frac{863}{362\,880}+\frac{275}{169\,344} +\frac{33\,953}{29\,030\,400} + \ldots\\[5mm]
\displaystyle &\displaystyle
= \sum_{n=1}^{\infty}\!\frac{\, (-1)^{n-1}\,}{n}\cdot
\frac{1}{\,n!\,}\!\sum_{l = 1 }^{n} \!\frac{\,S_1(n, l)\,}{l+1}
\end{array}
\ee
which was given by Lorenzo Mascheroni in 1790
\cite[p.~23]{mascheroni_01} and was subsequently rediscovered several times (in particular, by Ernst Schr\"oder in 1879 \cite[p.~115, Eq.~(25a)]{schroder_01},
by Niels E.~Nørlund in 1923 \cite[p.~244]{norlund_02},
by Jan C.~Kluyver in 1924 \cite{kluyver_02}, by Charles Jordan in 1929 \cite[p.~148]{jordan_02}, 
by Kenter in 1999 \cite{kenter_01}, by Victor Kowalenko in 2008 \cite{kowalenko_01}, \cite{kowalenko_02}).\footnote{The series itself was given by Gregorio Fontana, who, however, 
failed to find a constant to which it converges (he only proved that it should be lesser than 1). Mascheroni identified this \emph{Fontana's constant} 
and showed that it equals Euler's constant  \cite[pp.~21--23]{mascheroni_01}. Taking into account that both 
Fontana and Mascheroni did practically the equal work, series \eqref{fms} is called throughout the paper \emph{Fontana--Mascheroni's series}.
Coefficients of this series are usually witten in terms of \emph{Gregory's coefficients},
see \eqref{k2390234mrf}, \eqref{j38ndbr893r} and footnote \ref{jpqwcnqasd}.}
At large $n$ and moderate values of argument $z$, discovered series converge approximately at the same rate as 
$\sum\, (n\ln^{m} n)^{-2}\,$, where $m=1$ for $\ln\Gamma(z)$ and $\Psi(z)$, 
$m=2$ for $\Psi_1(z)$ and $\Psi_2(z)$,
$m=3$ for $\Psi_3(z)$ and $\Psi_4(z)$, \emph{etc.}
At the same time, first partial sums of the series may behave
quite irregularly, and the sign of its general term changes in a complex pattern.
However, in all cases, including small and moderate values of $n$, the absolute value of the $n$th general term remains bounded by
$\alpha n^{-2}$, where $\alpha$ does not depend on $n$.
Finally, in the manuscript, we also obtain a number of other series expansions 
containing Stirling numbers, Gregory's coefficients (logarithmic numbers, also known as Bernoulli numbers of the second kind), 
Cauchy numbers, ordinary and generalized Bernoulli numbers,
binomial coefficients and harmonic numbers, as well as provide the convergence analysis for some of them.

\subsection{Notations}\label{notations}
Throughout the manuscript, following abbreviated notations are used: $\,\gamma=0.5772156649\ldots~$ for
Euler's constant, 
$\binom{k}{n}$ denotes the binomial coefficient $C^n_k$, 
${B}_n$ stands for the $n$th Bernoulli number,\footnote{In particular 
${B}_0=+1$, ${B}_1=-\frac{1}{2}$, ${B}_2=+\frac{1}{6}$, 
${B}_3=0$, ${B}_4=-\frac{1}{30}$, ${B}_5=0$, ${B}_6=+\frac{1}{42}$, ${B}_7=0$,
${B}_8=-\frac{1}{30}$, ${B}_9=0$, ${B}_{10}=+\frac{5}{66}$, ${B}_{11}=0$, ${B}_{12}=-\frac{691}{2730}$, \emph{etc}., 
see \cite[Tab.~23.2, p.~810]{abramowitz_01}, \cite[p.~5]{krylov_01} or \cite[p.~258]{gelfond_01} for further values. Note also 
that some authors may use slightly different definitions for the Bernoulli numbers, see e.g.~\cite[p.~19, \no 138]{gunter_03_eng}
or \cite[pp.~3--6]{arakawa_01}.}
$\lfloor x\rfloor$ 
for the integer part of $x$, $\tg z$ for the tangent of $z$, 
$\ctg z$ for the cotangent of $z$, $\ch z$ for the hyperbolic cosine of $z$, $\sh z$ for the hyperbolic sine of $z$, 
$\th z$ for the hyperbolic tangent of $z$, $\cth z$ for the hyperbolic cotangent of $z$.
In order to avoid any confusion between compositional inverse and multiplicative inverse, 
inverse trigonometric and hyperbolic functions are denoted 
as $\arccos$, $\arcsin$, $\arctg, \ldots~$ and not as $\cos^{-1}$,
$\sin^{-1}$, $\tg^{-1}, \ldots~$
Writings $\Gamma(z)$, $\Psi(z)$, $\Psi_1(z)$, $\Psi_2(z)$, $\Psi_3(z)$, $\Psi_4(z)$, $\zeta(z)$ and $\zeta(z,v)$ denote
respectively the gamma, the digamma,
the trigamma, the tetragamma, the pentagamma, the hexagamma, the Riemann zeta and the Hurwitz zeta functions of argument $z$. 
The Pochhammer symbol $(z)_n$, which is also known as the generalized factorial function, is defined as the rising factorial
$(z)_n\equiv z(z+1)(z+2)\cdots (z+n-1)=\Gamma(z+n)/\Gamma(z)$.\footnote{For nonpositive and complex $n$, only the latter 
definition $(z)_n=\Gamma(z+n)/\Gamma(z)$ holds.}\up{,}\footnote{Note that some writers (mostly German--speaking)
call such a function \emph{faculté analytique} or \emph{Facultät}, see e.g.~\cite{schlomilch_04}, \cite[p.~186]{schlomilch_05}, 
\cite[vol.~II, p.~12]{schlomilch_06}, \cite[p.~119]{hagen_01}, \cite{kramp_01}. Other names and notations
for $(z)_n$ are briefly discussed in \cite[pp.~45--47]{jordan_01} and in \cite[pp.~47--48]{knuth_01}.} For sufficiently large $n$, not necessarily integer,
the latter can be given by this useful approximation
\be\label{lk2093mffmnjw}
\begin{array}{ll}
\displaystyle
(z)_n \;&\displaystyle =\,\frac{\,n^{n+z-\frac{1}{2}}\sqrt{2\pi} \,}{\Gamma(z)\,e^{n}} 
\left\{1+ \frac{\,6 z^2 - 6z + 1\,}{12 n} + \frac{\,36 z^4 - 120 z^3 + 120 z^2  - 36 z  + 1}{288 n^2} + O(n^{-3})\right\}\\[6mm]
&\displaystyle
\,=\,\frac{\,n^z\cdot \Gamma(n)\,}{\Gamma(z)}\left\{1+ \frac{\,z(z-1)\,}{2 n}
+ \frac{\,z(z-1)(z-2)(3z-1)\,}{24 n^2} +O(n^{-3})\right\}
\end{array}
\ee
which follows from the Stirling formula for the $\Gamma$--function.\footnote{A simpler 
variant of the above formula may be found in \cite{tricomi_01}, \cite{kratzer_01}.}
Writing $S_1(k,n)$ stands for the signed Stirling numbers of the first kind (see Sect.~\ref{jgf76fce}).
Kronecker symbol (Kronecker delta) of arguments $l$ and $k$ is denoted by $\,\delta_{l,k}\,$: $\,\delta_{l,k}=1\,$
if $l=k$ and $\,\delta_{l,k}=0\,$ if $l\neq k$.
$\Re{z}$ and $\Im{z}$ denote respectively real and imaginary parts 
of $z$. Natural numbers are defined in a traditional way as a set of positive integers, which is denoted by $\mathbbm{N}$.
Letter $i$ is never used as index and is $\sqrt{-1\,}$. Finally, by the relative error between the quantity $A$
and its approximated value $B$, we mean either $(A-B)/A$, or $|A-B|/|A|$, depending on the context.
Other notations are standard.

\section{Stirling numbers and their role in various expansions}\label{stirling}
\subsection{General information in brief}\label{jgf76fce}
Stirling numbers were introduced by the Scottish mathematician James Stirling in his famous treatise \cite[pp.~1--11]{stirling_01}, and 
were subsequently rediscovered in various forms by numerous authors, including Christian Kramp, Pierre--Simon Laplace,
Andreas von Ettingshausen, Ludwig Schläffli, Oskar Schl\"omilch, Paul Appel,
Arthur Cayley, George Boole, 
James Glaisher, Leonard Carlitz and many others \cite{hindenburg_01}, 
\cite{kramp_01}, \cite[Book I, part I]{laplace_02}, \cite{ettingshausen_01}, \cite{schlaffli_01},
\cite{schlaffli_02}, \cite{schlomilch_04}, 
\cite[pp.~186--187]{schlomilch_05}, \cite[vol.~II, pp.~23--31]{schlomilch_06}, \cite{appel_01},
\cite{cayley_00}, \cite{cayley_01}, \cite{cayley_02}, \cite{boole_01}, \cite{glaisher_02}, \cite[p.~129]{carlitz_02},
\cite[pp.~67--78]{nielsen_01}, \cite[p.~1]{jordan_01}, \cite{gould_03}, \cite{knuth_02}.\footnote{Although it is largely
accepted that these numbers were introduced by James Stirling in his famous treatise \cite{stirling_01} published in 1730,
Donald E.~Knuth \cite[p.~416]{knuth_02} rightly remarks that these numbers may be much older. In particular, the above--mentioned
writer found them in an old unpublished manuscript of Thomas Harriot, dating about 1600.}
Traditionally, Stirling numbers are devided
in two different ``kinds'': Stirling numbers of the first kind and those of the second kind, albeit there really is only one ``kind'' of
Stirling numbers [see the remark after Eq.~\eqref{io20323m3e}].\footnote{Within the 
framework of our study we are not concerned with the \emph{Stirling 
numbers of the second kind}; we, therefore, will not treat them here. By the way, it is interesting that Stirling himself,
first, introduced numbers of the second kind \cite[p.~8]{stirling_01}, and then, those of the first kind \cite[p.~11]{stirling_01}.} 
The Stirling numbers of the first kind appear in numerous occasions in
combinatorics, in calculus of finite differences,  in numerical analysis, in number theory and even in calculus of variations.
In combinatorics, Stirling numbers of the first kind, denoted $|S_1(n,l)|$, are defined
as the number of ways to arrange $n$ objects into $l$ cycles or cyclic arrangements ($|S_1(n,l)|$ is often verbalized ``$n$ cycle $l$''). 
These numbers are also called 
\emph{unsigned} (or \emph{signless}) \emph{Stirling numbers}, as opposed to $S_1(n,l)$, which are called \emph{signed Stirling numbers} and which 
are related to the former as $S_1(n,l)=(-1)^{n\pm l}|S_1(n,l)|$.
In the analysis and related disciplines, 
the unsigned/signed Stirling numbers of the first kind are usually defined as the coefficients in the expansion of rising/falling factorial
\be\label{x2l3dkkk03d}
\specialnumber{a,b}
\begin{cases}
\displaystyle 
\prod_{k=0}^{n-1} (z+k) 
\,=\,(z)_n\,=\,\frac{\Gamma(z+n)}{\Gamma(z)}\,=\,\sum_{l=1}^n |S_1(n,l)|\cdot z^l \\[5mm]
\displaystyle 
\prod_{k=0}^{n-1} (z-k) 
\,=\,(z-n+1)_n\,=\frac{\Gamma(z+1)}{\Gamma(z+1-n)}\,=\,\sum_{l=1}^n S_1(n,l)\cdot z^l 
\end{cases}
\ee
where $z\in\mathbbm{C}$ and $n\geqslant1$.
Stirling numbers of the first kind are also often introduced via their generating functions
\be\label{h37qp2b237ds}
\specialnumber{a,b}
\begin{cases}
\displaystyle
\sum_{n=l}^\infty\!\frac{|S_1(n,l)|}{n!}z^n =\,(-1)^l\frac{\ln^l(1-z)}{l!}\,, \qquad & l=0, 1, 2, \ldots \\[5mm]
\displaystyle
\sum_{n=l}^\infty\!\frac{S_1(n,l)}{n!}z^n =\,\frac{\ln^l(1+z)}{l!}\,, \qquad & l=0, 1, 2, \ldots
\end{cases}
\ee
both series on the left being uniformly and absolutely convergent on the disk $|z|<1$.\footnote{Remark that formally,
in \eqref{h37qp2b237ds}, 
the summation may be started not only from $n=l$, but from any $n$ in the range $[0,l-1]$, because $S_1(n,l)=0$ for such $n$.}
Signed Stirling numbers of the first kind may be calculated explicitly via the following formula 
\be\label{io20323m3e}
S_1(n,l)\,=\,
\begin{cases}
\displaystyle
\frac{(2n-l)!}{(l-1)!}
\sum_{k=0}^{n-l}\frac{1}{(n+k)(n-l-k)!(n-l+k)!}
\sum_{r=0}^{k}\frac{(-1)^{r} r^{n-l+k} }{r!(k-r)!} \,,\quad& l\in[1,n] \\[3mm]
1 \,,& n=0\,,\, l=0 \\[1mm] 
0 \,,& \text{otherwise} 
\end{cases}
\ee
where $\,S_1(0,0)=1\,$ by convention.\footnote{In the above, we always supposed that $n$ and $l$ are nonnegative,
although, this, strictly speaking, is not  necessary. 
In fact, for negative  arguments $n$ and $l$, Stirling numbers of the first kind
reduce to those of the second kind and \emph{vice--versa}, see e.g.~\cite[p.~412]{knuth_02}, 
\cite[p.~116]{gould_02}, \cite[p.~60 \emph{et seq.}]{hagen_01}.} From the above definitions, it is visible that numbers $S_1(n,l)$ are necessarily integers:
$S_1(1,1)=+1$, $S_1(2,1)=-1$, $S_1(2,2)=+1$, $S_1(3,1)=+2$, $S_1(3,2)=-3$, $S_1(3,3)=+1$, \ldots\,,
$S_1(8,5)=-1\,960$, \ldots\,, $S_1(9,3)=+118\,124$, \emph{etc.}

Stirling numbers of the first kind were studied in a large number of works and have many various properties
which we cannot describe in a small article. 
These numbers are of great utility, especially, for the summation of series, the fact which was noticed primarily by Stirling in
his marvellous treatise \cite{stirling_01} and which was later emphasized by numerous writers.
In particular, Charles Jordan, who worked a lot on Stirling numbers, see e.g.~\cite{jordan_01}, \cite{jordan_02}, \cite{jordan_00}, remarked that these
numbers may be even more important than Bernoulli numbers. 
In what follows we give only a small amount 
of the information necessary for the understanding of the rest of our work. 
Readers interested in a more deep study of these numbers are kindly invited to refer to the above--cited historical references, as well as to
the following specialized literature:
\cite[Chapt.~IV]{jordan_01}, \cite{jordan_02}, \cite{jordan_00}, \cite{nielsen_04}, \cite[pp.~67--78]{nielsen_01}, \cite{nielsen_03},
\cite{tweedie_01}, 
\cite[Sect.~6.1]{knuth_01}, \cite[pp.~410--422]{knuth_02}, \cite[Chapt.~V]{comtet_01}, \cite{dingle_01}, 
\cite[Chapt.~4, \S~3, \no 196--210]{polya_01_eng}, \cite[p.~60 \emph{et seq.}]{hagen_01}, \cite{netto_01},
\cite[p.~70 \emph{et seq.}]{riordan_01}, \cite[vol.~1]{stanley_01},
\cite{bender_01}, \cite[Chapt.~8]{charalambides_01}, \cite[\no 24.1.3, p.~824]{abramowitz_01}, \cite[Sect.~21.5-1, p.~824]{korn_01}, 
\cite[vol.~III, \S~19.7]{bateman_01}, \cite{norlund_02}, \cite{steffensen_02}, \cite[pp.~91--94]{conway_01}, \cite[pp.~2862--2865]{weisstein_04}, 
\cite[Chapt.~2]{arakawa_01}, \cite{mitrinovic_01}, \cite{gould_01}, \cite{gould_02}, \cite{gould_03},
\cite{wachs_01}, \cite{carlitz_02}, \cite{carlitz_03}, \cite[p.~642]{olson_01}, \cite{salmieri_01}, \cite{gessel_01},
\cite{wilf_01}, \cite{moser_01}, \cite{bellavista_01}, \cite{wilf_02}, \cite{temme_02}, 
\cite{howard_01}, \cite{butzer_02}, \cite{butzer_01}, \cite{hwang_01}, \cite{adamchik_03},
\cite{timashev_01}, \cite{grunberg_01}, \cite{louchard_01}, 
\cite{shen_01}, \cite{shirai_01}, \cite{sato_01},
\cite{rubinstein_01}, \cite{rubinstein_02}, \cite{hauss_01}, \cite{skramer_01}.
Note that many writers discovered
these numbers independently, without realizing that they deal with the Stirling numbers.
For this reason, in many sources, these numbers may appear under different names, different notations and even slightly different definitions.
Actually, only in the beginning of the XXth century, the name ``Stirling numbers'' appeared in mathematical literature
(mainly, thanks to Thorvald N.~Thiele
and Niels Nielsen \cite{nielsen_04}, \cite{tweedie_01}, \cite[p.~416]{knuth_02}).
 Other names for these numbers include: \emph{factorial coefficients}, \emph{faculty's coefficients} (\emph{Facultätencoefficienten},
\emph{coefficients de la faculté analytique}), \emph{differences of zero} and even \emph{differential coefficients of nothing}.  
The Stirling numbers are also closely connected to the \emph{generalized Bernoulli numbers} $B^{(s)}_n$, also known as 
\emph{Bernoulli numbers of higher order}, see e.g.~\cite[p.~129]{carlitz_02}, \cite[p.~449]{gould_01}, \cite[p.~116]{gould_02},
\cite[vol.~III, \S~19.7]{bateman_01}, \cite{young_02}, \cite{brychkov_01}, \cite{brychkov_02}; 
many of their properties may be, therefore, deduced from those of $B^{(s)}_n$. 
As concerns notations, there exist more than 50 notations for them, 
see e.g.~\cite{gould_02}, \cite[pp.~vii--viii, 142, 168]{jordan_01}, 
\cite[pp.~410--422]{knuth_02}, \cite[Sect.~6.1]{knuth_01}, and we do not insist on our particular notation, which may
seem for certain not properly chosen.
Lastly, we remark that there also are several slightly different definitions of the Stirling numbers of the first kind;
our definitions \eqref{x2l3dkkk03d}--\eqref{io20323m3e} agree with those
adopted by Jordan \cite[Chapt.~IV]{jordan_01}, \cite{jordan_02}, \cite{jordan_00}, Riordan \cite[p.~70 \emph{et seq.}]{riordan_01},
Mitrinovi\'c \cite{mitrinovic_01}, Abramowitz \& Stegun \cite[\no 24.1.3, p.~824]{abramowitz_01} and many others.\footnote{Modern CAS, such as \textsl{Maple} or \textsl{Mathematica}, also share
these definitions; in particular \texttt{Stirling1(n,l)} in the former and \texttt{StirlingS1[n,l]} in the latter correspond to our $S_1(n,l)$.} 
A quick analysis of several alternative definitions may be found in \cite{gould_02}, \cite{gould_03}, \cite[pp.~vii--viii and Chapt.~IV]{jordan_01},
\cite[pp.~410--422]{knuth_02}.

\subsection{MacLaurin series expansions of certain composite functions and some other series with Stirling numbers}\label{lkce2mcwndc}
Let now focus our attention on expansions \eqref{h37qp2b237ds}. 
An appropriate use of these series provides numerous fascinating formul\ae,
and especially, the series expansions of the MacLaurin--Taylor type for the composite functions involving logarithms and inverse 
trigonometric and hyperbolic functions. The technique is based
of the summation over $l$ of \eqref{h37qp2b237ds}, on the fact that $S_1(n,l)$ vanishes for $l\notin[1,n]$
and on the interchanging the order of summation\footnote{Series in question being absolutely convergent.}. 
For example, the trivial 
summation of the right part of (\ref{h37qp2b237ds}b) over $l\in[1,\infty)$, yields
\be\notag
\sum_{l=1}^\infty\frac{\ln^l(1+z)}{l!}=e^{\ln(z+1)}-1=z
\ee
since the sum in the left--hand side is simply the MacLaurin series of $\,e^{\ln(z+1)}$ without the first term. 
At the same time, the summation of the left part of (\ref{h37qp2b237ds}b) results in
\be\notag
\sum_{l=1}^\infty \sum_{n=l}^\infty\!\frac{S_1(n,l)}{n!}z^n = 
\sum_{l=1}^\infty \sum_{n=1}^\infty\!\frac{S_1(n,l)}{n!}z^n = 
\sum_{n=1}^\infty\!\frac{z^n}{n!} \underbrace{\sum_{l=1}^n S_1(n,l)}_{\delta_{n,1}} = z
\ee
where the last sum may be truncated at $l=n$ thanks to \eqref{io20323m3e} and equals 1 if $n=1$ and 0 otherwise. 
Let now consider more complicated cases.
Write in (\ref{h37qp2b237ds}b) $2l$ for $l$, and then, sum the result with respect to $l$ from $l=1$ to $l=\infty$. This gives
\be
\ch\ln(1+z)\,=\,1+\sum_{n=1}^\infty\!\frac{z^n}{n!} \cdot\!\!\sum_{l=1}^{\lfloor\frac{1}{2}n\rfloor}\!S_1(n,2l)=
1+\frac{1}{2}\sum_{n=2}^\infty (-1)^n z^n
\,,\qquad |z|<1
\ee
By the same line of reasoning, if we divide the right--hand side of (\ref{h37qp2b237ds}b) by $l+1$ and sum it over $l\in[1,\infty)$, then we get
\be\notag
\begin{array}{ll}
\displaystyle
\sum_{l=1}^\infty\frac{1}{l+1}\cdot\frac{\ln^l(1+z)}{l!} \,=\, \frac{1}{\ln(1+z)}\sum_{l=1}^\infty\! \frac{\ln^{l+1}(1+z)}{(l+1)!}\\[5mm]
\displaystyle\qquad
\,=\, \frac{1}{\ln(1+z)} \left[e^{\ln(1+z)}-\ln(1+z)-1\right] \,=\,\frac{z}{\ln(1+z)} -1
\end{array}
\ee
Applying the same operation to the left--hand side of (\ref{h37qp2b237ds}b) and comparing both sides yields
\be\label{k2390234mrf}
\frac{z}{\ln(1+z)}\,=\,1+\sum_{n=1}^\infty z^n\!\cdot\!
\underbrace{
\frac{1}{n!} \!\sum_{l=1}^{n} \frac{S_1(n,l)}{l+1} 
}_{G_n}\,
=\,1+\sum_{n=1}^\infty G_n\, z^n  \,,\quad
\qquad |z|<1\,,
\ee
the equality, which is more known as the generating equation for the  \emph{Gregory's coefficients} $G_n$ (in particular, $G_1=+\frac{1}{2}$, 
\mbox{$G_2=-\frac{1}{12}$}, $G_3=+\frac{1}{24}$, $G_4=-\frac{19}{720}$, 
$G_5=+\frac{3}{160}$, $G_6=-\frac{863}{60\,480}$, \ldots).\footnote{Coefficients 
$\,G_n=\frac{1}{n!} \sum\limits_{l=1}^{n} \frac{S_1(n,l)}{l+1} =\frac{1}{n!}\int\limits_0^1 (x-n+1)_n\, dx=-\frac{B_n^{(n-1)}}{(n-1)\,n!}
=\frac{C_{1,n}}{n!}\,$ are also called
\emph{(reciprocal) logarithmic numbers}, \emph{Bernoulli numbers of the second kind},
normalized \emph{generalized Bernoulli numbers} $B_n^{(n-1)}$, \emph{Cauchy numbers} and normalized \emph{Cauchy numbers
of the first kind} $C_{1,n}$. They were introduced by
James Gregory in 1670 in the context of area's interpolation formula (which is known nowdays as \emph{Gregory's interpolation formula})
and were subsequently rediscovered in various contexts by many famous mathematicians, including Gregorio Fontana, Lorenzo Mascheroni ,
Pierre--Simon Laplace, Augustin--Louis Cauchy, Jacques Binet, 
Ernst Schr\"oder, Oskar Schl\"omilch, Charles Hermite, Jan C.~Kluyver and Joseph Ser 
\cite[vol.~II, pp.~208--209]{rigaud_01},
\cite[vol.~1, p.~46, letter written on November 23, 1670 to John Collins]{newton_01}, \cite[pp.~266--267, 284]{jeffreys_02}, 
\cite[pp.~75--78]{goldstine_01},
\cite[pp.~395--396]{chabert_01}, 
\cite[pp.~21--23]{mascheroni_01}, \cite[T.~IV, pp.~205--207]{laplace_01}, \cite[pp.~53--55]{boole_01}, \cite{van_veen_01}, 
\cite[pp.~192--194]{goldstine_01},
\cite{lienard_01}, \cite{wachs_01}, \cite{schroder_01}, \cite{schlomilch_03}, \cite[pp.~65, 69]{hermite_01},
\cite{kluyver_02}, \cite{ser_01}.
For more information about these important coefficients, see
\cite[pp.~240--251]{norlund_02}, \cite{norlund_01}, \cite[p.~132, Eq.~(6), p.~138]{jordan_02}, \cite[p.~258, Eq.~(14)]{jordan_00}, \cite[pp.~266--267, 277--280]{jordan_01}, 
\cite{nielsen_01}, \cite{nielsen_03}, \cite{steffensen_00}, \cite{steffensen_01},  \cite[pp.~106--107]{steffensen_02}, \cite{davis_02},
\cite[p.~190]{weisstein_04},
\cite[p.~45, \no 370]{gunter_03_eng}, 
\cite[vol.~III, pp.~257--259]{bateman_01}, \cite{stamper_01}, \cite[p.~229]{krylov_01}, \cite[\no 600, p.~87]{proskuriyakov_01_eng}, 
\cite[p.~216, \no 75-a]{knopp_01}
\cite[pp.~293--294, \no13]{comtet_01}, \cite{carlitz_01}, \cite{howard_02}, \cite{young_01}, \cite{adelberg_01}, \cite{zhao_01}, 
\cite{candelpergher_01}, 
\cite[Eq.~(3)]{mezo_01}, \cite{merlini_01}, \cite{nemes_01},
\cite[pp.~128--129]{alabdulmohsin_01}, 
\cite[Chapt.~4]{arakawa_01}, \cite{skramer_01}.\label{jpqwcnqasd}}
Analogously, performance of same procedures with (\ref{h37qp2b237ds}a), written for $-z$ instead of $z$, results in 
\be\label{i2039dj23r}
\frac{z}{(1+z)\ln(1+z)}\,=\,1+\sum_{n=1}^\infty\!\frac{(-1)^n z^n}{n!} \cdot\!\underbrace{\sum_{l=1}^{n} \frac{|S_1(n,l)|}{l+1}}_{C_{2,n}}\,=
\,1+\sum_{n=1}^\infty\!\frac{(-1)^n C_{2,n} }{n!} \, z^n\,,
\quad\qquad |z|<1\,,
\ee
which is also known as the generating series for the \emph{Cauchy numbers of the second kind} $C_{2,n}$ (in particular,
$C_{2,1}=\frac{1}{2}$, \label{lkjfc02kme}
\mbox{$C_{2,2}=\frac{5}{6}$}, $C_{2,3}=\frac{9}{4}$, $C_{2,4}=\frac{251}{30}$, 
$C_{2,5}=\frac{475}{12}$, $C_{2,6}=\frac{19\,087}{84}$, \ldots).\footnote{These
numbers $\,C_{2,n}=\sum\limits_{l=1}^{n} \frac{|S_1(n,l)|}{l+1} = \int\limits_0^1 (x)_n\, dx=|B_n^{(n)}|\,$, 
called by some authors signless \emph{generalized Bernoulli numbers}
$|B_n^{(n)}|$ and signless \emph{Nørlund numbers}, are much less famous than
Gregory's coefficients $G_n$, but their study is also very interesting, see e.g.~\cite[pp.~150--151]{norlund_02}, 
\cite{norlund_01}, \cite[vol.~III, pp.~257--259]{bateman_01}, 
\cite[pp.~293--294, \no13]{comtet_01}, \cite{howard_03}, \cite{adelberg_01}, \cite{zhao_01}, \cite{qi_01}.\label{c23409jc}}
Dividing by $z$, integrating and determining the constant of integration yields another interesting series
\be\label{ki1d039dn321}
\ln\ln(1+z)\,=\,\ln z +\sum_{n=1}^\infty\!\frac{(-1)^n z^n}{n\cdot n!} \cdot\sum_{l=1}^{n} \frac{|S_1(n,l)|}{l+1}
\,=\,\ln z +\sum_{n=1}^\infty\!\frac{(-1)^n C_{2,n}}{n\cdot n!} \, z^n\,,\quad
\qquad |z|<1\,,
\ee
which is an ``almost MacLaurin series'' for $\ln\ln(1+z)$.
Asymptotic studies of general terms in series \eqref{k2390234mrf} and \eqref{ki1d039dn321} reveal that for $n\to\infty$
both terms decrease logarithmically:
\be\label{j38ndbr893r}
G_n\,=\,\frac{1}{\,n!\,}\!\sum_{l = 1 }^{n} \!\frac{\,S_1(n, l)\,}{l+1}\sim\frac{(-1)^{n-1}}{\,n\ln^2\! n\,}\qquad
\text{and}\qquad
\frac{C_{2,n}}{\,n\cdot n!\,}\,=\,
\frac{1}{n\cdot n!}\sum_{l=1}^{n} \frac{|S_1(n,l)|}{l+1} \sim \frac{1}{\,n\ln n\,}
\ee
respectively (see \ref{appendix2} and \ref{appendix1}),
and hence, series \eqref{k2390234mrf} and \eqref{ki1d039dn321} converge not only in $\,|z|<1\,$, but also at $z=\pm1$ 
at $z=1$ respectively.
Thus, putting $\,z=1\,$ into \eqref{k2390234mrf}, we have
\be\label{kjc0293nfr}
\frac{1}{\,\ln2\,}\,=\,1+\sum_{n=1}^\infty\!\frac{1}{\,n!\,}\sum_{l=1}^{n} \frac{S_1(n,l)}{\,l+1\,}\,=\,1+\sum_{n=1}^\infty G_n
\ee
while setting $z=1$ into \eqref{ki1d039dn321} gives a series for $\,\ln\ln2\,$, see \eqref{fmse}.
Moreover, application of Abel's theorem on power's series to \eqref{k2390234mrf} at $z\to-1^+$ yields
\emph{Fontana's series}\footnote{This series appears in 
a letter of Gregorio Fontana to which Lorenzo Mascheroni refers in \cite[pp.~21--23]{mascheroni_01}.}
\be\label{uih39hlj2983dh}
1\,=\,\sum_{n=1}^\infty\!\frac{(-1)^{n-1}}{n!}\sum_{l=1}^{n} \frac{S_1(n,l)}{l+1}\,=
\sum_{n=1}^\infty\!\big|G_n\big|
\ee
converging at the same rate as $\,\sum n^{-1}\ln^{-2}\!n$, see \eqref{j38ndbr893r}.

The use of the same and of similar 
techniques allows to readily derive expressions for the coefficients of the MacLaurin series for even more complicated functions. 
Further examples demonstrate better than words the powerfulness of the method (note
that some examples are actually the Laurent series in a neighborhood of $z=0$):
\begin{eqnarray}
&&\displaystyle\label{ld2k3memne}
\sh\ln(1+z)\,=\,\sum_{n=1}^\infty\!\frac{z^n}{n!} \cdot\!\!\sum_{l=0}^{\lfloor\frac{1}{2}n\rfloor}\!S_1(n,2l+1) = 
z-\frac{1}{2}\sum_{n=2}^\infty\! (-1)^n z^n \,,\qquad\quad   |z|<1\quad\\[4mm]
&&\displaystyle
\cos\ln(1+z)\,=\,1+\sum_{n=1}^\infty\!\frac{z^n}{n!}\cdot\!\!\sum_{l=1}^{\lfloor\frac{1}{2}n\rfloor}\! (-1)^lS_1(n,2l)\,,
\qquad\quad   |z|<1  \\[4mm]
&&\displaystyle
\sin\ln(1+z)\,=\,\sum_{n=1}^\infty\!\frac{z^n}{n!} \cdot\!\!\sum_{l=0}^{\lfloor\frac{1}{2}n\rfloor} \!(-1)^l S_1(n,2l+1)\,,
\qquad\quad   |z|<1 \\[4mm]
&&\displaystyle
\ln\!\big[1+\ln(1+z)\big]\,=\,\sum_{n=1}^\infty\!\frac{z^n}{n!} \cdot\sum_{l=0}^{n-1} (-1)^l l!\cdot S_1(n,l+1)
\,,\qquad\quad  |z|<1-e^{-1}\approx0.63\qquad\\[3mm]
&&\displaystyle\label{k2390234mrf2}
\frac{1}{\,\ln^2(1+z)\,}\,=\,\frac{1}{z^2}+\frac{1}{z}+\sum_{n=0}^\infty\!\frac{\,z^{n}}{(n+2)!} \cdot\sum_{l=1}^{n+1} 
\frac{\,1-n(l+1)\,}{(l+1)(l+2)}\cdot S_1(n+1,l)\,, \qquad\quad |z|<1 \\[4mm]
&&\displaystyle\label{k02eck432jr}
\frac{1}{\,\ln^m(1+z)\,}\,=\,\frac{1}{z}\cdot\!\sum_{k=1}^{m-1}\!\frac{\,1}{\,k!\cdot\ln^{m-k}(1+z)\,} +\frac{1}{\,m!\cdot z\,}
+\sum_{n=1}^\infty\!\frac{\,z^{n-1}}{n!} \cdot\sum_{l=1}^{n} 
\frac{\,S_1(n,l)\,}{(l+1)_m}\,,
\qquad 
\begin{array}{l}
m=2, 3, 4,\ldots \\[1mm]
|z|<1
\end{array}
\quad\\[4mm]
&&\displaystyle\label{dbe46gb}
\frac{\ln^m(1+z)}{1+z}\,=\,(-1)^m m!\cdot\!
\sum_{n=0}^\infty \frac{\,(-1)^n\big|S_1(n+1, m+1)\big|\,}{n!} \,z^n\,,
\qquad 
\begin{array}{l}
m=0, 1, 2,\ldots  \\[1mm]
|z|<1
\end{array}
\quad\\[4mm]
&&\displaystyle\label{jh9sqhhl}
\arctg\ln(1+z) \,=\sum_{n=1}^\infty\!\frac{z^{n}}{\,n!\,}\cdot\!\sum_{l=0}^{\lfloor\frac{1}{2}n\rfloor}\! (-1)^{l} 
(2l)!\cdot S_1(n,2l+1)\,, \qquad\quad |z|<2\sin\frac{1}{2}\approx0.96\\[4mm]
&&\displaystyle\label{jh9sqhhl}
\arcth\ln(1+z) \,=\sum_{n=1}^\infty\!\frac{z^{n}}{\,n!\,}\cdot\!\sum_{l=0}^{\lfloor\frac{1}{2}n\rfloor} \!
(2l)!\cdot S_1(n,2l+1)\,, \qquad\quad  |z|<1-e^{-1}\approx0.63\qquad\\[4mm]
&&\displaystyle\label{jh9sqhhl2}
\frac{\arcth^m z}{m!} \,=\sum_{n=m}^\infty\!z^{n}\cdot\sum_{l=m}^{n} \!\binom{n-1}{l-1}
\cdot\frac{\,2^{l-m}\cdot S_1(n,l)\,}{l!} \,,\qquad m=1, 2, 3,\ldots\,, \qquad\quad |z|<1 \\[4mm]
&&\displaystyle\label{jx20j2nm3d}
\tg\ln(1+z)\,=\sum_{n=1}^\infty\!\frac{z^n}{n!}\,\cdot\!\sum_{l=0}^{\lfloor\frac{1}{2}n\rfloor} \!
\frac{\,2^{2l+1}(2^{2l+2}-1)\cdot|{B}_{2l+2}|\cdot S_1(n,2l+1)\,}{l+1} 
\,, \qquad\quad  |z|<1-e^{-\pi/2}\approx0.79\qquad\quad   \\[4mm]
&&\displaystyle\notag
\th\ln(1+z)\,=\sum_{n=1}^\infty\!\frac{z^n}{n!}\,\cdot\!\sum_{l=0}^{\lfloor\frac{1}{2}n\rfloor} \!(-1)^l
\frac{\,2^{2l+1}(2^{2l+2}-1)\cdot|{B}_{2l+2}|\cdot S_1(n,2l+1)\,}{l+1}\\[1mm]
&&\displaystyle\qquad\qquad\quad\;\,\label{jx20jve548y}
=\sum_{n=0}^\infty(-1)^{n}\frac{z^{4n+1}}{2^{2n}}
-\sum_{n=0}^\infty(-1)^{n}\frac{z^{4n+2}}{2^{2n+1}}
+\sum_{n=0}^\infty(-1)^{n}\frac{z^{4n+4}}{2^{2n+2}}
\,, \qquad\quad  |z|<\sqrt2\approx1.41\qquad   
\end{eqnarray}
Derived expansions coincide with the corresponding MacLaurin series, converge everywhere where expanded functions 
are analytic\footnote{The radius of convergence of such series $r$ is conditioned by the singularities of expanded functions. For instance,
$\ln\!\big[1+\ln(1+z)\big]$ is analytic on the entire complex $z$--plane except points at which $1+z=0$ and $1+\ln(1+z)=0$, which are both branch points.
From the former we conclude that the radius of convergence cannot be greater than 1, and from the latter, it follows that it cannot
exceed  $1-e^{-1}$ which is even lesser than 1. Hence $r=1-e^{-1}\approx0.63$} and contain rational coefficients only.
The main advantage of this technique is that we do not need to ``mechanically''
compute the $n$th derivative of the composite function at $z=0$, which often may be a very laborious task.\footnote{Some other 
power series expansions 
involving Stirling numbers are also given in works of Wilf \cite{wilf_01}, Kruchinin \cite{kruchinin_01}, \cite{kruchinin_02}, \cite{kruchinin_03}
and Rz\c{a}dkowski \cite{rzadkowski_01}. Moreover, series expansions of certain composite functions, not necessarily containing Stirling numbers,
may be found in \cite[Chapt.~VI]{knopp_01}, \cite[p.~20 \& 63]{gunter_03_eng}, \cite{hansen_01} and \cite[vol.~I]{prudnikov_en}
(in the third reference, the author also provides a list of related references)}\up{,}\footnote{Since these expansions 
are not particularly difficult to obtain and also may be derived by other techniques, it is possible that some of them 
could appear in earlier works. The same remark also concerns
formul\ae~\eqref{ioj2c404jne}--\eqref{kljc204ijnd}. For instance, formula \eqref{uiewhc83he} may be found in other sources as well, 
see e.g.~\cite[p.~431, Eq.~(76)]{kowalenko_01}, \cite[p.~128, Eq.~(7.3.11)]{alabdulmohsin_01}, \cite[p.~4006]{young_03}
(the same series also appears in \cite[p.~14, Eq.~(2.39)]{coffey_07}, but the result is incorrect). Series \eqref{hgwibddw} is also known, 
see e.g.~\cite[p.~2952, Eq.~(1.3)]{young_01},  \cite[p.~20, Eq.~(3.6)]{coffey_07}, \cite[p.~307, Eq.~for $F_0(2)$]{candelpergher_01}.}

Generating equations for Stirling numbers of the first kind may be also successfully used for the derivation of more complicated and quite
unexpected results. 
For instance, it is known that 
\be\label{kljd023jdnr}
\zeta(k+1)\,=\,\sum_{n=k}^\infty\frac{|S_1(n,k)|}{n\cdot n!}\,,\qquad\quad k=1, 2, 3,\ldots
\ee
see e.g.~Jordan's book \cite[pp.~166, 194--195]{jordan_01}.\footnote{Jordan derives this formula and remarks that particular cases of it
were certainly known to Stirling \cite{stirling_01} (see also \cite[p.~302, Eq.~(36bis)]{nielsen_04},
\cite[p.~10]{tweedie_01}, and compare it to formul\ae~from \cite[p.~11]{stirling_01}).} 
This result was recently rediscovered by several modern writers, e.g.~by Shen \cite{shen_01} and Sato \cite{sato_01}; however, their proofs are exceedingly long.
Using \eqref{h37qp2b237ds} the whole procedure takes only a few lines:
\be\label{kljd023jdnr2}
\begin{array}{ll}
\displaystyle
\sum_{n=k}^\infty\frac{|S_1(n,k)|}{n\cdot n!} \,=\,
\sum_{n=k}^\infty\frac{|S_1(n,k)|}{ n!}\!\underbrace{\int\limits_0^1 \! x^{n-1} dx}_{1/n}\,=\int\limits_0^1 \!
\underbrace{\sum_{n=k}^\infty\frac{|S_1(n,k)|}{ n!}\,  x^n }_{\text{see \eqref{h37qp2b237ds}}}\frac{dx}{x}\,=\,
\frac{(-1)^k}{k!}\int\limits_0^1 \!\frac{\ln^{k}(1-x)}{x}\,dx\\[8mm]
\displaystyle\qquad\qquad\qquad\;
\,=\,\frac{(-1)^k}{k!}\int\limits_0^\infty\! \frac{(-t)^{k}}{\,e^t\big(1-e^{-t}\big)\,}\,dt\,=\,
\frac{1}{\Gamma(k+1)}\int\limits_0^\infty\! \frac{t^{k}}{\,e^t-1\,}\,dt\,=\,\zeta(k+1)
\end{array}
\ee
where in last integrals we made a change of variable $\,x=1-e^{-t}$. 
The above formula may be readily generalized to
\be\label{kljd023jdnr}
\zeta(k+1,v)\,=\,\sum_{n=k}^\infty\frac{|S_1(n,k)|}{n\cdot (v)_n}\,,\qquad\quad k=1, 2, 3,\ldots\,,
\qquad \Re v>0\,.
\ee
where at large $n$
\be
\frac{|S_1(n,k)|}{n\cdot (v)_n} \,\sim\, \frac{\Gamma(v)}{\,(k-1)!\,}\cdot\frac{\ln^{k-1}\! n}{n^{v+1}}
\,,\qquad\qquad n\to\infty\,,
\ee
in virtue of \eqref{lk2093mffmnjw} and a known asymptotics for the Stirling numbers 
\cite[p.~261]{jordan_00}, \cite[p.~161]{jordan_01}, \cite[\no24.1.3, p.~824]{abramowitz_01}, \cite[p.~348, Eq.~(8)]{wilf_02}.
Moreover, by a slight modification of the above technique,
we may also obtain the following results:
\be\label{ioj2c404jne}
\begin{array}{ll}
\displaystyle
\sum_{n=1}^{\infty}\!\frac{\, (-1)^{n-1}\,}{\,n\,}\cdot\frac{1}{\,n!\,}\sum_{l = 1 }^{n}  f(l)\,S_1(n, l)
\,=\,\sum_{l=1}^{\infty} (-1)^{l+1}f(l)\zeta(l+1)\,, \\[6mm]
\displaystyle
\sum_{n=1}^{\infty}\!\frac{\, (-1)^{n-1}H_n\,}{\,n\,}\cdot\frac{1}{\,n!\,}\sum_{l = 1 }^{n}  f(l)\,S_1(n, l)
\,=\,\sum_{l=1}^{\infty} (-1)^{l+1}(l+1)f(l)\zeta(l+2)\,,
\end{array}
\ee
where $f(l)$ is an arbitrary function ensuring the convergence and $H_n$ is the $n$th harmonic number;
\be
\begin{array}{ll}
\displaystyle
\sum_{n=1}^{\infty}\!\frac{\, (-1)^{n-1}\,}{\,n\,}\cdot\,\frac{1}{\,n!\,}\sum_{l = 1 }^{n} \!\frac{\,S_1(n, l)\,}{l+k}
\,=\sum_{l=2}^{\infty}\!\frac{(-1)^{l}\!\cdot\zeta(l)}{\,l+k-1\,}\,=\,\frac{\,1\,}{\,k-1\,} -\frac{\,\ln2\pi\,}{k}+\frac{\,\gamma\,}{2} 
\\[6mm]
\displaystyle\qquad
+\!\!\!\!\sum_{l=1}^{\lfloor\frac{1}{2}(k-1)\rfloor}\!\!\!\! (-1)^l \binom{k-1}{2l-1}\frac{\,(2l)!\cdot\zeta'(2l)\,}{l\cdot(2\pi)^{2l}}
+\!\!\sum_{l=1}^{\lfloor\frac{1}{2}k\rfloor-1}\!\! (-1)^l \binom{k-1}{2l}\frac{\,(2l)!\cdot\zeta(2l+1)\,}{2\cdot(2\pi)^{2l}}
\qquad
\end{array}
\ee
where $k=2, 3, 4,\,\ldots$\, and where the series on the left converges as $\,\sum n^{-2} \ln^{-k-1}\! n \,$;
for $k=1, 2, 3,\ldots$
\begin{eqnarray}
&&\displaystyle\sum_{n=1}^{\infty}\!\frac{\, (-1)^{n-1}\,}{\,n+k\,}\cdot\frac{1}{\,n!\,}\sum_{l = 1 }^{n} \frac{\,S_1(n, l)\,}{l+1}\,
=\sum_{n=1}^{\infty}\!\frac{\big|\,G_n\big|}{\,n+k\,}\,
=\,\frac{1}{k}+\sum_{m=1}^{k} (-1)^m \binom{k}{m}\ln(m+1)\notag\\[1mm]
&&\displaystyle\qquad\qquad\qquad\qquad\qquad\qquad\qquad\qquad\qquad\quad\;\,
\,=\,\frac{1}{k} + \Big.\Delta^k\ln(x)\Big|_{x=1}\,
\label{s2kj0sdw} \\[1mm]
&&\displaystyle\sum_{n=1}^{\infty}\!\frac{\, (-1)^{n-1}\,}{\,n\,}\cdot\frac{1}{\,n!\,}\sum_{l = 1 }^{n} \frac{\,S_1(n, l)\,}{l+1}
\,=\sum_{n=1}^{\infty} \frac{\big|\,G_n\big|}{n}\,=\,\gamma \label{fms2x}\\[1mm]
&&\displaystyle\sum_{n=2}^{\infty}\!\frac{\, (-1)^{n-1}\,}{\,n-1\,}\cdot\frac{1}{\,n!\,}\sum_{l = 1 }^{n} \frac{\,S_1(n, l)\,}{l+1}
\,=\sum_{n=2}^{\infty}\!\frac{\big|\,G_n\big|}{\,n-1\,} 
\,=\,  -\frac{1}{2} + \frac{\ln2\pi}{2} -\frac{\gamma}{2} \label{uiewhc83he}\ \\[1mm]
&&\displaystyle\sum_{n=3}^{\infty}\!\frac{\, (-1)^{n-1}\,}{\,n-2\,}\cdot\frac{1}{\,n!\,}\sum_{l = 1 }^{n} \frac{\,S_1(n, l)\,}{l+1}
\,=\sum_{n=3}^{\infty}\!\frac{\big|\,G_n\big|}{\,n-2\,}
\,=\,  -\frac{1}{8} + \frac{\ln2\pi}{12} - \frac{\zeta'(2)}{\,2\pi^2}\\[1mm]
&&\displaystyle\sum_{n=4}^{\infty}\!\frac{\, (-1)^{n-1}\,}{\,n-3\,}\cdot\frac{1}{\,n!\,}\sum_{l = 1 }^{n} \frac{\,S_1(n, l)\,}{l+1}
\,=\sum_{n=4}^{\infty}\!\frac{\big|\,G_n\big|}{\,n-3\,}
\,=\,  -\frac{1}{16} + \frac{\ln2\pi}{24} - \frac{\zeta'(2)}{\,4\pi^2} + \frac{\zeta(3)}{\,8\pi^2}  
\end{eqnarray}
where $\Delta^k$ is the $k$th finite difference, see e.g.~\cite[p.~270, Eq.~(14.17)]{vorobiev_01},
and where all series on the left converges as $\,\sum (n \ln n )^{-2}\,$;
\begin{eqnarray}
&&\displaystyle\label{hgwibddw}
\sum_{n=1}^{\infty}\!\frac{\, (-1)^{n-1}H_n\,}{\,n\,}
\cdot\frac{1}{\,n!\,}\sum_{l = 1 }^{n} \frac{\,S_1(n, l)\,}{l+1}\,=
\sum_{n=1}^{\infty}\!\frac{\, \big|\,G_n\big|\cdot H_n\,}{\,n\,}\!\,=
\,\frac{\pi^2}{6}-1 \\[1mm]
&&\displaystyle
\sum_{n=1}^{\infty}\!\frac{\, (-1)^{n-1}H_n\,}{\,n\,}\cdot\frac{1}{\,n!\,}\sum_{l = 1 }^{n}  
 \frac{\,S_1(n, l)\,}{\,(l+1)(l+2)\,} \,=\,\frac{\pi^2}{12}- \gamma \\[1mm]
&&\displaystyle\label{kljc204ijnd}
\sum_{n=1}^{\infty}\!\frac{\, (-1)^{n-1}H_n\,}{\,n\,}\cdot\frac{1}{\,n!\,}\sum_{l = 1 }^{n}  
 \frac{\,S_1(n, l)\,}{\,(l+1)(l+3)\,}\,=\,\frac{\pi^2}{18} +\frac{1}{2}\ln2\pi -\frac{\gamma}{2}- 1\!
\end{eqnarray}
which all converge as $\sum n^{-2} \ln^{-1}\! n$, and even this beautiful alternating series
\be\label{c234c2gf5g}
\sum_{n=1}^{\infty}\!\frac{\, 1\,}{n}\cdot
\frac{1}{\,n!\,}\!\sum_{l = 1 }^{n} \!\frac{\,S_1(n, l)\,}{l+1}\,=
\sum\limits_{n=1}^\infty \! \frac{\,G_n\,}{n}\,=\,\Ei(\ln2)-\gamma\,=\,\li(2)-\gamma\,,
\ee
where $\Ei(\cdot)$ and $\li(\cdot)$ denote exponential integral and logarithmic integral functions respectively.\footnote{Numbers $G_n$ 
are strictly alternating: $G_n=(-1)^{n-1} \big|G_n\big|\,$. The left side of \eqref{c234c2gf5g}
is, therefore, the alternating variant of Fontana--Mascheroni's series \eqref{fms}, \eqref{fms2x}, and from various points of view
the constant $\,\li(2)-\gamma\,=\,0.4679481152\ldots$ may be regarded 
as the \emph{alternating Euler's constant}, by analogy to $\,\ln\frac{4}{\pi}\,$,
which was earlier proposed as such by Jonathan Sondow in \cite{sondow_02}.}
Finally, Stirling numbers of the first kind may also appear in the evaluation of certain integrals, which, at first sight, have nothing to do 
with Stirling numbers. For instance, if $k$ is positive integer and $\Re{s}>k-1$, then
\be\label{kljc204ijnd}
\begin{array}{ll}
\displaystyle
\int\limits_0^1 \! \frac{\,\ln^s(1-x)\,}{x^k}\, dx \,=\, \frac{\,(-1)^s\!\cdot \Gamma(s+1)\,}{(k-1)!}\cdot
\begin{cases}
\,\zeta(s+1)  \,,\quad & k=1 \\[2mm]
\,\zeta(s)   \,,\quad  &  k=2\\[1.6mm]
\,\zeta(s-1)+\zeta(s) \,,\quad  & k=3\\[1.6mm]
\,\zeta(s-2)+3\zeta(s-1)+2\zeta(s) \,,\quad  & k=4 \\[1.6mm]
\displaystyle\sum_{r=1}^{k-1} S_1(k-1,r) \! \sum_{m=0}^r \!
\binom{r}{m}\!\cdot(k-2)^{r-m}\cdot\zeta(s+1-m)\,,\quad  & k\geqslant3
\end{cases}
\end{array}
\ee
The proofs of some of these results being quite long, we placed them in \ref{app_1}.

\subsection{An inspiring example for the derivation of the series for $\ln\Gamma(z)$}\label{j8y2bs4r}
\begin{figure}[!t]   
\centering
\includegraphics[width=0.8\textwidth]{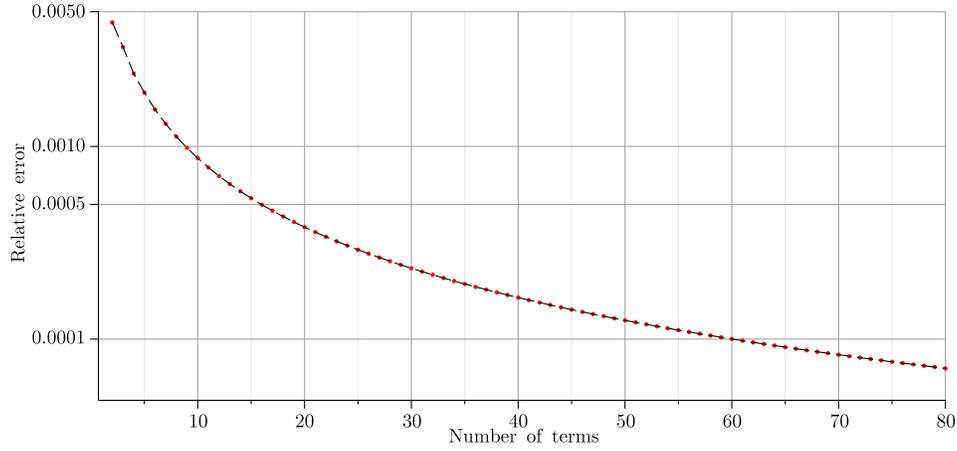}
\caption{Relative error of series expansion \protect\eqref{kx1230xk3ec}, logarithmic scale.}
\label{jhc3uiexn2}
\end{figure}

Let now consider the example
which was originally our inspiration for this work. 
In exercise \no 39-b in \cite[Sect.~4]{iaroslav_06} we established that
\be\label{oijh23098uyrb}
\int\limits_0^{\,1}\! \frac{\,\arctg\arcth{x}\,}{x}\, dx
=\,\pi\left\{\ln\Gamma\!\left(\!\dfrac{1}{\,\pi\,}\!\right)
- \ln\Gamma\!\left(\!\dfrac{1}{2}+\dfrac{1}{\,\pi\,}\!\right) 
-\frac{1}{2}\ln\pi\right\} =\, 1.025760510\ldots  
\ee
The arctangent of the hyperbolic arctangent is analytic in the whole disk $|x|<1$, 
and therefore, can be expanded into the MacLaurin series.\footnote{Function $\arctg\arcth x$ 
has branch points at $x=\pm 1$ and $x=\pm i\tg1\approx\pm1.56i$.} 
The coefficients of such an expansion require a careful watching, the law for their formation
being difficult to derive by inductive or semi--inductive methods.
So we resort again to the method employing Stirling numbers:
\be\notag
\begin{array}{l}
\displaystyle 
\arctg\arcth{x}\,=\sum_{l=0}^\infty (-1)^l (2l)!\cdot\frac{\arcth^{2l+1}\! x}{(2l+1)!} =
\sum_{n=1}^\infty x^n \!\cdot\!\sum_{k=1}^n\!\binom{n-1}{k-1}\frac{2^k}{k!}\cdot\!\! \sum_{l=0}^{\lfloor\frac{1}{2}n\rfloor}\!
(-1)^l \cdot\frac{(2l)!\cdot S_1(k,2l+1)}{2^{2l+1}}  \\[10mm]
\displaystyle \qquad\qquad\qquad
=\sum_{n=0}^\infty x^{2n+1} \!\cdot\!\underbrace{\sum_{k=1}^{2n+1}\!\binom{2n}{k-1}\frac{2^k}{k!}\cdot\! \sum_{l=0}^{n}
(-1)^l \cdot\frac{(2l)!\cdot S_1(k,2l+1)}{2^{2l+1}}}_{A_n} =\,
x+\frac{1}{15}x^5+\frac{\,1}{45}x^7\\[10mm]
\displaystyle \qquad\qquad\qquad\phantom{mm}
+\frac{\,64}{2835}x^9+\frac{\,71}{4725}x^{11}
+\frac{\,5209}{405\,405}x^{13} + \frac{2\,203\,328}{212\,837\,625} x^{15}+\ldots\,, \quad \qquad |x|<1\,,
\end{array}
\ee
where we used result \eqref{jh9sqhhl2}, as well as the oddness of the expanded function.\footnote{Note that although
the MacLaurin series for the arctangent is valid only in the unit circle, i.e.~formally only for such $x$ that \mbox{$|\arcth{x}|<1$},
the above expansion holds uniformly in the whole disk $|x|<1$ (in virtue of the Cauchy's theorem on the representation 
of analytic functions by power series, as well as of the principle of analytic continuation). Moreover, an advanced study of this series,
analogous to that performed in the next section, shows that it also converges for $x=1$. 
}
Inserting this expansion into \eqref{oijh23098uyrb} and performing the term--by--term integration,
we obtain
the following series for the difference of first two terms in curly brackets in \eqref{oijh23098uyrb}
\be\label{kx1230xk3ec}
\begin{array}{ll}
\displaystyle
\ln\Gamma\!\left(\!\dfrac{1}{\,\pi\,}\!\right)
&\displaystyle - \ln\Gamma\!\left(\!\dfrac{1}{2}+\dfrac{1}{\,\pi\,}\!\right) = \frac{1}{2}\ln\pi 
+ \frac{1}{\pi}\!\sum_{n=0}^\infty\frac{A_n}{2n+1} 
\,=\,\frac{1}{2}\ln\pi 
+ \frac{1}{\pi}\left\{1+ \frac{1}{75}+ \frac{1}{315}+ \frac{64}{25\,515} \right.
\\[8mm]
&\displaystyle \!\!\!\!\!\!\left.
+ \,\frac{71}{51\,975}+ \frac{5209}{5\,270\,265}+ \frac{2\,203\,328}{3\,192\,564\,375}+ \frac{132\,313}{253\,127\,875}
+\ldots \right\}=\,0.8988746544\ldots
\end{array}
\ee
with $A_n$ defined in the preceding equation. 
The derived series does not converge rapidly, see Fig.~\ref{jhc3uiexn2}, but the most remarkable 
is  that it contains rational coefficients only, which is quite unusual, especially
for the arguments related to $\pi^{-1}$.
This suggests that there might be some more general series similar in nature to \eqref{kx1230xk3ec}, 
which allows to expand the logarithm
of the $\Gamma$--function at certain points related to $\pi^{-1}$ into the series with rational coefficients only.
Such series expansions are the subject of our study in the next section.

\section{Series expansions for the logarithm of the $\Gamma$--function and polygamma functions}\label{oijxc32094djfd}
\subsection{First series expansion for the logarithm of the $\Gamma$--function} 
\subsubsection{Derivation of the series expansion}
\begin{figure}[!t]   
\centering
\includegraphics[width=0.8\textwidth]{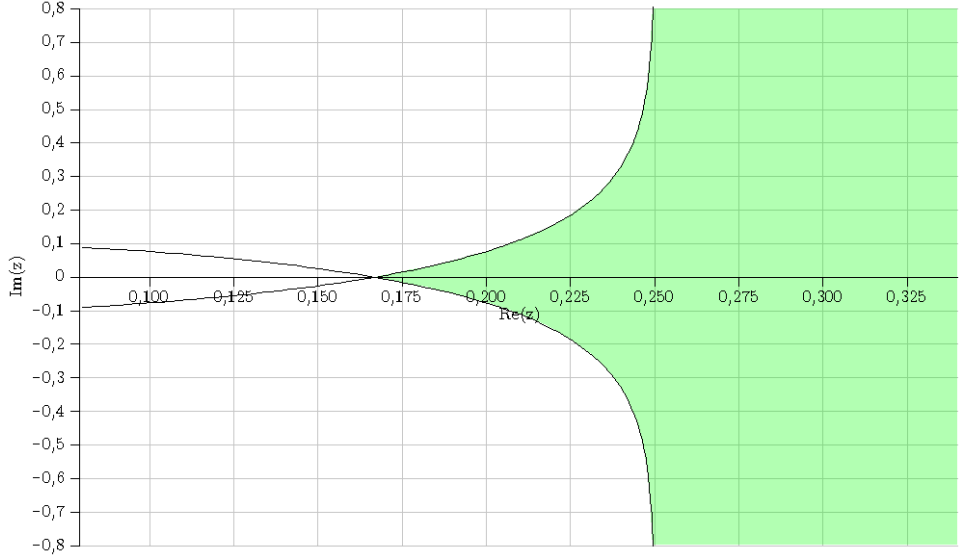}
\caption{The region of convergence of series \eqref{k034fjdfnsh2} and \eqref{kj20ejcn2dnd} in the complex $z$--plane for $|u|<1$
is the common part of two zones, each of which lying to the right of curves \eqref{kjf3204fjnr} [green zone]. Both curves start from the point $z=\frac{1}{6}$ 
and have vertical asymptotes in the complex $z$--plane at the line $\Re z=\frac{1}{4}$. 
The convergence of the series in the vertical strip $\frac{1}{6}<\Re z<\frac{1}{4}$ depends, therefore,
on the imaginary part of $z$. On the contrary, in the half--plane $\Re z\geqslant\frac{1}{4}$ both series converge everywhere
independently of the imaginary part of $z$. }
\label{kawfe3o84yg}
\end{figure}

Consider the general form of the second Binet's integral formula for the logarithm of the $\Gamma$--function
\be\label{pio2u092je3n4}
\int\limits_0^{\,\infty} \!\frac{\,\arctg{a x}\,}{e^{b x}-1}\, dx \,=\,
\frac{\,\pi\,}{b} \ln\Gamma\!\left(\!\dfrac{b}{\,2\pi a\,}\!\right) 
+\frac{1}{\,2a\,}\!\left(1-\ln\frac{b}{2\pi a} \right) + \frac{\,\pi\,}{\,2b\,}\ln\frac{b}{4\pi^2 a}
\ee
$a>0\,$ and $\Re{b}>0\,$, 
see e.g.~\cite[vol.~I, \no 2.7.5-6]{prudnikov_en},
\cite[pp.~335--336]{binet_01}, \cite[pp.~250--251]{whittaker_01}, \cite[vol.~I, p.~22, Eq.~1.9(9)]{bateman_01}
or \cite[Sect.~4, exercise \no 40]{iaroslav_06}. The general idea of the method consists in finding such a change of variable 
that reduces the integrand in the left--hand side of \eqref{pio2u092je3n4} to a function (probably, a composite function) which may be ``easily'' expanded into the MacLaurin series. 
In our case, this change of variable may be easily found by requiring, for example, that
\be\notag
\int\!\frac{\,dx\,}{e^{b x}-1}\,  = \,\alpha\!\int\!\frac{\,du\,}{u} 
\ee
where $u$ is the new variable and $\alpha$ is some normalizing coefficient, which can be chosen later at our convenience. 
Other changes of variables, of course, are possible as well (see, e.g., numerous examples in exercises 39 \& 45
\cite[Sect.~4]{iaroslav_06}), but this one is particularly successful, especially if we set $\alpha=1/b$.
Thus, putting \mbox{$\,x=-\frac{1}{b}\ln(1-u)\,$}  
and rewriting the result for $\,z=\frac{b}{2\pi a}\,$, Binet's formula takes the form
\be\label{k034fjdfnsh1}
\begin{array}{ll}
\displaystyle
-\!\int\limits_0^1 \!\arctg{\!\Big[\frac{1}{\,2\pi z\,}\ln(1-u)\Big]} \frac{\,du\,}{u} \,=\,
\pi \ln\Gamma(z)+ \pi z(1-\ln z) + \frac{\,\pi\,}{\,2\,}\ln\frac{z}{\,2\pi\,}
\end{array}
\ee
where $\Re z>0$.
The integrand on the left may be expanded into the MacLaurin series in powers of $u$ accordingly to the method
described in Section \ref{stirling}. This yields
\be\label{k034fjdfnsh2}
\begin{array}{ll}
\displaystyle
\arctg{\!\Big[\frac{1}{\,2\pi z\,}\ln(1-u)\Big]} &\displaystyle\,=\,
\sum_{l=0}^\infty (-1)^l (2l)!\cdot\frac{\left[\frac{1}{\,2\pi z\,}\ln(1-u)\right]^{2l+1}}{(2l+1)!}\\[5mm]
&\displaystyle 
=\,-\!\sum_{l=0}^\infty (-1)^{l} \frac{(2l)!}{\,(2\pi z)^{2l+1}}\,\cdot\!\!\!\!
\sum_{n=2l+1}^\infty\!\!\!\!\frac{\,|S_1(n,2l+1)|\,}{\,n!\,}\, u^{n}\\[5mm]
&\displaystyle 
=\,-\!\sum_{n=1}^\infty\!\frac{u^{n}}{\,n!\,}\cdot\!\sum_{l=0}^{\lfloor\frac{1}{2}n\rfloor} \!(-1)^{l} 
 \frac{(2l)!\cdot|S_1(n,2l+1)|}{\,(2\pi z)^{2l+1}}
\end{array}
\ee
This expansion converges in the disk $|u|<r$ in which $\,\arctg{\!\big[\frac{1}{\,2\pi z\,}\ln(1-u)\big]}$ is analytic.
The radius of this disk $r$ depends on the parameter $z$ and is conditioned by the singularities of the arctangent, 
which occur at $u=1-\exp(\pm 2\pi i z)$ [branch points], and by that of the logarithm, which is located at $u=1$
[branch point as well]. The latter restricts the value of $r$ to 1, and the unit radius of convergence 
corresponds to such $z$ that $\,2\cos(2\pi\Re z\big)=\exp(\pm2\pi\Im z)$.
The zone of convergence of series \eqref{k034fjdfnsh2} for $|u|<1$ consists, therefore, in 
the intersection of two zones, each of which lying to the right of curves
\be\label{kjf3204fjnr}
\Im{z}\,=\,\frac{1}{\,2\pi\,}\cdot
\begin{cases}
+\ln2 + \ln\cos\big(2\pi\Re z\big)\\[1mm]
-\ln2 - \ln\cos\big(2\pi\Re z\big)\\[1mm]
\end{cases}
\ee
respectively, see Fig.~\ref{kawfe3o84yg}.
Now, a close study of the general term of series \eqref{k034fjdfnsh2} reveals that it also converges for $u=1$.
Indeed, from \eqref{i8ud034jdn}, it follows that one can always find such a constant $C>1$ that for sufficiently large $n_0$,
inequality
\be\notag
\frac{1}{\,n!\,}\!\sum_{l=0}^{\lfloor\frac{1}{2}n\rfloor}\! (-1)^{l} 
\frac{(2l)!\cdot|S_1(n,2l+1)|}{\,(2\pi z)^{2l+1}} \,
\leqslant\,\frac{2\pi z \, C}{\,n\ln^2\! n\,}\,,
 \qquad\quad \text{} n>n_0\,,
\ee
holds.
Hence, since series $\,\sum n^{-1}\!\ln^{-\alpha}\!n$ converges for $\alpha>1$,
so does series \eqref{k034fjdfnsh2} at $u=1$.
An interesting consequence of the latter statement is this curious identity
\be\label{jd2903djd}
\frac{\pi}{2}\,=
\sum_{n=1}^\infty\!\frac{1}{\,n!\,}\cdot\!\sum_{l=0}^{\lfloor\frac{1}{2}n\rfloor}\! (-1)^{l} 
\frac{(2l)! \cdot|S_1(n,2l+1)|}{\,(2\pi z)^{2l+1}}\,,
\ee
which holds in the region of convergence of $z$.
Thus, expansion \eqref{k034fjdfnsh2} converges uniformly in each point of the disk $|u|<1$ and can be integrated term--by--term.\footnote{Another way to show that \eqref{k034fjdfnsh2} is uniformly convergent is to directly verify that 
\be\notag
\int\limits_0^1 \!\left[\sum_{n=N}^\infty\! \frac{u^n}{\,n!\,}\!\!\sum_{l=0}^{\lfloor\frac{1}{2}n\rfloor}\! (-1)^{l} 
\frac{(2l)!\cdot|S_1(n,2l+1)|}{\,(2\pi z)^{2l+1}}\right] du \,\to\,0  \qquad\quad \text{as }\, N\to\infty
\ee
see e.g.~\cite[pp.~161--162]{gunter_02_eng}. }
Substituting series \eqref{k034fjdfnsh2} into \eqref{k034fjdfnsh1} and performing the indicated term--by--term integration from $u=0\,$ to $\,u=1$,
we obtain the following series expansion for the logarithm of the $\Gamma$--function
\be\label{kj20ejcn2dnd}
\begin{array}{ll}
\displaystyle
\ln\Gamma(z)\; & \displaystyle =\left(z-\frac{1}{\,2\,}\right)\!\ln z -z +\frac{1}{\,2\,}\ln2\pi + \frac{1}{\,\pi\,}\!
\sum_{n=1}^\infty\!\frac{ 1}{\,n\cdot n!\,}\!\sum_{l=0}^{\lfloor\frac{1}{2}n\rfloor}\! (-1)^{l} 
\frac{\, (2l)!\cdot|S_1(n,2l+1)|\,}{(2\pi z)^{2l+1}} \\[4mm]
& \displaystyle 
= \left(z-\frac{1}{\,2\,}\right)\!\ln z -z +\frac{1}{\,2\,}\ln2\pi + \frac{1}{\,\pi\,}
\left\{\frac{1}{2\pi z}+\frac{1}{8\pi z}+\frac{1}{18}\left(\frac{1}{\pi z}-\frac{1}{4\pi^3 z^3}\right)\right.\\[4mm]
& \displaystyle \quad
+\frac{3}{96}\left(\frac{1}{\pi z}-\frac{1}{2\pi^3 z^3}\right)
+\frac{1}{600}\left(\frac{12}{\pi z}-\frac{35}{4 \pi^3 z^3}+\frac{3}{4 \pi^5 z^5}\right)\\[4mm]
& \displaystyle \quad\left.
+\frac{1}{4320}\left(\frac{60}{\pi z}-\frac{225}{4 \pi^3 z^3}+\frac{45}{4 \pi^5 z^5}\right)+\ldots\right\}
\end{array}
\ee
converging in the same region as series \eqref{k034fjdfnsh2}, see \eqref{kjf3204fjnr} and Fig.~\ref{kawfe3o84yg}.
In particular, if $z$ is real, it converges for $z>\frac{1}{6}$; on the contrary, if $z$ is complex,
then independently of its imaginary part, it converges everywhere in the right half--plane $\Re z\geqslant\frac{1}{4}$.
A quick analysis of the above series shows that for $z$ rational multiple of $\pi^{-1}$, it
contains rational coefficients only. Another important observation is that
this series, unlike the classic Stirling series \eqref{jxc2mxp2mxdx}, cannot be explicitly written in powers of $z$.
To illustrate this point, we write down its first 
2, 3 and 4 terms respectively:
\be\notag
\begin{array}{ll}
\displaystyle
\sum_{n=1}^N\!\frac{ 1}{\,n\cdot n!\,}\!\sum_{l=0}^{\lfloor\frac{1}{2}n\rfloor} \!(-1)^{l} 
\frac{\, (2l)!\cdot|S_1(n,2l+1)|\,}{(2\pi z)^{2l+1}} = \\[6mm]
\displaystyle
\qquad \,=\,
\begin{cases}
\displaystyle
\frac{1}{2\pi z}+\frac{1}{8\pi z}  \,=\, \frac{5}{8\pi z}\,, \qquad & N=2\\[4mm]
\displaystyle
\frac{1}{2\pi z}+\frac{1}{8\pi z}+\frac{1}{18}\left(\frac{1}{\pi z}-\frac{1}{4\pi^3 z^3}\right)
\,=\,\frac{49}{72\pi z}-\frac{1}{72\pi^3z^3}  \,, \qquad & N=3\\[4mm]
\displaystyle
\frac{1}{2\pi z}+\frac{1}{8\pi z}+\frac{1}{18}\left(\frac{1}{\pi z}-\frac{1}{4\pi^3 z^3}\right)
+\frac{3}{96}\left(\frac{1}{\pi z}-\frac{1}{2\pi^3 z^3}\right)
\,=\,\frac{205}{288\pi z}-\frac{17}{576\pi^3z^3}  \,, \qquad & N=4
\end{cases}
\end{array}
\ee
By the way, as concerns the divergent Stirling series \eqref{jxc2mxp2mxdx}, it can be readily derived from 
\eqref{kj20ejcn2dnd}.
By formally interchanging sum signs in \eqref{kj20ejcn2dnd},
which is obviously not permitted because series are not absolutely convergent, we have
\be\notag
\begin{array}{cc}
\displaystyle
\sum_{n=1}^\infty\!\frac{ 1}{\,n\cdot n!\,}\!\sum_{l=0}^{\lfloor\frac{1}{2}n\rfloor}\! (-1)^{l} 
\frac{\, (2l)!\cdot|S_1(n,2l+1)|\,}{(2\pi z)^{2l+1}} \,\asymp\,
\sum_{l=0}^{\infty}  (-1)^{l} \frac{\, (2l)!\,}{(2\pi z)^{2l+1}} \underbrace{
\sum_{n=1}^\infty\!\frac{|S_1(n,2l+1)|}{\,n\cdot n!\,}}_{\zeta(2l+2)} \\
\displaystyle
= \sum_{l=1}^{\infty}  (-1)^{l-1} \frac{\, (2l-2)!\,}{(2\pi z)^{2l-1}} \cdot\zeta(2l) \,=
\sum_{l=1}^\infty \frac{\pi\cdot{B}_{2l}}{2l(2l-1)z^{2l-1}} 
\end{array}
\ee
where we first used \eqref{kljd023jdnr} for $\zeta(2l+2)$, and then, Euler's formula 
\be\notag
\zeta(2l)=(-1)^{l+1}\frac{\,(2\pi)^{2l}\cdot{B}_{2l}\,}{2\cdot(2l)!}\,=\,
\frac{\,(2\pi)^{2l}\cdot|{B}_{2l}|\,}{2\cdot(2l)!}\,,\qquad\quad l=1, 2, 3,\ldots
\ee 
Further observations concern the convergence of the derived series and 
are treated in details in the next section.

\subsubsection{Convergence analysis of the derived series}\label{wkc2me23dwe}
The complete study of the convergence of \eqref{kj20ejcn2dnd} is quite long and complicated,
that is why we split it in two stages. First, we 
obtain the upper bound for the general term of \eqref{kj20ejcn2dnd},
and then, derive an accurate approximation for it when $n$ becomes sufficiently large.
In what follows, we may suppose, without essential loss of generality, that $z$ is real and positive.
The general term of series \eqref{kj20ejcn2dnd} is given by the finite sum over $l$.
This truncated sum has only odd terms, and hence, by elementary transformations, may be reduced 
to that containing both odd and even terms
\be\label{lkh908gb9gi8vityr}
\begin{array}{ll}
\displaystyle
& \displaystyle 
\!\sum_{l=0}^{\lfloor\frac{1}{2}n\rfloor}\! (-1)^{l} 
\frac{\, (2l)!\cdot|S_1(n,2l+1)|\,}{(2\pi z)^{2l+1}} = \!\sum_{l=0}^{\lfloor\frac{1}{2}n\rfloor}\! (-1)^{\frac{1}{2}(2l+1)-\frac{1}{2}}
\frac{\, (2l+1)!\cdot|S_1(n,2l+1)|\,}{(2l+1)\cdot(2\pi z)^{2l+1}}\\[5mm]
& \displaystyle \qquad\qquad\qquad
= \,\frac{1}{2}\!\sum_{l=1}^{n} \big[1-(-1)^{l}\big] \cdot(-1)^{\frac{1}{2}(l-1)}\cdot\frac{(l-1)!\cdot|S_1(n, l)|}{(2\pi z)^l} \,=\,\ldots
\end{array}
\ee
Now, from Legendre's integral for the Euler $\Gamma$--function,\footnote{Namely, 
$\,(l+1)!=\int\! x^{l+1} e^{-x} dx\,$ taken over $x=[0, \infty)$.} it follows that
\be\notag
\left\{
\begin{array}{ll}
\displaystyle
(-1)^{\frac{1}{2}(l-1)}\cdot\frac{\,(l-1)!\,}{(2\pi z)^l}\,=\,-i\!\int\limits_0^\infty\! \left[\frac{i\,x}{2\pi z}\right]^l\!\cdot\frac{\,e^{-x}\,dx\,}{x}\\[6mm]
\displaystyle
(-1)^l\cdot(-1)^{\frac{1}{2}(l-1)}\cdot\frac{\,(l-1)!\,}{(2\pi z)^l}\,=\,-i\!\int\limits_0^\infty\! \left[-\frac{i\,x}{2\pi z}\right]^l\!\cdot\frac{\,e^{-x}\,dx\,}{x}
\end{array}
\right.
\ee
Hence, expression \eqref{lkh908gb9gi8vityr} may be continued as follows
\be\label{lkhfveovede}
\begin{array}{ll}
\displaystyle
\ldots\, & \displaystyle 
=\, \frac{i}{2}\!\int\limits_0^\infty\! \sum_{l=1}^{n} \left\{ \left[-\frac{i\,x}{2\pi z}\right]^l - 
 \left[\frac{i\,x}{2\pi z}\right]^l \right\} |S_1(n, l)|\cdot \frac{\,e^{-x}\,dx\,}{x}   \\[6mm]
& \displaystyle 
=\, \frac{i}{2}\!\int\limits_0^\infty\! \left\{\! \left(-\frac{i\,x}{2\pi z}\right)_{\!n} - 
 \left(\frac{i\,x}{2\pi z}\right)_{\!n} \!\right\} \frac{\,e^{-x}\,dx\,}{x}   
\\[5mm]
& \displaystyle 
=\, \frac{i}{\,4\pi^2 z\,}\!\int\limits_0^\infty\! \sh\frac{x}{2z}\left\{ 
\Gamma\!\left(\frac{i\,x}{2\pi z}\right) \Gamma\!\left(n-\frac{i\,x}{2\pi z}\right) 
- \Gamma\!\left(-\frac{i\,x}{2\pi z}\right) \Gamma\!\left(n+\frac{i\,x}{2\pi z}\right)  \!\right\} 
e^{-x}\,dx  \\[6mm]
\displaystyle
& \displaystyle 
=\, -\frac{1}{\,2\pi^2 z\,}\!\int\limits_0^\infty\! \sh\frac{x}{2z}\cdot e^{-x}\cdot\Im\left[ 
\Gamma\!\left(\frac{i\,x}{2\pi z}\right) \Gamma\!\left(n-\frac{i\,x}{2\pi z}\right)  \!\right] dx 
\end{array}
\ee
where at the final stage we, first, replaced Pochhammer symbols by $\Gamma$--functions, and
then, used the well--known relationship $\Gamma(z)\Gamma(-z)\,=\,-(\pi/z)\csc\pi z\,$.
The last integral in \eqref{lkhfveovede} is difficult to evaluate in a closed form, but its upper bound may be readily obtained. 
In view of the fact that $\,|\Im\Gamma(v)|\leqslant|\Gamma(v)|\leqslant|\Gamma(\Re v)|\,$,
we have
\be\label{jlhoiybnlnb}
\begin{array}{ll}
\displaystyle
& \displaystyle 
\frac{1}{\,2\pi^2 z\,}\left|\int\limits_0^\infty\! \sh\frac{x}{2z}\cdot e^{-x}\cdot\Im\left[ 
\Gamma\!\left(\frac{i\,x}{2\pi z}\right) \Gamma\!\left(n-\frac{i\,x}{2\pi z}\right)  \!\right] dx  \right|\,\leqslant  \\[6mm]
& \displaystyle \qquad
\leqslant\,\frac{\Gamma(n)}{\,2\pi^2 z\,}\int\limits_0^\infty\! \sh\frac{x}{2z}\cdot e^{-x}\cdot 
\left|\Gamma\!\left(\frac{i\,x}{2\pi z}\right) \!\right| dx \,
=\, \frac{\,(n-1)!\,}{\,\pi\sqrt{2z}\,}\int\limits_0^\infty\! e^{-x}\sqrt{\sh\frac{x}{2z}\,} \cdot\frac{\,dx\,}{\sqrt{x}} 
\end{array}
\ee 
Whence, by making a change of variable in the latter integral $\,x=2zt$, we have for any positive integer $n$ (not necessarily large) 
\be\label{k039dm3dmedc}
\frac{1}{\,n\cdot n!\,}\left|\sum_{l=0}^{\lfloor\frac{1}{2}n\rfloor}\! (-1)^{l} 
\frac{\, (2l)!\cdot|S_1(n,2l+1)|\,}{(2\pi z)^{2l+1}}\right| \leqslant\, \frac{1}{\,n^2\,}\cdot\frac{\,1\,}{\,\pi\,}\!
\int\limits_0^\infty\! \sqrt{\frac{\sh t\,}{\,t\,}}\cdot e^{-2zt}\,dt
\ee
where the latter integral converges uniformly in the half--plane $z$ which lies to the right of the line
\mbox{$\,\Re z=\frac{1}{4}$} (imaginary part of $z$ contributes only to the bounded oscillations
of the integrand). Consequently, 
series \eqref{kj20ejcn2dnd} converges at least in $\Re z>\frac{1}{4}$, and this at the same rate 
or better than Euler's series $\,\sum n^{-2}$.

\begin{figure}[!t]   
\centering
\includegraphics[width=0.8\textwidth]{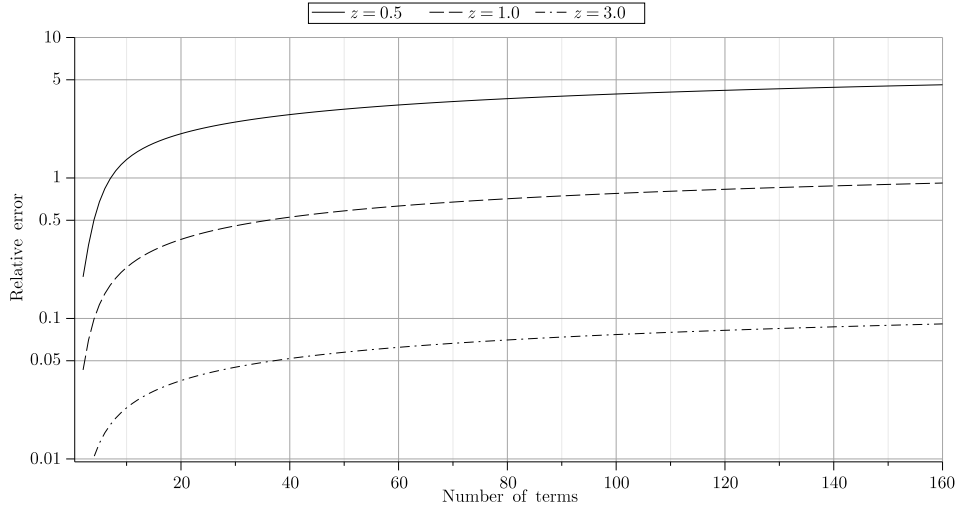}
\caption{Relative error between the upper bound and the left--hand side in \eqref{k039dm3dmedc} as a function of $n$ 
for three different values of argument $z$, logarithmic scale.}
\label{cj2394chcned}
\end{figure}

Numerical simulations show, however, that the greater $n$, the greater the relative difference between the upper bound and the left--hand side
in \eqref{k039dm3dmedc}, see Fig.~\ref{cj2394chcned},
and thus, this upper bound is relatively rough.\footnote{The error is mainly
due to the use of inequality  $\,|\Im\Gamma(v)|\leqslant|\Gamma(v)|\,$.}
A more accurate description of the behavior of
sum \eqref{lkh908gb9gi8vityr} at large $n$ may be obtained by seeking its asymptotics. In order to find it, we proceed as follows.
We first rewrite the second line of \eqref{lkhfveovede} as
\be\label{jd2099dj23lnd}
\frac{i}{2}\!\int\limits_0^\infty\! \left\{\! \left(-\frac{i\,x}{2\pi z}\right)_{\!n} - 
 \left(\frac{i\,x}{2\pi z}\right)_{\!n} \!\right\} \frac{\,e^{-x}\,dx\,}{x}   
 =\, \int\limits_0^\infty\! \Im\left[\!  \left(\frac{i\,x}{2\pi z}\right)_{\!n}\right]\frac{\,e^{-x}\,dx\,}{x} 
\ee
Now, it is well--known that function 
$1/\Gamma(z)$ is regular
on the entire complex $z$--plane, and therefore, may be expanded into the MacLaurin's series
\be\label{poi2d293dm}
\frac{1}{\Gamma(z)}\,=\,z+\gamma z^2+ \left(\! \frac{\gamma^2}{2}-\frac{\pi^2}{12}\right)\!z^3 +
\ldots\,\equiv\sum_{k=1}^\infty z^k a_k \,,\qquad |z|<\infty\,,
\ee
where
\be\label{jkhe9843hdhr}
a_k\,\equiv\,\frac{1}{k!}\cdot \left[\frac{1}{\Gamma(z)}\right]^{(k)}_{z=0} \!=\, 
\frac{(-1)^k}{\,\pi \, k!\,}\cdot \Big[\sin\pi x\cdot\Gamma(x)\Big]^{(k)}_{x=1} 
\ee 
the last representation for coefficients $a_k$, which follows from the reflection formula for the $\Gamma$--function,
being often more suitable for computational purposes.\footnote{On the computation of $a_k$,
see also \cite[p.~256, \no 6.1.34]{abramowitz_01}, \cite[pp.~344 \& 349]{wilf_02}, \cite{hayman_01}.} 
Using approximation \eqref{lk2093mffmnjw} for the Pochhammer symbol and the above MacLaurin series 
for $1/\Gamma(z)$, we have for sufficiently large $n$
\be\label{jlhoiybnlnb}
\begin{array}{ll}
\displaystyle
\Im\left[\!  \left(\frac{i\,x}{2\pi z}\right)_{\!n}\right]\sim\,\Im\left[ \frac{n^{\frac{ix}{2\pi z}}\cdot\Gamma(n)}
{\Gamma\left(\frac{ix}{2\pi z}\right)}\right]  = \, (n-1)!\Im\left[\!\left(\cos\frac{x\ln n}{2\pi z}+i\sin\frac{x\ln n}{2\pi z}\right)
\!\cdot\!\sum_{k=1}^\infty a_k \cdot\!\left(\frac{i\,x}{2\pi z}\right)^{\!k}\right]\\[8mm]
\displaystyle\qquad\qquad
=\,(n-1)! \left[\cos\frac{x\ln n}{2\pi z}\cdot\!\sum_{k=0}^\infty (-1)^k  a_{2k+1} \left(\frac{x}{2\pi z}\right)^{\!2k+1}  
\!+ \, \sin\frac{x\ln n}{2\pi z}\cdot\!\sum_{k=1}^\infty  (-1)^k a_{2k} \left(\frac{x}{2\pi z}\right)^{\!2k} \right]
\end{array}
\ee 
the error due to considering only the first term in \eqref{lk2093mffmnjw} being negligible with respect to logarithmic terms
which will appear later.
Inserting this expression into \eqref{jd2099dj23lnd}, performing the term--by--term integration and taking into account that\footnote{These 
equalities are valid wherever the integrals on the left converge, see e.g.~\cite[p.~130]{cauchy_01}, 
\cite[p.~12]{malmsten_01}, \cite{kratzer_01}.}
\be
\begin{cases}
\displaystyle
\int\limits_0^\infty \! x^{s-1} e^{-zx} \cos ux\,dx\,=\, \frac{\Gamma(s)}
{\big(z^2+u^2\big)^{s/2}}\cdot\cos\left[s\arctg\frac{u}{z}\right] \\[6mm]
\displaystyle
\int\limits_0^\infty \! x^{s-1} e^{-zx} \sin ux\,dx\,=\, \frac{\Gamma(s)}
{\big(z^2+u^2\big)^{s/2}}\cdot\sin\left[s\arctg\frac{u}{z}\right]
\end{cases}
\ee
yields
\be\notag
\begin{array}{ll}
\displaystyle
\int\limits_0^\infty\! \Im\left[\!  \left(\frac{i\,x}{2\pi z}\right)_{\!n}\right]\frac{\,e^{-x}\,dx\,}{x} \,\sim\,
(n-1)! \sum_{k=0}^\infty (-1)^k  a_{2k+1} \frac{(2k)!}{(2\pi z)^{2k+1}} 
\cdot\frac{\cos\left[(2k+1)\arctg\dfrac{\ln n}{2\pi z}\right]}{\left[1+\dfrac{\ln^2\! n}{4\pi^2 z^2}\right]^{k+\frac{1}{2}}} \\[12mm]
\displaystyle\qquad\qquad\qquad\qquad\qquad\qquad\qquad
+ (n-1)! \sum_{k=1}^\infty (-1)^k  a_{2k} \frac{(2k-1)!}{(2\pi z)^{2k}} 
\cdot\frac{\sin\left[2k\arctg\dfrac{\ln n}{2\pi z}\right]}{\left[1+\dfrac{\ln^2\! n}{4\pi^2 z^2}\right]^k}
\end{array}
\ee
Whence, the required asymptotics is
\be\label{jxc2803dn}
\begin{array}{ll}
\displaystyle
\frac{1}{\,n\cdot n!\,}\!\!\sum_{l=0}^{\lfloor\frac{1}{2}n\rfloor}\! (-1)^{l} 
\frac{\, (2l)!\cdot|S_1(n,2l+1)|\,}{(2\pi z)^{2l+1}}\,=   \\[6mm]
\displaystyle\qquad\qquad\quad
= \frac{1}{\,n^2\,}\sum_{k=0}^\infty  \frac{(-1)^k}{\,(2k+1)\big(4\pi^2 z^2+\ln^2\! n\big)^{k+\frac{1}{2}}\,} 
\cdot\cos\left[(2k+1)\arctg\dfrac{\ln n}{2\pi z}\right]\cdot\left[\frac{1}{\Gamma(x)}\right]^{(2k+1)}_{x=0} \\[10mm]
\displaystyle\qquad\qquad\qquad\quad
+ \frac{1}{\,n^2\,}\sum_{k=1}^\infty  \frac{(-1)^k}{\,2k\big(4\pi^2 z^2+\ln^2\! n\big)^k\,} 
\cdot\sin\left[2k\arctg\dfrac{\ln n}{2\pi z}\right]\cdot\left[\frac{1}{\Gamma(x)}\right]^{(2k)}_{x=0}
+O\!\left(\frac{1}{n^3}\right)
\end{array}
\ee
for sufficiently large $n$. Retaining first few terms, we have 
\be\label{i8ud034jdn}
\begin{array}{ll}
\displaystyle
\frac{1}{\,n\cdot n!\,}\!\!\sum_{l=0}^{\lfloor\frac{1}{2}n\rfloor} \!(-1)^{l} 
\frac{\, (2l)!\cdot|S_1(n,2l+1)|\,}{(2\pi z)^{2l+1}}\,=  \,\frac{\,2\pi z\,}{n^2}\left\{ \frac{1}{\,4\pi^2 z^2+\ln^2\! n\,}
- \frac{2\gamma\ln n}{\,\big(4\pi^2 z^2+\ln^2\! n\big)^2\,} \right\}  \\[6mm]
\displaystyle\qquad\qquad\qquad\qquad\qquad\qquad\qquad\qquad\qquad
+\, O\!\left(\frac{1}{n^2\ln^4\! n}\right)
\,, \qquad n\to\infty
\end{array}
\ee
Thus, for moderate values of $z$, series \eqref{kj20ejcn2dnd} converges approximately at the same rate as $\,\sum (n\ln n)^{-2}$,
i.e.~at the same rate as, for example, Fontana--Mascheroni's series \eqref{fms}, \eqref{fms2x}, see asymptotics \eqref{j38ndbr893r}.

\subsubsection{Some important particular cases of the derived series}
\begin{figure}[!t]   
\centering
\includegraphics[width=0.8\textwidth]{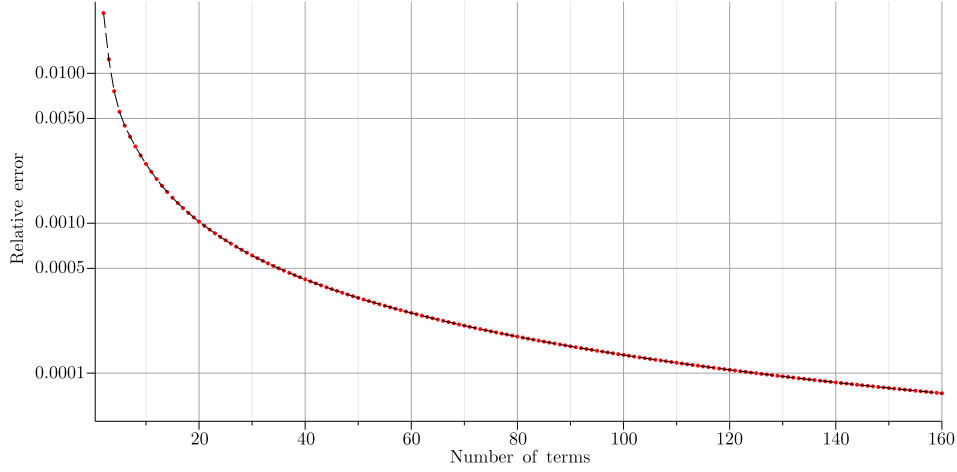}
\caption{Relative error of the series expansion for $\ln\Gamma(\pi^{-1})$ given by \protect\eqref{lkidj109jmndx3}, logarithmic scale.}
\label{jkerg45hgs}
\end{figure}
\begin{figure}[!t]   
\centering
\includegraphics[width=0.8\textwidth]{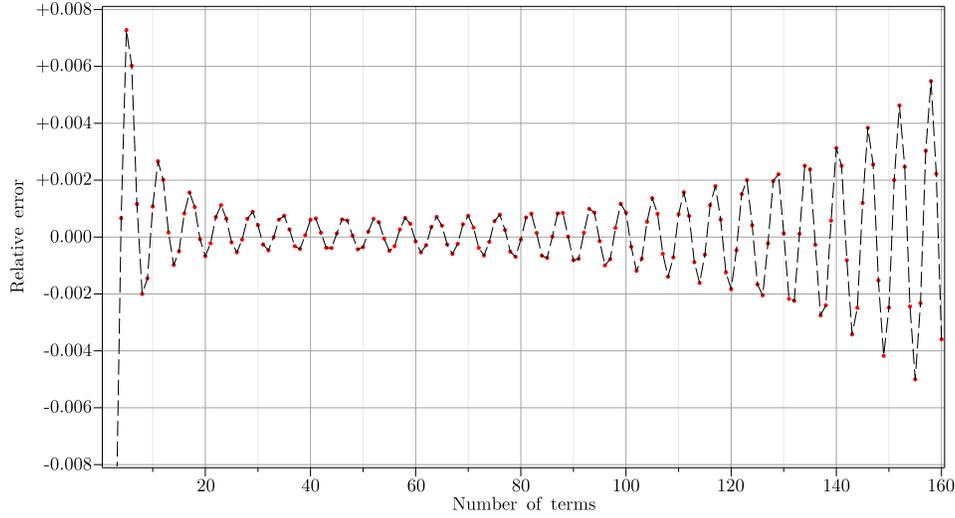}
\caption{Relative error of the series expansion for $\ln\Gamma\big(\frac{1}{2}\pi^{-1}\big)$ given by \protect\eqref{po30im2d}, linear scale.}
\label{jkcdh73b9s}
\end{figure}

Let now consider some applications of the formula \eqref{kj20ejcn2dnd}. In the first instance, it is natural to obtain 
a series expansion for 
\be\label{lkidj109jmndx3}
\begin{array}{ll}
\displaystyle
\ln\Gamma\!\left(\!\dfrac{1}{\,\pi \,}\!\right) \, &\displaystyle =
\left(\!1-\frac{1}{\pi}\!\right)\!\cdot\ln\pi -\frac{1}{\pi} +\frac{1}{2}\ln2
+\frac{1}{2\pi} \!\sum_{n=1}^\infty\!\frac{1}{\,n\cdot n!\,}\cdot\Biggl\{\!\sum_{l=0}^{\lfloor\frac{1}{2}n\rfloor} \!(-1)^{l} 
\frac{(2l)!\cdot |S_1(n,2l+1)|}{2^{2l}}\Biggr\} \\[6mm]
&\displaystyle \!\!\!\!\!\!\!\!\!\!\!\!\!\!\!\!\!\!\!\!\!
=\left(\!1-\frac{1}{\pi}\!\right)\!\cdot\ln\pi -\frac{1}{\pi} +\frac{1}{2}\ln2 +\frac{1}{2\pi}
\left\{1+ \frac{1}{4}+ \frac{1}{12}+  \frac{1}{32}+ \frac{1}{75}+ \frac{1}{144} 
+ \frac{13}{2880} + \frac{157}{46\,080} +\ldots \right\}
\end{array}
\ee
The graphical illustration of the convergence of this series is given in Fig.~\ref{jkerg45hgs}.
With equal ease, we derive 
\be\label{jnx0q93j2onix2l3dcn}
\begin{array}{ll}
\displaystyle
\ln\Gamma\!\left(\!\dfrac{2}{\,\pi \,}\!\right) \, &\displaystyle =
\left(\!1-\frac{2}{\pi}\!\right)\!\cdot\ln\pi -\frac{2}{\pi} +\frac{2}{\pi}\ln2
+\frac{1}{4\pi} \!\sum_{n=1}^\infty\!\frac{1}{\,n\cdot n!\,}\cdot\Biggl\{\!\sum_{l=0}^{\lfloor\frac{1}{2}n\rfloor}\! (-1)^{l} 
\frac{(2l)!\cdot |S_1(n,2l+1)|}{2^{4l}}\Biggr\} \\[6mm]
&\displaystyle \!\!\!\!\!\!\!\!\!\!\!\!\!\!\!\!\!\!\!\!\!\!\!\!\!\!\!\!\!\!
=\left(\!1-\frac{2}{\pi}\!\right)\!\cdot\ln\pi -\frac{2}{\pi} +\frac{2}{\pi}\ln2+\frac{1}{4\pi}
\left\{1+ \frac{1}{4}+ \frac{5}{48}+  \frac{7}{128}+ \frac{631}{19\,200}+ \frac{199}{9216} 
+ \frac{19\,501}{1\,290\,240} + \frac{32\,707}{2\,949\,120} +\ldots\right\}
\end{array}
\ee
Many other similar expansions may be derived analogously.
Let now see how the series behaves outside the region of convergence.
For this aim, we take $z=\frac{1}{2}\pi^{-1}$. Formula \eqref{kj20ejcn2dnd} yields 
\be\label{po30im2d}
\begin{array}{ll}
\displaystyle
\ln\Gamma\!\left(\!\dfrac{1}{\,2\pi \,}\!\right) \, &\displaystyle \stackrel{?}{=}
\left(\!1-\frac{1}{2\pi}\!\right)\!\cdot\ln2\pi -\frac{1}{2\pi} +\frac{1}{\pi}
\sum_{n=1}^\infty\!\frac{1}{\,n\cdot n!\,}\cdot\Biggl\{\!\sum_{l=0}^{\lfloor\frac{1}{2}n\rfloor}\! (-1)^{l} 
(2l)!\cdot |S_1(n,2l+1)|\Biggr\}    \\[6mm]
&\displaystyle 
=\left(\!1-\frac{1}{2\pi}\!\right)\!\cdot\ln2\pi -\frac{1}{2\pi} +\frac{1}{\pi}
\left\{ 1+ \frac{1}{4}- \frac{1}{16}-  \frac{11}{300}+ \frac{1}{144}+ \frac{17}{630} 
+ \frac{101}{5760} - \frac{311}{102\,060} - \ldots \right\}
\end{array}
\ee
At first sight, it might seem that this alternating series slowly converges to $\ln\Gamma(\frac{1}{2}\pi^{-1})\approx
1.765383194$: the summation 
of its first 3 terms gives the value $1.764207893\ldots$ which corresponds to the relative accuracy $6.6\!\times\!10^{-4}$,
that of 18 terms gives $1.765525087\ldots$, i.e.~the relative accuracy $8.0\!\times\!10^{-5}$, 
the summation of first 32 terms yields $1.765392783\ldots$ which corresponds to the relative error $5.4\!\times\!10^{-6}$.\footnote{We do 
not count the third term which is zero.}
Notwithstanding, further numerical simulations, see Fig.~\ref{jkcdh73b9s},
leave no doubts: this series is divergent.

\subsection{Second series expansion for the logarithm of the $\Gamma$--function}
Rewrite formula \eqref{kj20ejcn2dnd} for $2z$ instead of $z$, and subtract
the result from \eqref{kj20ejcn2dnd}. In virtue of Legendre's duplication formula for the $\Gamma$--function
\mbox{$\ln\Gamma(2z)=(2z-1)\ln2-\frac{1}{2}\ln\pi +\ln\Gamma(z)+\ln\Gamma(z+\frac{1}{2})$}, we have
\be\label{jkhc92h8cn23d}
\begin{array}{ll}
\displaystyle
\ln\Gamma\!\left(\frac{1}{2}+z\right)\: & \displaystyle =\,z\ln{z}- z +\frac{1}{\,2\,}\ln2\pi - \frac{1}{\,\pi\,}\!
\sum_{n=1}^\infty\!\frac{ 1}{\,n\cdot n!\,}\!\sum_{l=0}^{\lfloor\frac{1}{2}n\rfloor}\! (-1)^{l} 
\frac{\, (2l)!\cdot(2^{2l+1}-1)\cdot|S_1(n,2l+1)|\,}{(4\pi z)^{2l+1}} \\[5mm]
& \displaystyle \!\!\!\!\!\!\!\!\!\!\!\!\!\!
= \,z\ln{z}- z +\frac{1}{\,2\,}\ln2\pi -  \frac{1}{\,\pi\,}
\left\{\frac{1}{4\pi z}+\frac{1}{16\pi z}+\frac{1}{18}\left(\frac{1}{2\pi z}-\frac{7}{32\pi^3 z^3}\right)
+\frac{3}{96}\left(\frac{1}{2\pi z}-\frac{7}{16\pi^3 z^3}\right)\right.\\[5mm]
& \displaystyle \!\!\!\!\!\!\!\!\!\!\!\!\!\! \quad\left.
+\frac{1}{600}\left(\frac{6}{\pi z}-\frac{245}{32 \pi^3 z^3}+\frac{93}{128 \pi^5 z^5}\right)
+\frac{1}{4320}\left(\frac{30}{\pi z}-\frac{1575}{32 \pi^3 z^3}+\frac{1395}{128 \pi^5 z^5}\right)+\ldots\right\}
\end{array}
\ee
which holds in the green zone shown in Fig.~\ref{kawfe3o84yg}. This expression allows to expand any value of the form
$\ln\Gamma(\frac{1}{2}+\alpha\pi^{-1})$ into the series 
with rational coefficients if $\alpha$ is rational greater than $\frac{1}{6}\pi$. For example, putting $z=\pi^{-1}$, we have
\begin{eqnarray}
\displaystyle\notag
 \ln\Gamma\!\left(\!\dfrac{1}{2}+\dfrac{1}{\,\pi\,}\!\right)\; && \displaystyle = \, -\frac{1+\ln\pi}{\pi} +\frac{1}{\,2\,}\ln2\pi - \frac{1}{\,4\pi\,}\!
\sum_{n=1}^\infty\!\frac{ 1}{\,n\cdot n!\,}\!\sum_{l=0}^{\lfloor\frac{1}{2}n\rfloor}\! (-1)^{l} 
\frac{\, (2l)!\cdot(2^{2l+1}-1)\cdot|S_1(n,2l+1)|\,}{2^{4l}} \\[1mm]
&& \displaystyle\notag\!\!\!\!\!\!\!\!\!\!\!\!\!\!\!\!\!\!\!\!\!\!\!\!\!\!
=\,-\frac{1+\ln\pi}{\pi} +\frac{1}{\,2\,}\ln2\pi  - \frac{1}{\,4\pi\,}\!
\left\{1+\frac{1}{4}+\frac{1}{16}+\frac{1}{128}-\frac{119}{19\,200}
-\frac{71}{9216}-\frac{7853}{1\,290\,240} -\frac{12\,611}{2\,949\,120}  -\ldots\right\}
\end{eqnarray}
Furthermore, both series expansions \eqref{kj20ejcn2dnd} and \eqref{jkhc92h8cn23d}, used together with 
the reflection formula and the recurrence relationship for the $\Gamma$--function, yield
series with rational coefficients for any values of the form
$\ln\Gamma(\frac{1}{2}n\pm\alpha\pi^{-1})$,
where $n$ is integer.

As a final remark, we note that expression \eqref{jkhc92h8cn23d}, written for $z$ instead of $\frac{1}{2}+z$,
straightforwardly produces
another series expansion for the logarithm of the $\Gamma$--function
\begin{eqnarray}
\displaystyle\notag 
\ln\Gamma(z) = &&\displaystyle\left(z-\frac{1}{\,2\,}\right)\!\ln\!{\left(z-\frac{1}{\,2\,}\right)}
- z +\frac{1}{\,2\,}+\frac{1}{\,2\,}\ln2\pi  \\
\displaystyle&&\displaystyle\label{lj023od230dend}
- \frac{1}{\,\pi\,}\!
\sum_{n=1}^\infty\!\frac{ 1}{\,n\cdot n!\,}\!\sum_{l=0}^{\lfloor\frac{1}{2}n\rfloor}\! (-1)^{l} 
\frac{\, (2l)!\cdot(2^{2l+1}-1)\cdot|S_1(n,2l+1)|\,}{(4\pi)^{2l+1}\cdot\big(z-\frac{1}{2}\big)^{2l+1}}
\end{eqnarray}
which converges in the green zone given in Fig.~\ref{kawfe3o84yg} shifted by $\frac{1}{2}$ to the right. 
In particular, if $z$ is real, it converges for any $z>\frac{2}{3}$.\\

{\small\noindent\textbf{Remark}
Expansion \eqref{jkhc92h8cn23d} may be also derived if we replace in \eqref{pio2u092je3n4} 
Binet's formula by its analog with ``conjugated'' denominator
\begin{equation}\notag
\int\limits_0^{\,\infty} \!\frac{\,\arctg{a x}\,}{e^{b x}+1}\, dx \,=\,
-\frac{\,\pi\,}{b} \ln\Gamma\!\left(\!\dfrac{1}{2}+\dfrac{b}{\,2\pi a\,}\!\right) 
-\frac{1}{\,2a\,}\!\left(1+ \ln\frac{2\pi a}{b} \right) + \frac{\,\pi\,}{\,2b\,}\ln2\pi .
\end{equation}
where $a>0$ and $\Re{b}>0$, see \cite[Sect.~4, exercise \no 40-a]{iaroslav_06}, or if we replace it by the following formula
\begin{equation}\notag
\int\limits_0^{\,\infty}\!\frac{\,\arctg{ax}\,}{\sh{bx}}\, dx \,=\, \frac{\,\pi\,}{b}\left\{
\ln\Gamma\!\left(\!\dfrac{b}{\,2\pi a\,}\!\right)
- \ln\Gamma\!\left(\!\dfrac{1}{2}+\dfrac{b}{\,2\pi a\,}\!\right) 
-\frac{1}{2}\ln\frac{2\pi a}{b}\right\} 
\end{equation}
derived in \cite[Sect.~4, exercise \no 39-e]{iaroslav_06}.
Making a change of 
variable $\,x=-\frac{2}{b}\arcth u\,$, and then, proceeding analogously to \eqref{k034fjdfnsh1}--\eqref{kj20ejcn2dnd}, yields 
\be
\begin{array}{l}
\displaystyle
\ln\Gamma(z) - \ln\Gamma\!\left(\frac{1}{2}+z\right) = \,-\frac{1}{2}\ln{z} + \frac{1}{\pi}
\sum_{n=0}^\infty \frac{1}{\,2n+1\,}\!
\sum_{k=1}^{2n+1}\frac{2^k}{k!} \binom{2n}{k-1} \sum_{l=0}^{n}
(-1)^l \frac{(2l)!\cdot S_1(k,2l+1)}{(2\pi z)^{2l+1}}
\end{array}
\ee
which, being combined with \eqref{jkhc92h8cn23d}, leads to a rearranged version of \eqref{kj20ejcn2dnd}.}\normalsize

\subsection{Series expansion for the polygamma functions} 
By differentiating expressions \eqref{kj20ejcn2dnd} and \eqref{lj023od230dend}, one may easily deduce similar series expansions for the polygamma functions. 
Differentiating the former expansion yields the following series representations
for the digamma 
and trigamma 
functions
\be\label{odi2je09cnhg}
\begin{array}{ll}
\displaystyle
\Psi(z)\; & \displaystyle=\,\ln z  - \frac{1}{\,2z\,} - \frac{1}{\,\pi z\,}\!
\sum_{n=1}^\infty\!\frac{ 1}{\,n\cdot n!\,}\!\sum_{l=0}^{\lfloor\frac{1}{2}n\rfloor}\! (-1)^{l} 
\frac{\, (2l+1)!\cdot|S_1(n,2l+1)|\,}{(2\pi z)^{2l+1}} \\[5mm]
& \displaystyle 
= \,\ln z  - \frac{1}{\,2z\,} - \frac{1}{\,\pi z\,}
\left\{\frac{1}{2\pi z}+\frac{1}{8\pi z}+\frac{1}{18}\left(\frac{1}{\pi z}-\frac{3}{4\pi^3 z^3}\right)
+\frac{3}{96}\left(\frac{1}{\pi z}-\frac{3}{2\pi^3 z^3}\right)\right.\\[5mm]
& \displaystyle \quad\left.
+\frac{1}{600}\left(\frac{12}{\pi z}-\frac{105}{4 \pi^3 z^3}+\frac{15}{4 \pi^5 z^5}\right)
+\frac{1}{4320}\left(\frac{60}{\pi z}-\frac{675}{4 \pi^3 z^3}+\frac{225}{4 \pi^5 z^5}\right)+\ldots\right\}
\end{array}
\ee
and 
\be\label{ojhfd83fh4hfhg}
\begin{array}{ll}
\displaystyle
\Psi_1(z)\; & \displaystyle=\,\frac{1}{\,2z^2\,} + \frac{1}{\,z\,} +  \frac{1}{\,\pi z^2\,}\!
\sum_{n=1}^\infty\!\frac{ 1}{\,n\cdot n!\,}\!\sum_{l=0}^{\lfloor\frac{1}{2}n\rfloor}\! (-1)^{l} 
\frac{\, (2l+2)!\cdot|S_1(n,2l+1)|\,}{(2\pi z)^{2l+1}} \\[5mm]
& \displaystyle 
= \,\frac{1}{\,2z^2\,} + \frac{1}{\,z\,} +  \frac{1}{\,\pi z^2\,}
\left\{\frac{1}{\pi z}+\frac{1}{4\pi z}+\frac{1}{18}\left(\frac{2}{\pi z}-\frac{3}{\pi^3 z^3}\right)
+\frac{3}{96}\left(\frac{2}{\pi z}-\frac{6}{\pi^3 z^3}\right)\right.\\[5mm]
& \displaystyle \quad\left.
+\frac{1}{600}\left(\frac{24}{\pi z}-\frac{105}{\pi^3 z^3}+\frac{45}{2\pi^5 z^5}\right)
+\frac{1}{4320}\left(\frac{120}{\pi z}-\frac{675}{\pi^3 z^3}+\frac{675}{2\pi^5 z^5}\right)+\ldots\right\}
\end{array}
\ee
respectively.
More generally, by differentiating $k$ times with respect to $z$ the above series for $\Psi(z)$, we obtain a series expansion
for the $k$th polygamma function
\be\label{094udn0sd4rtf}
\begin{array}{ll}
\displaystyle
\Psi_k(z)\; & \displaystyle=\,(-1)^{k+1}\frac{k!}{\,2z^{k+1}\,} + (-1)^{k+1}\frac{(k-1)!}{\,z^{k}\,} +  \frac{(-1)^{k+1}}{\,\pi z^{k+1}\,}\!
\sum_{n=1}^\infty\!\frac{ 1}{\,n\cdot n!\,}\!\sum_{l=0}^{\lfloor\frac{1}{2}n\rfloor}\! (-1)^{l} 
\frac{\, (2l+k+1)!\cdot|S_1(n,2l+1)|\,}{(2\pi z)^{2l+1}} \\[5mm]
& \displaystyle 
= \,(-1)^{k+1}\frac{k!}{\,2z^{k+1}\,} + (-1)^{k+1}\frac{(k-1)!}{\,z^{k}\,} +  \frac{(-1)^{k+1}}{\,\pi z^{k+1}\,}\!
\left\{\frac{(k+1)!}{2\pi z}+\frac{(k+1)!}{8\pi z}+\frac{1}{18}\left[\frac{(k+1)!}{\pi z}-\frac{(k+3)!}{8\pi^3 z^3}\right]
\right.\\[5mm]
& \displaystyle \quad\left.
+\frac{3}{96}\left[\frac{(k+1)!}{\pi z}-\frac{(k+3)!}{4\pi^3 z^3}\right]
+\frac{1}{600}\left[\frac{12\,(k+1)!}{\pi z}-\frac{35\,(k+3)!}{8\pi^3 z^3}+\frac{(k+5)!}{32\pi^5 z^5}\right]
+\ldots\right\}
\end{array}
\ee
where $k=1, 2, 3, \ldots$ 
Convergence analysis of these series is analogous to that performed in Section \ref{wkc2me23dwe}, and
we omit the details because the calculations are a little bit long. This analysis
reveals that the general term of these series may be always bounded by $\alpha_k(z) n^{-2}$, where $\alpha_k(z)$
depends solely on $z$ and on the order of the polygamma function $k$. At large $n$, the general term of these series is of the same order as
$\big(n\ln^{m}\! n\big)^{-2}$, where $m=1$ for $\Psi(z)$, $m=2$ for $\Psi_1(z)$ and $\Psi_2(z)$, $m=3$ for $\Psi_3(z)$ and $\Psi_4(z)$,
and so on.

\begin{figure}[!t]   
\centering
\includegraphics[width=0.8\textwidth]{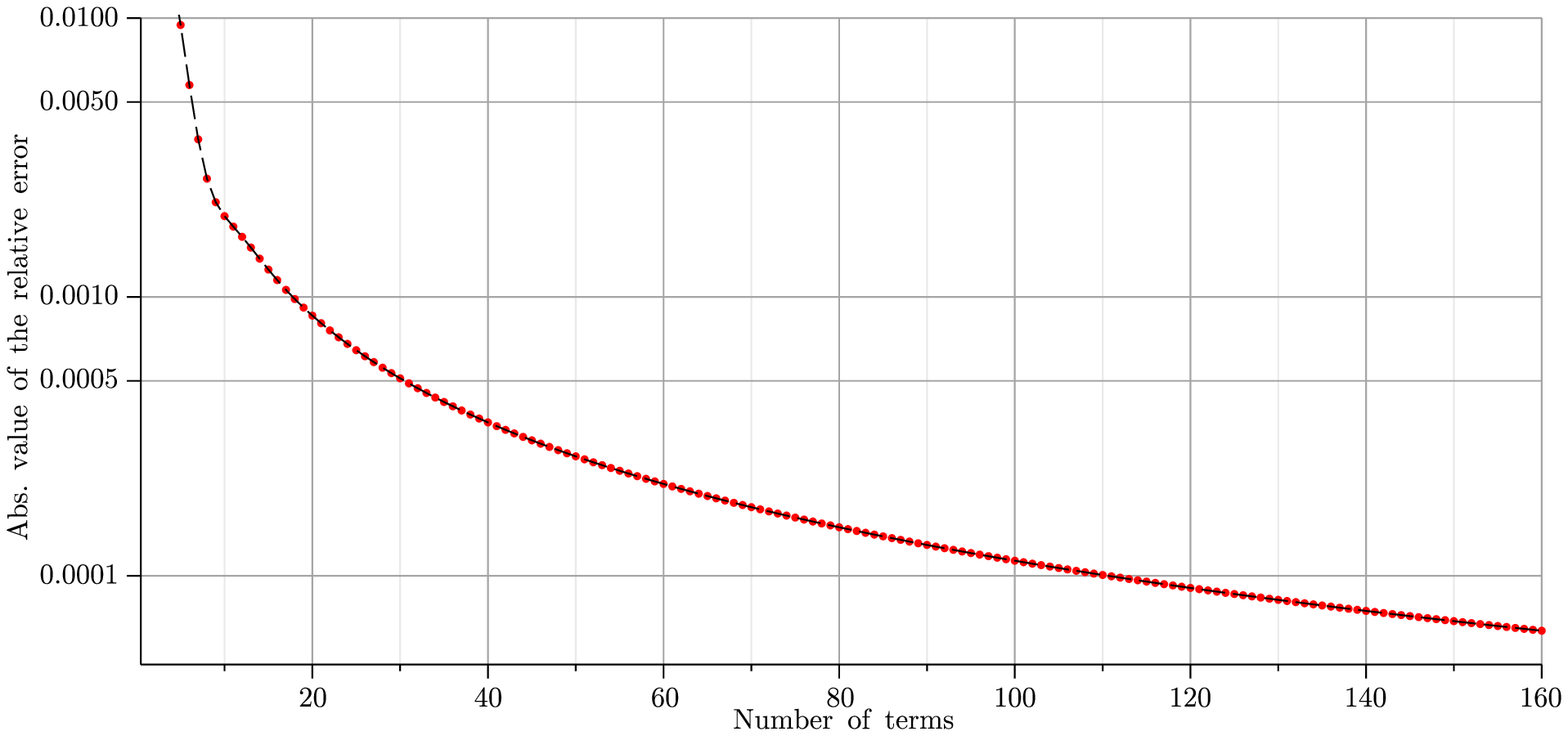}\\[-3mm]
\includegraphics[width=0.8\textwidth]{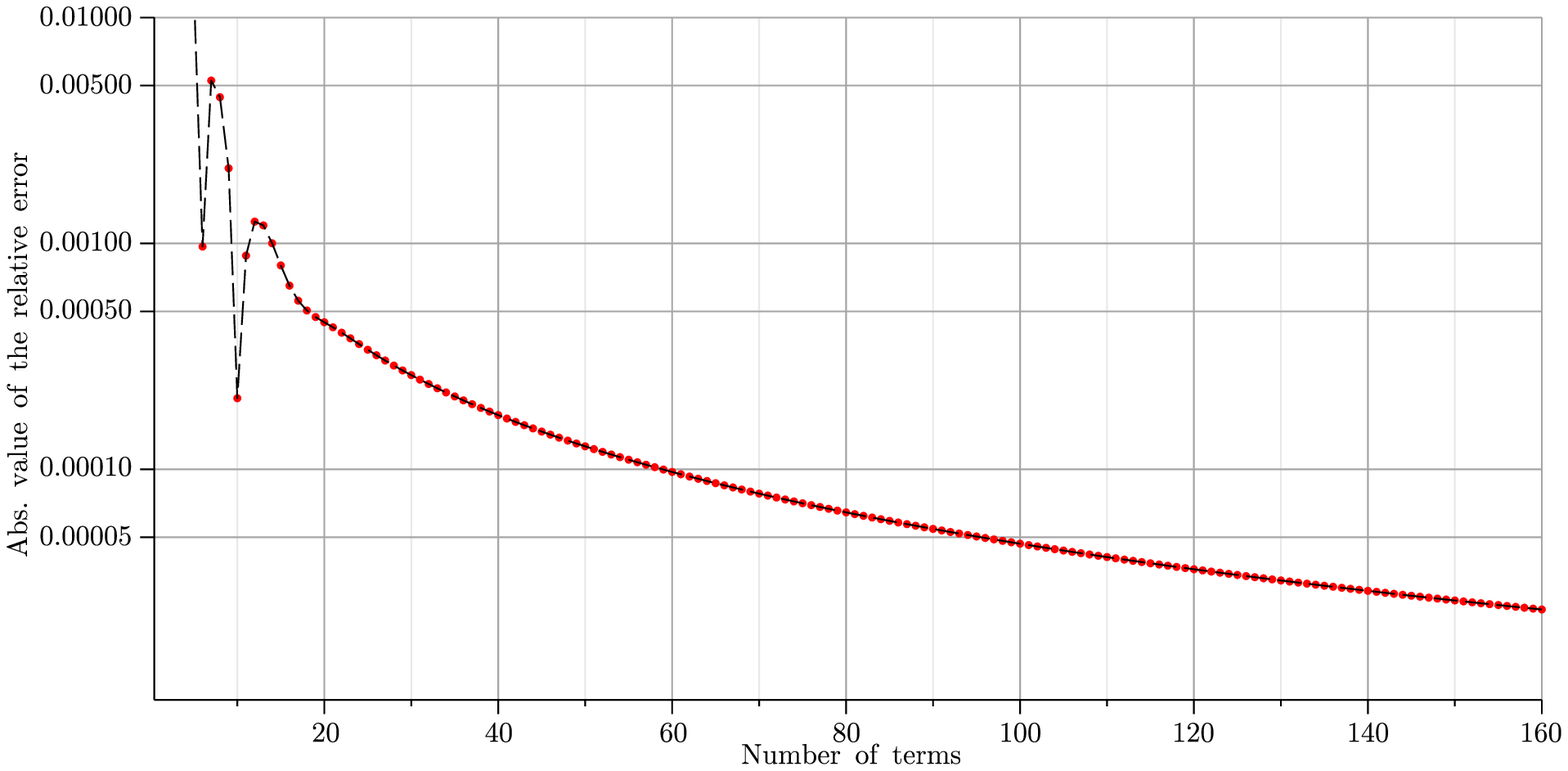}\\[-8mm]
\caption{\emph{Top:} Relative error of the series expansion for $\,\Psi(\pi^{-1})\,$ given by \protect\eqref{podi09123dmmx3}. 
\emph{Bottom:} Relative error of the series expansion for $\,\Psi_1(\pi^{-1})\,$ given by \protect\eqref{lkx23k0xe3c}. For better visibility, 
both errors are presented in absolute values and logarithmic scales. 
}
\label{hgd38cbda}
\end{figure}

We now give several particular cases of the above expansions. 
From \eqref{odi2je09cnhg}--\eqref{094udn0sd4rtf}, 
it follows that at $\pi^{-1}$, the polygamma functions have the following series representations
\be\label{podi09123dmmx3}
\begin{array}{ll}
\displaystyle
\Psi\!\left(\!\dfrac{1}{\,\pi \,}\!\right) \:& \displaystyle=\,-\ln \pi  - \frac{\pi}{\,2\,} - 
\sum_{n=1}^\infty\!\frac{ 1}{\,n\cdot n!\,}\!\sum_{l=0}^{\lfloor\frac{1}{2}n\rfloor}\! (-1)^{l} 
\frac{\, (2l+1)!\cdot|S_1(n,2l+1)|\,}{2^{2l+1}} \\[6mm]
& \displaystyle
=-\ln \pi  - \frac{\pi}{\,2\,} - \frac{1}{2}- \frac{1}{8}-\frac{1}{72}+\frac{1}{64}+\frac{7}{400}
+\frac{7}{576}+\frac{643}{94\,080}+\frac{103}{30\,720} + \ldots  
\end{array}
\ee
\be\label{lkx23k0xe3c}
\begin{array}{ll}
\displaystyle
\Psi_1\!\left(\!\dfrac{1}{\,\pi \,}\!\right) \:& \displaystyle=\,\frac{\pi^2}{\,2\,} + \pi +  \pi\!
\sum_{n=1}^\infty\!\frac{ 1}{\,n\cdot n!\,}\!\sum_{l=0}^{\lfloor\frac{1}{2}n\rfloor}\! (-1)^{l} 
\frac{\, (2l+2)!\cdot|S_1(n,2l+1)|\,}{2^{2l+1}} \\[6mm]
& \displaystyle
=\, \frac{\pi^2}{\,2\,} +\pi + \pi\left\{1+\frac{1}{4}-\frac{1}{18}-\frac{1}{8}-\frac{39}{400}-
\frac{29}{576}-\frac{353}{23\,520}+\frac{11}{3840}+\ldots  \right\} 
\end{array}
\ee
and, more generally, for $k=1, 2, 3, \ldots\,$
\be\label{lkx23k0xe3c}
\begin{array}{ll}
\displaystyle
\Psi_k\!\left(\!\dfrac{1}{\,\pi \,}\!\right) \:& \displaystyle=\,(-1)^{k+1}\pi^k\cdot\left\{\frac{\pi \, k!}{\,2\,} + (k-1)! +  
\sum_{n=1}^\infty\!\frac{ 1}{\,n\cdot n!\,}\!\sum_{l=0}^{\lfloor\frac{1}{2}n\rfloor} \!(-1)^{l} 
\frac{\, (2l+k+1)!\cdot|S_1(n,2l+1)|\,}{2^{2l+1}} \right\}\\[6mm]
& \displaystyle
=\, (-1)^{k+1}\pi^k\cdot\left\{\frac{\pi \, k!}{\,2\,} + (k-1)! + \frac{(1+k)!}{2}+ \frac{(1+k)!}{8}
+\frac{1}{18}\left[(k+1)!-\frac{(k+3)!}{8}\right]
\right.\\[6mm]
& \displaystyle\quad
\left.
+\frac{3}{96}\left[(k+1)!-\frac{(k+3)!}{4}\right]
+\frac{1}{600}\left[12\,(k+1)!-\frac{35\,(k+3)!}{8}+\frac{(k+5)!}{32}\right] +\ldots
\right\} 
\end{array}
\ee
Figure \ref{hgd38cbda} shows the rate of convergence of first two series. 

Second variant of the series expansions for the polygamma functions follows from \eqref{lj023od230dend}.
Differentiating the latter with respect to $z$ yields
\be\label{okjhe038nf3n}
\Psi(z) \, = \,\ln\!{\left(z-\frac{1}{\,2\,}\right)}
+ \frac{1}{\,\pi\,}\!
\sum_{n=1}^\infty\!\frac{ 1}{\,n\cdot n!\,}\!\sum_{l=0}^{\lfloor\frac{1}{2}n\rfloor}\! (-1)^{l} 
\frac{\, (2l+1)!\cdot(2^{2l+1}-1)\cdot|S_1(n,2l+1)|\,}{(4\pi)^{2l+1}\cdot\big(z-\frac{1}{2}\big)^{2l+2}}\vspace{-1.5mm}
\ee
and\vspace{-1.5mm}
\be\label{oif34jf3fnmnf}
\Psi_k(z) \,= \,\frac{\,(-1)^{k+1}(k-1)!\,}{\big(z-\frac{1}{2}\big)^k}
+ \frac{(-1)^k}{\,\pi\,}\!
\sum_{n=1}^\infty\!\frac{ 1}{\,n\cdot n!\,}\!\sum_{l=0}^{\lfloor\frac{1}{2}n\rfloor}\! (-1)^{l} 
\frac{\, (2l+k+1)!\cdot(2^{2l+1}-1)\cdot|S_1(n,2l+1)|\,}{(4\pi)^{2l+1}\cdot\big(z-\frac{1}{2}\big)^{2l+k+2}}
\ee
In particular,
\be
\begin{array}{ll}
\displaystyle
\Psi\!\left(\!\dfrac{1}{\,2\,}+\dfrac{1}{\,\pi \,}\!\right) \:& \displaystyle=\,-\ln \pi  +
\sum_{n=1}^\infty\!\frac{ 1}{\,n\cdot n!\,}\!\sum_{l=0}^{\lfloor\frac{1}{2}n\rfloor}\! (-1)^{l} 
\frac{\, (2l+1)!\cdot(2^{2l+1}-1)\cdot|S_1(n,2l+1)|\,}{2^{4l+2}} \\[6mm]
& \displaystyle
=-\ln \pi  + \frac{1}{4} + \frac{1}{16}-\frac{5}{576}-\frac{13}{512}-\frac{569}{25\,600}
-\frac{539}{36\,864} -\frac{98\,671}{12\,042\,240}-  \ldots
\end{array}
\ee
Similarly to expansions for $\ln\Gamma(z)$, expansions \eqref{094udn0sd4rtf} and \eqref{oif34jf3fnmnf}, combined with
the reflection formula and the recurrence relationship for polygamma functions,
give series with rational coefficients for any polygamma function of the argument
\mbox{$\,\frac{1}{2}n\pm\alpha\pi^{-1}\,$},
where $\alpha$ is rational greater than $\frac{1}{6}\pi$ and $n$ is integer.\\

\noindent\textbf{Final remark}\label{ckj3w0emnd}
Series which we discovered in the present work are very interesting especially because of the implication of combinatorial numbers $S_1(n,l)$.
In this context, it seems appropriate to note that series of a similar nature for $\ln\Gamma(z)$ and $\Psi_k(z)$ were already obtained, but remain
little--known and practically not mentioned in modern literature. 
For instance, in 1839, Jacques Binet \cite[pp.~231--236, 257, 237, 235, 335--339, Eqs.~(63), (81)]{binet_01} obtained several, rapidly convergent
for large $|z|$, expansions
\be
\begin{array}{l}
\displaystyle
\ln\Gamma(z) \,=\, \left(z-\frac{1}{\,2\,}\right)\!\ln z -z +\frac{1}{\,2\,}\ln2\pi 
+\frac{1}{2}\!\sum_{n=1}^\infty \frac{I(n)}{n}\cdot \frac{1}{ (z+1)_n} \\[5mm]
\displaystyle
\ln\Gamma(z) \,=\, \left(z-\frac{1}{\,2\,}\right)\!\ln z -z +\frac{1}{\,2\,}\ln2\pi 
-\frac{1}{2}\!\sum_{n=1}^\infty \frac{I'(n)}{n}\cdot \frac{1}{ (z)_n} \\[5mm]
\displaystyle
\Psi(z) \,=\, \ln z - \frac{1}{2z}  
-\frac{1}{2}\!\sum_{n=2}^\infty\!\frac{K(n)}{n}\cdot \frac{1}{ (z+1)_n} \\[5mm]
\displaystyle
\Psi(z) \,=\, \ln z - \frac{1}{2z}  
-\frac{1}{2}\!\sum_{n=2}^\infty\!\frac{K(n)-nK(n-1)}{n}\cdot \frac{1}{ (z)_n}
\end{array}
\ee
where numbers $I(n)$, $I'(n)$ and $K(n)$ are rational and may be given 
in terms of the Stirling numbers of the first kind. Binet recognized these numbers, denoted them by a capital $H$,
referenced Stirling's treatise \cite{stirling_01} and even corrected Striling's error: in the table on p.~11 \cite{stirling_01} 
he noticed that the value of $|S_1(9,3)|=118\,124$ and not $105\,056$. 
In our notations, Binet's formul\ae~for $I(n)$, $I'(n)$ and $K(n)$ read
\be
\begin{array}{l}
\displaystyle
I(n)\,= \int\limits_{0}^1 \!(2x-1)(x)_{n}\, dx 
=\sum_{l=1}^{n} \frac{l\cdot \big|S_1(n, l)\big|}{(l+1)(l+2)} \,= 
\,2\!\sum_{l=1}^{n} \frac{\big|S_1(n, l)\big|}{\,l+2\,} \,- \ \, C_{2,n} \\[5mm]
\displaystyle
I'(n)\,= \int\limits_{0}^1 \!(2x-1)(1-x)(x)_{n-1}\, dx 
=\sum_{l=1}^{n-2} \frac{l\cdot \big|S_1(n-1, l+1)\big|}{(l+2)(l+3)(l+4)} \\[5mm]
\displaystyle
K(n)\,= \, n!- 2\!\int\limits_{0}^1 \!(x)_{n}\, dx 
= \, n!- 2\!\sum_{l=1}^{n} \frac{\big|S_1(n, l)\big|}{\,l+1\,} \,=\,n!-2C_{2,n}
\end{array}
\ee
In particular, the first few coefficients are $I(1)=\frac{1}{6}$, $I(2)=\frac{1}{3}$, $I(3)=\frac{59}{60}$, 
$I(4)=\frac{58}{15}$, $I(5)=\frac{533}{28}$, $I(6)=\frac{1577}{14}$,
\ldots\,, $I'(1)=-\frac{1}{6}$, $I'(2)=0$, $I'(3)=+\frac{1}{60}$, $I'(4)=+\frac{1}{15}$, $I'(5)=+\frac{25}{84}$, $I'(6)=+\frac{11}{7}$,
\ldots\, and $K(2)=\frac{1}{3}$, $K(3)=\frac{3}{2}$, $K(4)=\frac{109}{15}$, 
$K(5)=\frac{245}{6}$, $K(6)=\frac{11\,153}{42}$, $K(7)=\frac{23\,681}{12}$, \ldots\footnote{Values $I'(1)=-\frac{1}{6}$, $I'(2)=0$ are found from the integral
formula, their definition via the sum with the Stirling numbers of the first kind being valid only for $n\geqslant3$.} Strictly speaking, Binet only found 
first four coefficients for each of these series and incorrectly calculated some of them (e.g.~for $I(4)$ he took $\frac{227}{60}$ instead 
of $\frac{232}{60}=\frac{58}{15}$, for $K(5)$ he took $\frac{245}{3}$ instead of $\frac{245}{6}$), 
but otherwise his method and derivations are correct.
Binet also remarked that
\be
\frac{I'(n)}{n}\cdot \frac{1}{ (z)_n} \,=\,O\left(\frac{1}{n^{z+1} \ln n}\right)\,,\qquad n\to\infty\,,
\ee
which implies that he knew the first--order approximation for the Cauchy numbers of the second kind as early as 1839.\footnote{Note, however,
that Binet stated this result without 
proof (he wrote \emph{Je ne développe pas ici ces résultats, parce que les détails sont un peu longs}).}
In 1923 Niels E.~Nørlund \cite[pp.~243--244]{norlund_02}, \cite[p.~335]{van_veen_01} obtained two series 
of a similar nature for the polygamma functions. In particular, for the Digamma function, he provided following expressions 
\be
\begin{array}{l}
\displaystyle
\Psi(z) \,=\, \ln z - \frac{1}{z}  
+\!\sum_{n=1}^\infty \frac{C_{2,n}}{n}\cdot \frac{1}{ (z+1)_n}   \\[5mm]
\displaystyle
\Psi(z) \,=\, \ln z - \frac{1}{2z}  
-\!\sum_{n=2}^\infty \frac{\big|G_n\big| \cdot(n-1)! }{ (z)_n} 
\end{array}
\ee
A careful inspection of both formul\ae~reveals that they actually are rewritten versions of the foregoing expansions for $\Psi(z)$ given by 
Binet 84 years earlier.\footnote{In order to reduce first Binet's series to first Nørlund's series, it suffices to remark that
$\sum\limits_{n=2}^\infty\! \frac{(n-1)!}{ (z+1)_n} 
= \sum\limits_{n=1}^\infty\! \frac{1\cdot2\cdots n}{ (z+1)(z+2)\cdots(z+n+1)}= \frac{1}{\,z\,} - \frac{1}{\,z+1\,}$ and
$C_{2,1}=\frac{1}{2}$. The equivalence between second Binet's series and second Nørlund's series follows from the fact
that $\frac{1}{2}\big[K(n) -nK(n-1)\big]=n C_{2,n-1}-C_{2,n}=\big|C_{1,n}\big|\equiv \big|G_n\big|\cdot n!\,$,
where $K(1)=0$ and $n=2, 3, 4,\ldots$} One may also notice that 
Fontana--Masceroni series \eqref{fms}, \eqref{fms2x}, is a particular case of the latter formula when $z=1$.
In contrast, the former expression at $z=1$ yields a not particularly well--known series for Euler's constant
\be\label{c43f13ffv1qe}
\gamma\,=\,1 -  \sum_{n=1}^\infty\!\frac{C_{2,n}}{\,n\cdot(n+1)!\,} =\,
1-\frac{1}{4}-\frac{5}{72}-\frac{1}{32}-\frac{251}{14\,400}-\frac{19}{1728} -
\frac{19\,087}{2\,540\,160} - \ldots
\ee
which is, in fact, closely related to the above--mentioned Fontana--Masceroni series and may be 
reduced to the latter by means of the recurrence relation
\be
n C_{2,n-1}-C_{2,n}=\big|C_{1,n}\big|\equiv \big|G_n\big|\cdot n!\,,\qquad  C_{2,0}=1\,,\quad
n=1, 2, 3,\ldots
\ee
Namely, by partial fraction decomposition, \eqref{c43f13ffv1qe} becomes
\be\notag
\begin{array}{l}
\displaystyle
\gamma\,=\,1 -  \sum_{n=1}^\infty\!\frac{C_{2,n}}{\,n!\,}\cdot\frac{1}{\,n(n+1)\,} \,=\,
1 -  \sum_{n=1}^\infty\!\frac{C_{2,n}}{\,n\cdot n!\,} + \sum_{n=1}^\infty\!\frac{C_{2,n}}{\,(n+1)!\,} \,=\\[6mm]
\displaystyle \phantom{\gamma\,} 
=\,1 - \sum_{n=1}^\infty\!\frac{C_{2,n-1} - \big|G_n\big|\cdot (n-1)!}{\,n!\,} 
+ \sum_{n=1}^\infty\!\frac{C_{2,n}}{\,(n+1)!\,}
\,=\sum_{n=1}^\infty\!\frac{\,\big|G_n\big|\,}{\,n\,}
\end{array}
\ee
It is interesting to note that series \eqref{c43f13ffv1qe} converges
at the same rate as $\sum n^{-2}\ln^{-1} n$, while Fontana--Masceroni series \eqref{fms}, \eqref{fms2x}
converges slightly faster, as $\sum n^{-2}\ln^{-2} n$, see \eqref{j38ndbr893r}.
It seems also appropriate to note here, that apart from Nørlund, the series expansions equivalent or similar to those derived by Binet in 1839, were 
also obtained (sometimes simply rediscovered, sometimes generalized) by various contemporaneous 
writers, see e.g.~\cite[p.~2052, Eq.~(1.17)]{coffey_08}, \cite[p.~11]{coffey_07}, \cite{nemes_03}, \cite[pp.~4005--4007]{young_03}.

\section*{Acknowledgments} 
The author is grateful to Oliver Oloa for providing several useful references regarding Stirling numbers,
to Vladimir V.~Kruchinin for sending electronic versions of \cite{kruchinin_02} and \cite{kruchinin_03} and to Nico Temme
for sending a scanned version of \cite{van_veen_01}.
The author also would like to thank Victor Kowalenko, Stefan Kr\"amer, Gergő Nemes,
Larry Glasser and an anonymous reviewer for their remarks and suggestions. Finally, the author is grateful to 
Vladimir V.~Reshetnikov for the numerical verification of some of the results obtained in this work and for correcting several misprints 
in the draft version of this paper.

\section*{A note on the history of this article}
Various internet searches may indicate that this article was first expected to be published
by the journal ``Mathematics of Computation''
(article submitted 
on 18 August 2014 and accepted for publication on 3 December 2014).\footnote{A more complete 
description of the publication history may be traced by consulting the arXiv version of this paper arXiv:1408.3902}
However, due to a 
disagreement with the managing editor of this journal during the production of this paper, I decided to withdraw it. 

Note also that the present article
was written before the recently published paper \cite{iaroslav_09}, which is an extension of the
present work to generalized Euler's constants (Stieltjes constants).

\appendix
\section{Some series expansions involving Stirling numbers, Gregory's coefficients, Cauchy numbers,  
Bernoulli numbers, binomial coefficients and harmonic numbers}\label{app_1}
In this first supplementary part of our work, we give the proofs of some of the results
given in Section \ref{lkce2mcwndc}. 
First of all, we show that Fontana--Mascheroni's formula \eqref{fms} may be readily derived from the 
generating equations for the Stirling numbers. The proof is similar to \eqref{kljd023jdnr2}:
\be\label{jkd293ndnw}
\begin{array}{c}
\displaystyle
\sum_{n=1}^{\infty}\!\frac{\, (-1)^{n-1}\,}{n}\cdot\frac{1}{\,n!\,}\!\sum_{l = 1 }^{n} \!\frac{\,S_1(n, l)\,}{l+1}\,=\,
-\sum_{n=1}^{\infty}\!\frac{1}{\,n!\,}\cdot\!\sum_{l = 1 }^{n} \!\frac{\,(-1)^l\big|S_1(n, l)\big|\,}{l+1}
\!\int\limits_0^1 \! x^{n-1} dx \\[8mm]
\displaystyle
\,=\,\sum_{l=1}^{\infty}\!\frac{ (-1)^{l+1}}{\,l+1\,}\!\int\limits_0^1 \!\underbrace{\left[
\sum_{n = 1}^{\infty} \!\frac{\,\big|S_1(n, l)\big|\,}{n!} x^n\right] }_{\text{see \eqref{h37qp2b237ds}}}\!\frac{dx}{x} 
\,=\,-\sum_{l=1}^{\infty}\!\frac{ 1}{\,(l+1)!\,}\!\underbrace{\int\limits_0^1 \!\frac{\ln^{l}(1-x)}{x}\,dx}_{\text{see \eqref{kljd023jdnr2}}} \\[8mm]
\displaystyle
\,=\,\sum_{l=1}^{\infty}\!\frac{(-1)^{l+1}\!\cdot\zeta(l+1)}{\,l+1\,}\,=\,\gamma
\end{array}
\ee
in virtue of known Euler's representation for $\gamma$, see e.g.~\cite[vol.~I, p.~45, Eq.~1.17(3)]{bateman_01}
or put simply $z=1$ into Legendre's series \eqref{098u43hfver}. The above derivation may be readily generalized to
\be\notag
\sum_{n=1}^{\infty}\!\frac{\, (-1)^{n-1}\,}{\,n\cdot n!\,}\!\sum_{l = 1 }^{n}  f(l)\,S_1(n, l)
\,=\,\sum_{l=1}^{\infty} (-1)^{l+1}f(l)\zeta(l+1)\,,
\ee
where $f(l)$ is an arbitrary function providing the convergence. In particular, taking $f(l)=1/(l+k)$,
the above equation reduces to
\be\label{d23dvygvs}
\sum_{n=1}^{\infty}\!\frac{\, (-1)^{n-1}\,}{\,n\cdot n!\,}\!\sum_{l = 1 }^{n} \!\frac{\,S_1(n, l)\,}{l+k}
\,=\,\sum_{l=2}^{\infty}\!\frac{(-1)^{l}\!\cdot\zeta(l)}{\,l+k-1\,}\,, \qquad\qquad k=1, 2, 3,\ldots
\ee
the series on the left converging as $\,\sum n^{-2} \ln^{-k-1}\! n \,$, see \eqref{kf3094fmf}. As concerns the series
in the right--hand side, it can be always expressed in terms of elementary functions,
Euler's constant $\gamma$, the $\zeta$--function and its first--order derivatives. 
For instance, integrating Legendre's series \eqref{098u43hfver}
over $z\in[0,1]$ and using Raabe's formula, see e.g.~\cite[vol.~I, p.~24, Eq.~1.9.1(18)]{bateman_01}, \cite[p.~552, Eq.~(30)]{iaroslav_07}, we get
\be
\sum_{l=2}^{\infty}\!\frac{(-1)^{l}\!\cdot\zeta(l)}{\,l(l+1)\,}\,=\,\frac{\,\gamma+\ln2\pi\,}{2} -1 
\ee
By partial fraction decomposition and using again above Euler's representation for $\gamma$, we find
\be\label{ij0e49ufnerndfh}
\sum_{l=2}^{\infty}\!\frac{(-1)^{l}\!\cdot\zeta(l)}{\,l+1\,}\,=\,\frac{\,\gamma-\ln2\pi\,}{2} +1 
\ee
so that for $k=2$
\be
\sum_{n=1}^{\infty}\!\frac{\, (-1)^{n-1}\,}{\,n\cdot n!\,}\!\sum_{l = 1 }^{n} \!\frac{\,S_1(n, l)\,}{l+2}
\,=\,\frac{\,\gamma-\ln2\pi\,}{2} +1 
\ee
which converges as $\,\sum n^{-2} \ln^{-3}\! n \,$, see \eqref{klxj2103xnd}.
Proceeding analogously and replacing Raabe's formula by the $(k-2)$th--order moment of $\ln\Gamma(z)$, 
see \cite[p.~175, (6.14)]{espinosa_01},\footnote{Expression for the first--order moment may be also found in
\cite[p.~552]{iaroslav_07}.} 
we obtain for $k=2, 3, 4,\ldots$
\be\label{c3e94mnd1}
\begin{array}{ll}
\displaystyle
\sum_{n=1}^{\infty}\!\frac{\, (-1)^{n-1}\,}{\,n\cdot n!\,}&\displaystyle\sum_{l = 1 }^{n} \!\frac{\,S_1(n, l)\,}{l+k}
\,=\sum_{l=2}^{\infty}\!\frac{(-1)^{l}\!\cdot\zeta(l)}{\,l+k-1\,}\,=\,\frac{\,1\,}{\,k-1\,} -\frac{\,\ln2\pi\,}{k}+\frac{\,\gamma\,}{2} 
\\[6mm]
&\displaystyle\quad
+\!\!\!\!\sum_{l=1}^{\lfloor\frac{1}{2}(k-1)\rfloor}\!\!\!\! (-1)^l \binom{k-1}{2l-1}\frac{\,(2l)!\cdot\zeta'(2l)\,}{l\cdot(2\pi)^{2l}}
+\!\!\sum_{l=1}^{\lfloor\frac{1}{2}k\rfloor-1}\!\! (-1)^l \binom{k-1}{2l}\frac{\,(2l)!\cdot\zeta(2l+1)\,}{2\cdot(2\pi)^{2l}}
\qquad
\end{array}
\ee
In particular, for $k=3, 4$ and $5$, the above formula yields following series
\be\label{c3e94mnd}
\begin{array}{l}
\displaystyle
\sum_{n=1}^{\infty}\!\frac{\, (-1)^{n-1}\,}{\,n\cdot n!\,}\!\sum_{l = 1 }^{n} \!\frac{\,S_1(n, l)\,}{l+3}
\,=\sum_{l=2}^{\infty}\!\frac{(-1)^{l}\!\cdot\zeta(l)}{\,l+2\,}\,=\,\frac{\,1\,}{2} -\frac{\,\ln2\pi\,}{3}+\frac{\,\gamma\,}{2} 
-\frac{\,\zeta'(2)\,}{\pi^2}  \\[6mm]
\displaystyle
\sum_{n=1}^{\infty}\!\frac{\, (-1)^{n-1}\,}{\,n\cdot n!\,}\!\sum_{l = 1 }^{n} \!\frac{\,S_1(n, l)\,}{l+4}
\,=\sum_{l=2}^{\infty}\!\frac{(-1)^{l}\!\cdot\zeta(l)}{\,l+3\,}\,=\,\frac{\,1\,}{3}  -\frac{\,\ln2\pi\,}{4}+\frac{\,\gamma\,}{2}
-\frac{\,3\zeta'(2)\,}{2\pi^2} -\frac{\,3\zeta(3)\,}{4\pi^2}  \\[6mm]
\displaystyle
\sum_{n=1}^{\infty}\!\frac{\, (-1)^{n-1}\,}{\,n\cdot n!\,}\!\sum_{l = 1 }^{n} \!\frac{\,S_1(n, l)\,}{l+5}
\,=\sum_{l=2}^{\infty}\!\frac{(-1)^{l}\!\cdot\zeta(l)}{\,l+4\,}\,=\,\frac{\,1\,}{4} -\frac{\,\ln2\pi\,}{5}+\frac{\,\gamma\,}{2} 
-\frac{\,2\zeta'(2)\,}{\pi^2} -\frac{\,3\zeta(3)\,}{2\pi^2} + \frac{\,3\zeta'(4)\,}{\pi^4} 
\end{array}
\ee
which converge at the same rate as $\,\sum n^{-2} \ln^{-4}\! n \,$, 
$\,\sum n^{-2} \ln^{-5}\! n \,$ and $\,\sum n^{-2} \ln^{-6}\! n \,$ respectively, see \eqref{kf3094fmf}.  
Note that numerically,
series with Stirling numbers
converge more rapidly than their analogs with the $\zeta$--function, because $\zeta(l)\sim1$ for $l\to\infty$.
It is interesting that if in \eqref{jkd293ndnw} we replace in the denominator $n$ by $n+1$,
and more generally by $n+k+1$, the resulting series reduces to elementary functions
\begin{eqnarray}
\displaystyle\notag
\sum_{n=1}^{\infty}\!\frac{\, (-1)^{n+1}\,}{\,n+k+1\,}\cdot\frac{1}{\,n!\,}\sum_{l = 1 }^{n} \frac{\,S_1(n, l)\,}{l+1}\;\,
&& =\,-\!\sum_{l=1}^{\infty} \frac{1}{\,(l+1)!\,}\! \int\limits_0^1 \!\! x^{k} \ln^{l}(1-x)\,  dx \\[2mm]
\displaystyle\label{89h98h78y6}
&&=\sum_{m=0}^k (-1)^m \binom{k}{m}\left\{\frac{1}{m+1}-\ln\frac{m+2}{m+1} \right\} \\[2mm]
\displaystyle\notag
&&=\,\frac{1}{k+1}+\sum_{m=1}^{k+1} (-1)^m \binom{k+1}{m}\ln(m+1)
\,,\qquad k=0, 1, 2,\ldots
\end{eqnarray}
since 
\be\notag
\sum\limits_{m=0}^k \frac{(-1)^m}{m+1} \binom{k}{m}\,=\,\frac{1}{k+1}\,,
\qquad\text{and}
\qquad\binom{k}{m}+\binom{k}{m-1}\,=\,\binom{k+1}{m}\,.
\ee
see e.g~\cite[p.~300, \no 30.12--30.13]{evgrafov_01_eng}.
In particular, for $k=1, 2$ and $3$, we obtain following series
\begin{eqnarray}
&&\displaystyle\sum_{n=1}^{\infty}\!\frac{\, (-1)^{n+1}\,}{\,n+1\,}\cdot\frac{1}{\,n!\,}\sum_{l = 1 }^{n} \frac{\,S_1(n, l)\,}{l+1}\,
=\, 1- \ln2   \label{hc2394ch23js}\\[1mm]
&&\displaystyle\sum_{n=1}^{\infty}\!\frac{\, (-1)^{n+1}\,}{\,n+2\,}\cdot\frac{1}{\,n!\,}\sum_{l = 1 }^{n} \frac{\,S_1(n, l)\,}{l+1}\,
=\, \frac{1}{2}-2\ln2 +\ln3  \\[1mm]
&&\displaystyle\sum_{n=1}^{\infty}\!\frac{\, (-1)^{n+1}\,}{\,n+3\,}\cdot\frac{1}{\,n!\,}\sum_{l = 1 }^{n} \frac{\,S_1(n, l)\,}{l+1}\,
=\, \frac{1}{3}-5\ln2+3\ln3 
\end{eqnarray}
first of which may be also found in \cite[Chapt.~V, p.~277]{jordan_01}, \cite[p.~426, Eq.~(50)]{kowalenko_01}, \cite[p.~129, Eq.~(7.3.14)]{alabdulmohsin_01}, \cite{skramer_01}.
On the contrary, if in Fontana--Mascheroni's series \eqref{jkd293ndnw} we replace $n$ by $n-1, n-2, \ldots\,$, and start the summation from 
an adequate value of $n$, then we again arrive at special constants
\begin{eqnarray}
&&\displaystyle\sum_{n=2}^{\infty}\!\frac{\, (-1)^{n+1}\,}{\,n-1\,}\cdot\frac{1}{\,n!\,}\sum_{l = 1 }^{n} \frac{\,S_1(n, l)\,}{l+1}
\,=\,  -\frac{1}{2} + \frac{\ln2\pi}{2} -\frac{\gamma}{2} \label{uy98yuhpuh78a}\\[1mm]
&&\displaystyle\sum_{n=3}^{\infty}\!\frac{\, (-1)^{n+1}\,}{\,n-2\,}\cdot\frac{1}{\,n!\,}\sum_{l = 1 }^{n} \frac{\,S_1(n, l)\,}{l+1}
\,=\,  -\frac{1}{8} + \frac{\ln2\pi}{12} - \frac{\zeta'(2)}{\,2\pi^2} \label{uy98yuhpuh78b}\\[1mm]
&&\displaystyle\sum_{n=4}^{\infty}\!\frac{\, (-1)^{n+1}\,}{\,n-3\,}\cdot\frac{1}{\,n!\,}\sum_{l = 1 }^{n} \frac{\,S_1(n, l)\,}{l+1}
\,=\,  -\frac{1}{16} + \frac{\ln2\pi}{24} - \frac{\zeta'(2)}{\,4\pi^2} + \frac{\zeta(3)}{\,8\pi^2}  \label{uy98yuhpuh78c}
\end{eqnarray}
the second series providing a comparatively simple representation for $\zeta'(2)/\pi^2$ in terms of rational coefficients.
The proof of these two formul\ae~is carried out in the same way as \eqref{jkd293ndnw}, with the additional help of \eqref{2iedu20j3hg8},  
\eqref{ij0e49ufnerndfh}, \eqref{c3e94mnd} [or \eqref{c3e94mnd1} for the more general case] and of these auxiliary formul\ae
\be\notag
\begin{array}{l}
\displaystyle
\sum_{n=2}^\infty \frac{\,(-1)^{n-1}S_1(n,1)\,}{(n-1)\cdot n!}\,=
\sum_{n=2}^\infty \frac{1}{\,n(n-1)\,}\,=\,1 \\[6mm]
\displaystyle
\sum_{n=3}^\infty \frac{\,(-1)^{n-1}S_1(n,1)\,}{(n-2)\cdot n!}\,=
\sum_{n=3}^\infty \frac{1}{\,n(n-2)\,}\,=\,\frac{3}{4}  \\[6mm]
\displaystyle
\sum_{n=3}^\infty \frac{\,(-1)^{n-1}S_1(n,2)\,}{(n-2)\cdot n!}\,=\,-\!
\sum_{n=3}^\infty \frac{H_{n-1}}{\,n(n-2)\,}\,=\,-\frac{3}{4}-\frac{\,\pi^2}{12} \\[6mm]
\displaystyle
\sum_{n=4}^\infty \frac{\,(-1)^{n-1}S_1(n,1)\,}{(n-3)\cdot n!}\,=
\sum_{n=4}^\infty \frac{1}{\,n(n-3)\,}\,=\,\frac{11}{18}  \\[6mm]
\displaystyle
\sum_{n=4}^\infty \frac{\,(-1)^{n-1}S_1(n,2)\,}{(n-3)\cdot n!}\,=\,-\!
\sum_{n=4}^\infty \frac{H_{n-1}}{\,n(n-3)\,}\,=\,-\frac{11}{12}-\frac{\,\pi^2}{18} \\[6mm]
\displaystyle
\sum_{n=4}^\infty \frac{\,(-1)^{n-1}S_1(n,3)\,}{(n-3)\cdot n!}\,=
\sum_{n=4}^\infty \frac{\,H^2_{n-1}-H^{(2)}_{n-1}\,}{\,2n(n-3)\,}\,=\,\frac{11}{36}+\frac{\,\pi^2}{12} + \frac{\,\zeta(3)\,}{3}
\end{array}
\ee
where $H_n$ stands for the $n$th harmonic number and $H^{(2)}_n$ denotes the $n$th generalized 
harmonic number of the second order.\footnote{Note that 
$\,S_1(n,1)=(-1)^{n-1} (n-1)!\,$, $\,S_1(n,2)=(-1)^{n} (n-1)!\cdot H_{n-1}$ and 
$\,S_1(n,3)=\frac{1}{2}(-1)^{n-1} (n-1)!\cdot \big\{H^2_{n-1}-H^{(2)}_{n-1}\big\}$,
see e.g.~\cite[p.~217]{comtet_01}, \cite[p.~1395]{shen_01}, \cite[p.~425, Eq.~(43)]{kowalenko_01}.} 
In virtue of asymptotics \eqref{ojc02m3d}, series 
\eqref{89h98h78y6}--\eqref{uy98yuhpuh78b} converge at the same rate as $\,\sum (n\ln n)^{-2}$.
More general formul\ae~of the same nature may be obtained with the aid of integral \eqref{o23id23j3}, which we come to evaluate later.
The reasoning similar to \eqref{jkd293ndnw} may be successfully applied to the evaluation of certain
series containing harmonic numbers in combination with 
Stirling numbers.\footnote{In this context, it may be useful to remark that $\int\limits_0^1 \! x^{n-1}\!\ln(1-x)\, dx \,=\, -{H_n}/{n}\,$
for $\,n=1, 2,3,\ldots$} For example
\be\label{cd230jd}
\sum_{n=1}^{\infty}\!\frac{\, (-1)^{n-1}H_n\,}{\,n\cdot n!\,}\!\sum_{l = 1 }^{n}  f(l)\,S_1(n, l)
\,=\,\sum_{l=1}^{\infty} (-1)^{l+1}(l+1)f(l)\zeta(l+2)\,,
\ee
where $f(l)$, as before, should be chosen so that the convergence is guaranteed. In particular, putting $f(l)=1/[(l+1)(l+2)]$ and using Euler's representation for $\gamma$
employed in the last line of \eqref{jkd293ndnw}
yields a relatively simple series with rational terms for Euler's constant
\be\label{kjl2039jds2}
\begin{array}{r}
\displaystyle
\gamma\,=\,\frac{\pi^2}{12}- \sum_{n=1}^{\infty}\!\frac{\, (-1)^{n-1}H_n\,}{\,n\cdot n!\,}\!\sum_{l = 1 }^{n}  
 \frac{\,S_1(n, l)\,}{\,(l+1)(l+2)\,}
\,=\,\frac{\pi^2}{12}
-\frac{1}{6}-\frac{1}{32}-{\frac {11}{810}}-{\frac {35}{4608}}\ldots \\[6mm]
\displaystyle
-{\frac {14\,659}{3\,024\,000}}- {\frac {1393}{414\,720}}
-{\frac {729\,751}{296\,352\,000}}
-{\frac {4\,368\,901}{2\,322\,432\,000}}-\ldots
\end{array}
\ee
In view of the fact that harmonic numbers grow logarithmically, see e.g.~\cite[p.~84, Eq.~(12)]{bromwich_02}, 
\cite[p.~46, \no 377]{gunter_03_eng}, 
and accounting for \eqref{ojc02m3d}, we establish that 
\be\notag
\frac{\, (-1)^{n-1}H_n\,}{\,n\cdot n!\,}\!\sum_{l = 1 }^{n}  
\frac{\,S_1(n, l)\,}{\,(l+1)(l+2)\,} \,\sim\,\frac{1}{\,n^2\ln n\,}\,,\qquad\qquad
 n\to\infty\,,
\ee
so that \eqref{kjl2039jds2} converges at the same rate as $\sum n^{-2} \ln^{-1}\! n$. Setting
$f(l)=1/[(l+1)(l+3)]$ into \eqref{cd230jd} yields yet another series for $\gamma$
\be
\begin{array}{r}
\displaystyle
\gamma\,=\,\ln2\pi +\frac{\pi^2}{9}-2 - 2\!\sum_{n=1}^{\infty}\!\frac{\, (-1)^{n-1}H_n\,}{\,n\cdot n!\,}\!\sum_{l = 1 }^{n}  
 \frac{\,S_1(n, l)\,}{\,(l+1)(l+3)\,}
\,=\,\ln2\pi +\frac{\pi^2}{9}-2  \\[6mm]
\displaystyle
- \frac{1}{4}-{\frac {7}{160}}-{\frac {121}{6480}}-{\frac 
{125}{12\,096}}-{\frac {39\,593}{6\,048\,000}}-{\frac {
140\,287}{31\,104\,000}}-{\frac {325\,127}{98\,784\,000}}-\ldots
\end{array}
\ee
converging at the same rate as \eqref{kjl2039jds2}. Analogously, setting $f(l)=1/(l+1)$ into \eqref{cd230jd}, one can show that\footnote{Note that the
corresponding series with the $\zeta$--function is semi--convergent.}
\be\label{98yh9gtv}
\sum_{n=1}^{\infty}\!\frac{\, (-1)^{n-1}H_n\,}{\,n\,}\!
\cdot\frac{1}{\,n!\,}\sum_{l = 1 }^{n} \frac{\,S_1(n, l)\,}{l+1}\,=\,\frac{\pi^2}{6}-1
\ee
which also converges as $\sum n^{-2} \ln^{-1}\! n$.

Finally, as we noticed in Section II, Stirling numbers of the first kind may also appear in the evaluation of certain integrals, 
which, at first sight, have nothing to do 
with Stirling numbers. For instance, it is known 
that some particular cases of $\int\ln^s(1-x) \, x^{-k} \,dx\,$ taken between $x=0$ and $x=1$
reduce to $\zeta$--function and to elementary functions.\footnote{Various particular cases of this integral appear in numerous texbooks,
see e.g. \cite[\no 2147]{gunter_02_eng}, }
By using Stirling numbers, we can show that more generally, if $k$ is integer and the integral is convergent, then the latter
may be always reduced to a finite set of $\zeta$--functions with coefficients containing Stirling numbers of the first kind and binomial coefficients. 
This
can be shown as follows. By differentiating $k$ times the well--known expansion $\,(1-y)^{-1}=\sum y^n$, 
we get
\be\notag
\frac{1}{(1-y)^{k+1}}\,=\,\frac{1}{k!}\sum_{n=k}^\infty n(n-1)\cdots(n-k+1) \, y^{n-k} =
\,\frac{1}{k!}\sum_{n=k}^\infty \sum_{r=1}^k S_1(k,r) \, n^r y^{n-k}\,,\qquad|y|<1\,,
\ee
in virtue of (\ref{x2l3dkkk03d}b). Rewriting this result for $y=e^{-t}$, multiplying it by $e^{-t}$ and putting $k-1$ instead of $k$, yields
\be\notag
\frac{e^{-t}}{\big(1-e^{-t}\big)^k}\,=\,\frac{1}{(k-1)!}\!\!
\sum_{n=k-1}^\infty \!\sum_{r=1}^{k-1} S_1(k-1,r) \, n^r e^{-(n-k+2)t}\,,\qquad
0<t<\infty
\ee
Consider now the integral in question. Making a change of variable $\,x=1-e^{-t}$, employing the above expansion and performing
the term--by--term integration, we have
\be\notag
\begin{array}{ll}
\displaystyle
\int\limits_0^1 \! \frac{\,\ln^s(1-x)\,}{x^k}\, dx \; &\displaystyle \,=(-1)^s\!\!\int\limits_0^\infty
\frac{\,t^s e^{-t}dt\,}{\big(1-e^{-t}\big)^k}  \,=\,\frac{(-1)^s}{(k-1)!}\!\!
\sum_{n=k-1}^\infty \!\sum_{r=1}^{k-1} S_1(k-1,r) \, n^r \!\!\!
\underbrace{\int\limits_0^\infty\! t^s e^{-(n-k+2)t} dt}_{\Gamma(s+1)\cdot(n-k+2)^{-s-1}} \\[6mm]
&\displaystyle 
=\,\frac{\,(-1)^s\!\cdot \Gamma(s+1)\,}{(k-1)!} \sum_{r=1}^{k-1} S_1(k-1,r) \!\!\!\sum_{n=k-1}^\infty \!\!
\frac{n^r}{(n-k+2)^{s+1}}
\end{array}
\ee
Changing the summation index $n=p+k-2$ and using the binomial expansion, the latter series becomes
\be\notag
\sum_{n=k-1}^\infty \!\!\frac{n^r}{(n-k+2)^{l+s}}\,=\,
\sum_{p=1}^\infty \frac{(p+k-2)^r}{p^{s+1}}\,=\,\sum_{m=0}^r \!
\binom{r}{m}(k-2)^{r-m}\zeta(s+1-m)
\ee
Whence we obtain
\be\label{o23id23j3}
\int\limits_0^1 \! \frac{\,\ln^s(1-x)\,}{x^k}\, dx  \,
=\,\frac{\,(-1)^s\!\cdot \Gamma(s+1)\,}{(k-1)!} \sum_{r=1}^{k-1} S_1(k-1,r) \! \sum_{m=0}^r \!
\binom{r}{m}\!\cdot(k-2)^{r-m}\cdot\zeta(s+1-m)
\ee
This formula holds for any for $k=3, 4, 5,\ldots\,$ and $\Re{s}>k-1$.
For lower values of $k$, this integral equals
\be\label{2iedu20j3hg8}
\int\limits_0^1 \! \frac{\,\ln^s(1-x)\,}{x^k}\, dx \,=\,
\begin{cases}
(-1)^s\Gamma(s+1)\zeta(s+1)  \,,\quad & k=1\,,\quad \Re{s}>0 \\[2mm]
 (-1)^s\Gamma(s+1)\zeta(s)   \,,\quad  &  k=2\,,\quad \Re{s}>1 \\[1.6mm]
\frac{1}{2} (-1)^s\Gamma(s+1)\big\{\zeta(s-1)+\zeta(s) \big\}\,,\quad  & k=3\,,\quad \Re{s}>2 \\[1.6mm]
\frac{1}{6} (-1)^s\Gamma(s+1)\big\{\zeta(s-2)+3\zeta(s-1)+2\zeta(s) \big\}\,,\quad  & k=4\,,\quad \Re{s}>3 
\end{cases}
\quad
\ee

\section{Bounds and full asymptotical expansions for some special numbers which appeared in Section II}\label{appendix}
\subsection{Bounds and full asymptotical expansion for the Cauchy numbers of the second kind $C_{2,n}$,
also known as generalized Bernoulli numbers $\big|B_n^{(n)}\big|\,$}\label{appendix1}
\begin{figure}[!t]   
\centering
\includegraphics[width=0.8\textwidth]{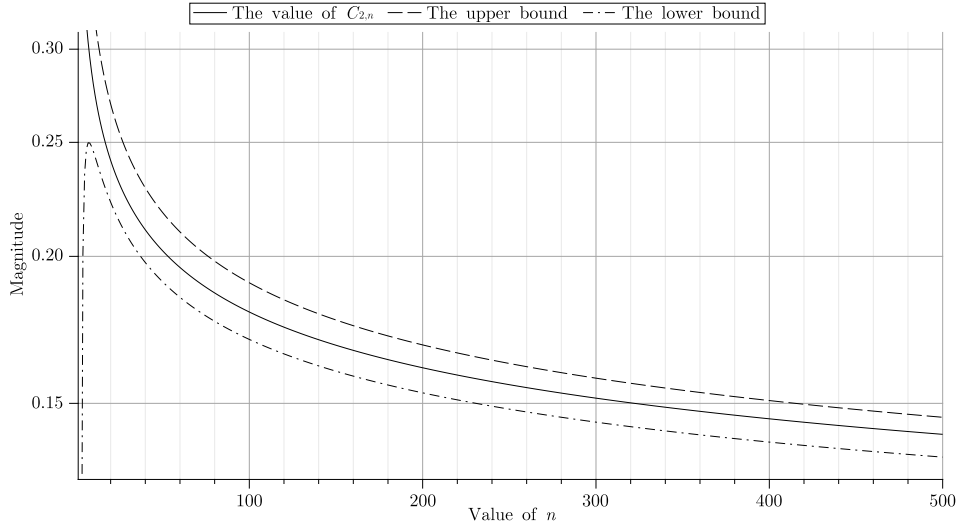}
\caption{Lower and upper bounds for the Cauchy numbers of the second kind $C_{2,n}$ given by inequalities \eqref{hu765f2dcefc}, 
logarithmic scale.}
\label{big8ghi87}
\end{figure}

Consider the Cauchy numbers of the second kind, which appear in series expansions for $\ln\ln2$, $\frac{z}{(1+z)\ln(1+z)}$ and
$\ln\ln(1+z)$, 
equations \eqref{fmse}, \eqref{i2039dj23r} and \eqref{ki1d039dn321} respectively.
Using definition \eqref{x2l3dkkk03d},
we may reduce it to a definite integral 
\be\label{hjc239fhsw}
\begin{array}{ll}
\displaystyle
C_{2,n} =\big|B_n^{(n)}\big|\;& \displaystyle=\,
\sum_{l=1}^{n} \frac{|S_1(n,l)|}{l+1} \,= \sum_{l=1}^{n} |S_1(n,l)|\!\int\limits_0^1  \! x^l \, dx\,=\\[6mm]
& \displaystyle
=\int\limits_0^1 \!\sum_{l=1}^{n} x^l  |S_1(n,l)| \, dx\, 
=\int\limits_0^1 \! (x)_n \, dx \,=\int\limits_0^1 \! \frac{\,\Gamma(x+n)\,}{\Gamma(x)} \, dx 
\end{array}
\ee
Since $\Gamma(x+n)$ for $x\in[0,1]$ and $n\geqslant2$ is positive and monotonically increases, the
following trivial inequalites are always true: $\,(n-1)!\leqslant\Gamma(x+n)\leqslant n!\,$ Therefore,
the Cauchy numbers of the second kind satisfy
\be\label{j87hf6fjh}
\frac{A}{\,n}\,\leqslant\,
\frac{C_{2,n}}{\,n!\,}
\,\leqslant \, A\,,\qquad
n=2, 3,4,\ldots
\ee
where
\be\notag
A\,\equiv\int\limits_0^1 \!\! \frac{\,dx\,}{\Gamma(x)} \,  =\,0.5412357343\ldots 
\ee
Numerical simulations indicate, however, that both bounds are very rough.
Better results may be obtained if we resort to more accurate estimations for the $\Gamma$--function; for instance, it is known that
\be\label{jd83dwsw}
(n+1)^{x-1}n!\leqslant\Gamma(x+n)\leqslant n^{x-1}n!\,\qquad\quad 
\begin{array}{l}
0\leqslant x\leqslant1\\[1mm] 
n=1, 2, 3,\ldots\,
\end{array}
\ee
where the left--hand side is strictly positive.\footnote{From the fact that the function 
$f(x,n)\,\equiv\,\frac{1}{1-x}\ln\!\left\{\Gamma(x+n) - \Gamma(n+1)\right\}\,$,
where $\,n=1, 2, 3,\ldots\,$ and $\,x\in[0,1]\,$,
is nonpositive and monotonically decreases, and because $\,f(0,n)\,=\,-\ln n\,$
and $\,\lim_{x\to1} f(x,n)\,=\,-\Psi(n+1) \,$,
it follows that $\,-\Psi(n+1)\leqslant f(x,n)\leqslant-\ln n\,$.
Therefore $\,e^{(x-1)\Psi(n+1)}n!\leqslant\Gamma(x+n)\leqslant n^{x-1}n!\,$ with the same conditions.
Since $\,\ln n<\Psi(n+1)<\ln(n+1)\,$, the latter also implies a weaker relation
$\,(n+1)^{x-1}n!\leqslant\Gamma(x+n)\leqslant n^{x-1}n!\,$ 
These inequalities are comparatively sharp, and albeit
they are quite elementary, they are usually attributed to Walter Gautschi who derived them in 1958 \cite[Eqs.~(6)--(7), Fig.~2]{gautschi_01}.}
Remarking that on the unit interval the function $1/\Gamma(x)$ is nonnegative and may be bounded
from below and from above as 
\be\label{d3rdf34s}
\,x\leqslant \frac{1}{\Gamma(x)}\leqslant (\gamma-1)x^2+(2-\gamma)x\,,\qquad 0\leqslant x\leqslant1\,,\quad\footnotemark\up{,}\footnotemark
\ee
as well as using \eqref{djoiejdnss}, we conclude that, on the one hand
\addtocounter{footnote}{-1}
\footnotetext{On the unit interval the function $1/\Gamma(x)$ is bounded and concave from $x=0$ to $x=0.3021417247\ldots\,$
(its first inflexion point on $\mathbbm{R}^+$)
and convex from $x=0.3021417247\ldots$ to $x=1$. It is therefore possible to construct a comparatively good
upper bound for  $1/\Gamma(x)$ on $x\in[0,1]$ with the aid of an arbitrary parabola 
$\,y(x)=a x^2+b x +c\,$ convex on $x\in[0,1]$ and whose coefficients may be found as follows: we first identify $\,c=0\,$ from the condition $\,y(0)=0\,$, then $\,b=1-a\,$ from the condition $\,y(1)=1\,$,
and finally $\,a=\gamma-1\,$ from the requirement that $\,\frac{1}{\Gamma(x)}-y(x)\,$ touches, without intersecting,
the coordinate axis at $x=1$ (i.e.~$\,\frac{1}{\Gamma(x)}-y(x)\,$ is minimum at $x=1$). The obtained upper bound is quite accurate: 
the maximum difference between the latter and $1/\Gamma(x)$ is $\epsilon\equiv0.05451361692\ldots$ and occurs at $x=0.2980329689\ldots$
(more precisely, the exact value of $\epsilon$ is found from the following formula: $\,\epsilon= (\gamma-1)x_0^2+(2-\gamma)x_0-
\frac{1}{\Gamma(x_0)}\,$, where $x_0=0.2980329689\ldots\,$ is the first positive root
of the equation $\,2(\gamma-1)x\Gamma(x) +(2-\gamma)\Gamma(x) + \Psi(x)=0\,$).
The lower bound is also accurate: the maximum error $\varepsilon$
equals $0.07218627968\ldots\,$ and occurs at $x=0.6388787411\ldots\,$ (see footnote \ref{kjco3enc}).}
\be\notag
\begin{array}{ll}
\displaystyle
\frac{1}{\,n!\,}\!\int\limits_0^1 \! \frac{\,\Gamma(x+n)\,}{\Gamma(x)} \, dx \,
\geqslant \int\limits_0^1 \!  \frac{\,(n+1)^{x-1}\,}{\Gamma(x)} \, dx  \,\geqslant
\int\limits_0^1 \! x \, (n+1)^{x-1}\, dx\,=  \\[7mm]
\qquad\qquad\qquad\quad
\displaystyle
=\,\frac{\,1\,}{\,\ln(n+1)\,}
- \frac{\,1\,}{\,\ln^2 (n+1)\,} + \frac{\,1\,}{\,(n+1)\ln^2 (n+1)\,} 
\end{array}
\ee
but on the other hand
\addtocounter{footnote}{1}
\footnotetext{In some cases, 
more simple estimations may be also useful, for example: $\,x\leqslant \frac{1}{\Gamma(x)}\leqslant 1\,$
or $\,x\leqslant \frac{1}{\Gamma(x)}\leqslant x+\varepsilon\,$, where $\varepsilon=0.07218627968\ldots\,$ (the exact value 
of $\varepsilon$ is given by the formula $\,\varepsilon=\frac{1}{\Gamma(x_0)}-x_0\,$, where 
$x_0=0.6388787411\ldots\,$ is the unique positive root
of the equation $\,\Gamma(x)+\Psi(x)=0\,$).
These estimations are, however, less accurate than \eqref{d3rdf34s}.\label{kjco3enc}}
\be\notag
\begin{array}{ll}
\displaystyle
\frac{1}{\,n!\,}\!\int\limits_0^1 \! \frac{\,\Gamma(x+n)\,}{\Gamma(x)} \, dx \,
\leqslant \int\limits_0^1 \!  \frac{\,n^{x-1}\,}{\Gamma(x)} \, dx  \,\leqslant
\int\limits_0^1 \! \Big[(\gamma-1)x^2+(2-\gamma)x\Big]  n^{x-1}\, dx\,=  \\[7mm]
\qquad\qquad\qquad\quad
\displaystyle
=\,\frac{\,1\,}{\,\ln n\,} - \frac{\,\gamma\,}{\,\ln^2\! n\,} 
 - \frac{\,2(1-\gamma)\,}{\,\ln^3\! n\,} +  \frac{\,2-\gamma \,}{\,n\ln^2\! n\,} 
 + \frac{\,2(1-\gamma)\,}{\,n\ln^3\! n} 
\end{array}
\ee
so that 
\begin{equation}\label{hu765fc}
\begin{array}{ll}
\displaystyle
\frac{\,1\,}{\,\ln(n+1)\,}
- \frac{\,1\,}{\,\ln^2 (n+1)\,} + \frac{\,1\,}{\,(n+1)\ln^2 (n+1)\,} \leqslant \frac{C_{2,n}}{\,n!\,}
\leqslant \frac{\,1\,}{\,\ln n\,} - \frac{\,\gamma\,}{\,\ln^2\! n\,} -  \\[6mm]
\displaystyle\qquad\qquad\qquad\qquad\qquad\qquad
\qquad
 - \frac{\,2(1-\gamma)\,}{\,\ln^3\! n\,} +  \frac{\,2-\gamma \,}{\,n\ln^2\! n\,} 
 + \frac{\,2(1-\gamma)\,}{\,n\ln^3\! n} \,,\qquad n\geqslant2\,,
\end{array}
\end{equation}
Taking into account that last three terms are negative for $n>2.621876631\ldots\,$, the above inequalities also imply that
\be\label{hf38h2d2d}
\,\frac{\,1\,}{\,\ln(n+1)\,}
- \frac{\,1\,}{\,\ln^2 (n+1)\,} 
\leqslant\,\frac{C_{2,n}}{\,n!\,}\,
=\,\frac{\big|B_n^{(n)}\big|}{\,n!\,}\,
\leqslant\,\frac{\,1\,}{\,\ln n\,}
- \frac{\,\gamma\,}{\ln^2\! n} \,,
\qquad\quad n\geqslant3\,.
\ee
Furthermore, since 
\be\notag
\frac{\,1\,}{\,\ln(n+1)\,}
- \frac{\,1\,}{\,\ln^2 (n+1)\,} + \frac{\,1\,}{\,(n+1)\ln^2 (n+1)\,} >
\frac{\,1\,}{\,\ln n\,}
- \frac{\,1\,}{\,\ln^2 \! n\,} \,,\qquad n\geqslant2\,,
\ee
we also have
\be\label{hu765f2dcefc}
\frac{\,1\,}{\,\ln n\,}
- \frac{\,1\,}{\,\ln^2 \! n\,} 
\leqslant\,\frac{C_{2,n}}{\,n!\,}\,
=\,\frac{\big|B_n^{(n)}\big|}{\,n!\,}\,
\leqslant\,\frac{\,1\,}{\,\ln n\,}
- \frac{\,\gamma\,}{\ln^2\! n} \,,
\qquad\quad n\geqslant3\,.
\ee
These bounds are sharper than \eqref{hf38h2d2d}, but weaker than \eqref{hu765fc};
Fig.~\ref{big8ghi87} illustrates their quality.
As $n\to\infty$, the gap between the lower and upper bounds becomes infinitely small,
and hence
$\,
\frac{C_{2,n}}{\,n!\,}
\,\sim\,\frac{\,1\,}{\,\ln n\,}
\,$,
the result which was announced in \eqref{j38ndbr893r}. 
It therefore also follows that series \eqref{fmse} for $\,\ln\ln2\,$ 
converges at the same rate as $\,\sum (-1)^n (n \ln  n)^{-1}\,$.

Let now derive the full asymptotics for $C_{2,n}$ at $n\to\infty$.
For this aim, we first rewrite the last integral in \eqref{hjc239fhsw} in the asymptotical form
\be\label{jec20jdjs}
\int\limits_0^1 \! \frac{\,\Gamma(x+n)\,}{\Gamma(x)} \, dx \,=\,(n-1)!\!
\int\limits_0^1 \!\frac{\,n^x\,}{\Gamma(x)}
\left\{1+ \frac{\,x(x-1)\,}{2 n} +O(n^{-2})\right\}  dx\,,\qquad n\to\infty\,,
\ee
see \eqref{lk2093mffmnjw}.
Then, by expanding the factor $\,1/\Gamma(x)\,$ into the MacLaurin series,
see \eqref{poi2d293dm}--\eqref{jkhe9843hdhr}, and by taking into account that for positive integer $k$
\be
\begin{array}{ll}
\displaystyle
\int\limits_0^1  \! n^{x} x^k \, dx \;& \displaystyle=\,\frac{\,n\cdot k!\,}{\,\ln n\,}\!
\sum_{l=0}^k \frac{(-1)^l\!}{\,(k-l)!\cdot\ln^{l}\! n\,} \,-\,  \frac{(-1)^k k!}{\ln^{k+1}\! n} \\[7mm]
&\displaystyle\label{djoiejdnss}
=\,n\left\{\frac{1}{\ln n}-\frac{k}{\ln^2 \!n}+\frac{k(k-1)}{\ln^3\! n}
-\frac{k(k-1)(k-2)}{\ln^4\! n} + \ldots+ \frac{(-1)^k k!}{\ln^{k+1}\! n}\right\}  -  \frac{(-1)^k k!}{\ln^{k+1}\! n} 
\end{array}
\ee
we have for the last integral in \eqref{jec20jdjs}
\begin{eqnarray}
\displaystyle\notag 
\int\limits_0^1 &&\displaystyle \!\frac{\,n^x\,}{\Gamma(x)}
\left\{1+ \frac{\,x(x-1)\,}{2 n} +O(n^{-2})\right\} dx\,=  \int\limits_0^1 \! n^x
\left\{1- \frac{\,x\,}{2 n}+ \frac{\,x^2\,}{2 n} +O(n^{-2})\right\} \sum_{k=1}^\infty x^k a_k \, dx\\[3mm]
&&\displaystyle\notag
=  \sum_{k=1}^\infty  a_k \!\int\limits_0^1  \! n^x x^k \, dx \, 
-\frac{1}{2n}\sum_{k=1}^\infty  a_k \!\int\limits_0^1  \! n^x x^{k+1} \, dx 
\, +\frac{1}{2n}\sum_{k=1}^\infty  a_k \!\int\limits_0^1  \! n^x x^{k+2} \, dx \,=
\end{eqnarray}
\begin{eqnarray}
&&\displaystyle\notag
=  \frac{n}{\ln n} \underbrace{\sum_{k=1}^\infty  a_k}_{\frac{1}{\Gamma(1)}}  
- \frac{n}{\ln^2\! n} \underbrace{\sum_{k=1}^\infty k a_k}_{\left[\frac{1}{\Gamma(x)}\right]'_{x=1}} 
+ \frac{n}{\ln^3\! n} \underbrace{\sum_{k=1}^\infty k(k-1) a_k}_{\left[\frac{1}{\Gamma(x)}\right]''_{x=1}} 
- \ldots \\[3mm]
&&\displaystyle\quad\notag
-\frac{1}{2n}\Biggl\{  \frac{n}{\ln n} \underbrace{\sum_{k=1}^\infty  a_k}_{\frac{1}{\Gamma(1)}}  
- \frac{n}{\ln^2\! n} \underbrace{\sum_{k=1}^\infty (k+1) a_k}_{\left[\frac{x}{\Gamma(x)}\right]'_{x=1}} 
+ \frac{n}{\ln^3\! n} \underbrace{\sum_{k=1}^\infty (k+1)k a_k}_{\left[\frac{x}{\Gamma(x)}\right]''_{x=1}} 
- \ldots \Biggr\} \\[3mm]
&&\displaystyle\quad\notag
+\frac{1}{2n}\Biggl\{  \frac{n}{\ln n} \underbrace{\sum_{k=1}^\infty  a_k}_{\frac{1}{\Gamma(1)}}  
- \frac{n}{\ln^2\! n} \underbrace{\sum_{k=1}^\infty (k+2) a_k}_{\left[\frac{x^2}{\Gamma(x)}\right]'_{x=1}} 
+ \frac{n}{\ln^3\! n} \underbrace{\sum_{k=1}^\infty (k+2)(k+1) a_k}_{\left[\frac{x^2}{\Gamma(x)}\right]''_{x=1}} 
- \ldots \Biggr\}+O\!\left(\frac{1}{n\ln n}\right) \\[3mm]
&&\displaystyle\notag
= \, \frac{n}{\ln n} + n\sum\limits_{l=1}^\infty \frac{(-1)^{l}}{\ln^{l+1}\!n}\cdot\left[\frac{1}{\Gamma(x)}\right]^{(l)}_{x=1}
-\frac{1}{2}\sum\limits_{l=1}^\infty \frac{(-1)^{l}}{\ln^{l+1}\!n}\cdot\left[\frac{x}{\Gamma(x)}\right]^{(l)}_{x=1}\\[3mm]
&&\displaystyle\quad\label{lc0239jnmc}
+\frac{1}{2}\sum\limits_{l=1}^\infty \frac{(-1)^{l}}{\ln^{l+1}\!n}\cdot\left[\frac{x^2}{\Gamma(x)}\right]^{(l)}_{x=1}
+O\!\left(\frac{1}{n\ln n}\right)\
\end{eqnarray}
in virtue of the uniform convergence of series \eqref{poi2d293dm}.
Retaining only first significant terms, which are all contained in the first sum, we finally establish that 
\begin{eqnarray}
\displaystyle\label{xe21fd4}
\frac{C_{2,n}}{\, n!\,}\,
=\,\frac{\big|B_n^{(n)}\big|}{\,n!\,}\,=\,
\frac{1}{\,n!\,}\sum_{l=1}^{n} \frac{\big|S_1(n,l)\big|}{l+1} \;&&\displaystyle =\,  
\frac{1}{\,\ln n\,} + \sum\limits_{l=1}^\infty \frac{(-1)^{l}}{\ln^{l+1}n}\cdot\left[\frac{1}{\Gamma(x)}\right]^{(l)}_{x=1} 
+  O\!\left(\frac{1}{\,n\ln^2\! n\,}\right)		\\[3mm]
&&\displaystyle\notag
=\,  \frac{1}{\,\ln n\,} - \frac{\gamma}{\,\ln^2 \!n\,}
-\frac{\,\pi^2- 6\gamma^2\,}{\,6\ln^3 \!n\,} +  O\!\left(\frac{1}{\,\ln^4\!n\,}\right)\,,\qquad n\to\infty\,,
\end{eqnarray}

\noindent\textbf{Historical remark} The first--order asymptotics for $C_{2,n}$ was known to Binet as early as 1839 
(see the final remark on p.~\pageref{ckj3w0emnd}), and it may be also found in \cite[p.~294]{comtet_01}.
As regards the higher--order terms, as well as upper and lower bounds for $C_{2,n}$, inequalities \eqref{hu765fc}--\eqref{hu765f2dcefc}, 
we have no found them in previously published literature.

\subsection{Bounds and full asymptotical expansion for Gregory's coefficients $G_n$, also known as reciprocal logarithmic numbers,
Cauchy numbers of the first kind $C_{1,n}$ and generalized Bernoulli numbers $B_n^{(n-1)}\,$}\label{appendix2}
\begin{figure}[!t]   
\centering
\includegraphics[width=0.8\textwidth]{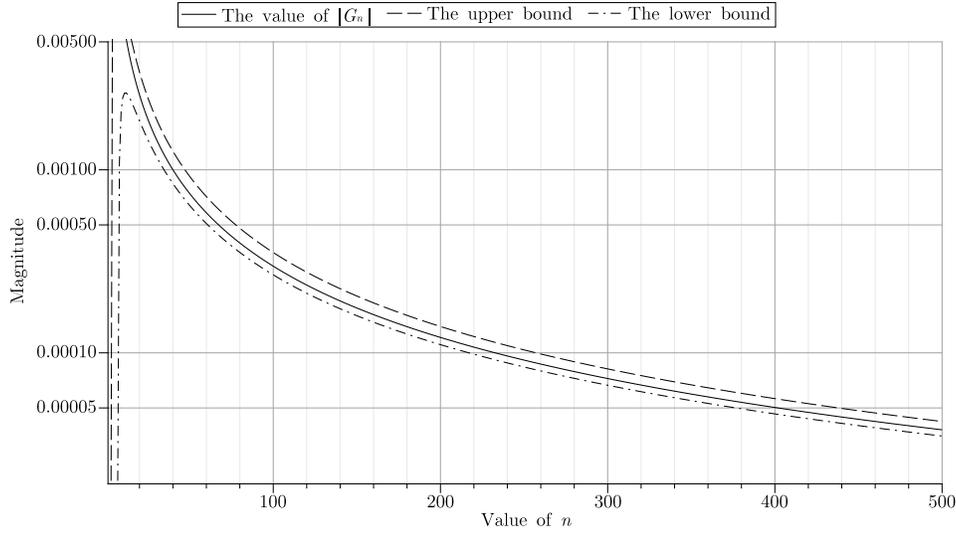}
\caption{Lower and upper bounds for numbers $G_n$ given by inequalities \eqref{uih297hdr2}, 
logarithmic scale.}
\label{oih98hhz}
\end{figure}

A method analogous to that we just employed may also provide equivalent results for Gregory's coefficients $G_n$,
which appear in equations \eqref{fms}, \eqref{fms2x}, \eqref{k2390234mrf}, \eqref{kjc0293nfr}, 
\eqref{uih39hlj2983dh}, \eqref{s2kj0sdw}--\eqref{hgwibddw} and \eqref{c234c2gf5g}. 
First, reducing the signed Stirling numbers to the unsigned ones, and then, performing the same procedure 
as in \eqref{hjc239fhsw}, we have
\be
\begin{array}{ll}
\displaystyle
G_n\, n!\;& \displaystyle=\,C_{1,n}\,=\,-\frac{B^{(n-1)}_n}{n-1}\,=
\sum_{l=1}^{n} \frac{S_1(n,l)}{l+1}  = \, (-1)^n\!\sum_{l=1}^{n} |S_1(n,l)|\!\int\limits_0^1  \! (-x)^l \, dx \,=\\[6mm]
& \displaystyle
=\, (-1)^n\!\int\limits_0^1 \!\sum_{l=1}^{n} (-x)^l  |S_1(n,l)| \, dx\,
=\, (-1)^n\!\int\limits_0^1 \! (-x)_n \, dx \,
=\, (-1)^n\!\int\limits_0^1 \! (-x)_n \, dx \,=\\[8mm]
& \displaystyle
=\, (-1)^n\!\int\limits_0^1 \! \frac{\,\Gamma(n-x)\,}{\Gamma(-x)} \, dx \,
=\, (-1)^{n-1}\!\int\limits_0^1 \! \frac{\,(1-z)\,\Gamma(n-1+z)\,}{\Gamma(z)} \, dz	\label{jhohbhgf3}
\end{array}
\ee
where, at the last stage, we first made a change of variable $x=1-z$, and then used the recurrence
relationship $\Gamma(z-1)=\Gamma(z)/(z-1)$. Using bounds \eqref{jd83dwsw} and \eqref{d3rdf34s},
as well as formula \eqref{djoiejdnss}, we deduce 
two following inequalities
\be\notag
\begin{array}{ll}
\displaystyle
\frac{1}{\,n!\,}\!\int\limits_0^1 \! \frac{\,(1-z)\,\Gamma(n-1+z)\,}{\Gamma(z)} \, dz \,
\geqslant \int\limits_0^1 \!  \frac{\,(1-z)\,n^{z-2}\,}{\Gamma(z)} \, dz  \,\geqslant
\int\limits_0^1 \! z\,(1-z) \, n^{z-2}\, dz\,=  \\[7mm]
\qquad\qquad\qquad\quad
\displaystyle
=\,\frac{\,1\,}{\,n\ln^2\! n\,} \,-\, \frac{\,2\,}{\,n\ln^3\! n\,}
\,+\, \frac{\,1\,}{\,n^2\ln^2\! n\,}\,+\, \frac{\,2\,}{\,n^2\ln^3\! n\,}
\end{array}
\ee
and
\be\notag
\begin{array}{ll}
\displaystyle
\frac{1}{\,n!\,}\!\int\limits_0^1 \! \frac{\,(1-z)\,\Gamma(n-1+z)\,}{\Gamma(z)} \, dz \,
\leqslant \,\frac{1}{\,n\,}\!\int\limits_0^1 \!  \frac{\,(1-z)\,(n-1)^{z-1}\,}{\Gamma(z)} \, dz \,\leqslant \\[7mm]
\qquad\quad
\displaystyle
\,\leqslant
\,\frac{1}{\,n\,}\!\int\limits_0^1 \! \Big[(\gamma-1)z^2+(2-\gamma)z\Big]  \,(1-z) \,(n-1)^{z-1}\, dz\,=  \\[7mm]
\qquad\quad
\displaystyle
=\,\frac{\,1\,}{\,(n-1)\ln^2(n-1)\,} - \frac{\,2\gamma\,}{\,(n-1)\ln^3(n-1)\,}
\,-\frac{\,6\,(1-\gamma)\,}{\,(n-1)\ln^4(n-1)\,}\, + \,\\[7mm]
\qquad\qquad
\displaystyle
+\,\frac{\,1-\gamma\,\,}{\,n(n-1)\ln^2(n-1)\,}
+\,\frac{\,2\,(3-\gamma)\,\,}{\,n(n-1)\ln^3(n-1)\,}
+\,\frac{\,12\,(1-\gamma)\,\,}{\,n(n-1)\ln^4(n-1)\,}
\end{array}
\ee
Whence, taking into account that numbers $G_n$ are strictly alternating $G_n=(-1)^{n-1} \big|G_n\big|$, 
we deduce that Gregory's coefficients enjoy these bounds
\begin{eqnarray}\label{jhged28gd}
\displaystyle
&&\frac{\,1\,}{\,n\ln^2\! n\,} \,-\, \frac{\,2\,}{\,n\ln^3\! n\,}
\,+\, \frac{\,1\,}{\,n^2\ln^2\! n\,}\,+\, \frac{\,2\,}{\,n^2\ln^3\! n\,}
 \leqslant \big| G_n\big| \leqslant \frac{\,1\,}{\,(n-1)\ln^2(n-1)\,} - \frac{\,2\gamma\,}{\,(n-1)\ln^3(n-1)\,}
\,+\qquad\\[3mm]
&&\qquad
\displaystyle\notag
+\left\{-\frac{\,6\,(1-\gamma)\,}{\,(n-1)\ln^4(n-1)\,}\, 
+\,\frac{\,1-\gamma\,\,}{\,n(n-1)\ln^2(n-1)\,}
+\,\frac{\,2\,(3-\gamma)\,\,}{\,n(n-1)\ln^3(n-1)\,}
+\,\frac{\,12\,(1-\gamma)\,\,}{\,n(n-1)\ln^4(n-1)\,}\right\}
\end{eqnarray}
Since the contribution of the terms in curly brackets is negative for $n>4.921304199\ldots$\,
(thanks to the leading first term), the above bounds also imply a weaker relation
\be\label{uih297hdr}
\frac{\,1\,}{\,n\ln^2\! n\,} \,-\, \frac{\,2\,}{\,n\ln^3\! n\,}
\leqslant\,\big|G_n\big|\,
\leqslant\, \frac{\,1\,}{\,(n-1)\ln^2(n-1)\,} \,-\, \frac{\,2\gamma\,}{\,(n-1)\ln^3(n-1)\,} \,,
\qquad\quad n\geqslant5\,.
\ee
Moreover, the detailed study of the right part in \eqref{jhged28gd} leads to another inequalities
\be\label{uih297hdr2}
\frac{\,1\,}{\,n\ln^2\! n\,} \,-\, \frac{\,2\,}{\,n\ln^3\! n\,}
\leqslant\,\big|G_n\big|\,
\leqslant\, \frac{\,1\,}{\,n\ln^2\! n\,} \,-\, \frac{\,2\gamma\,}{\,n\ln^3\! n\,} \,,
\qquad\quad n\geqslant5\,,
\ee
which are slightly stronger and simpler than \eqref{uih297hdr}, but weaker than parent inequalities \eqref{jhged28gd}.
Graphical simulations, see Fig.~\ref{oih98hhz}, show that bounds \eqref{uih297hdr2} are very sharp and their accuracy 
should be sufficient for most of the situations..
Making $n\to\infty$ yields the first--order approximation for Gregory's coefficients: $|G_n|\sim\frac{1}{n\ln^2 n}\,$.
Higher--order terms of this asymptotics may be found as follows. By \eqref{lk2093mffmnjw}, 
the last integral in \eqref{jhohbhgf3} becomes
\be\notag
\begin{array}{ll}
\displaystyle
&\displaystyle (-1)^{n-1}\!\int\limits_0^1 \! \frac{\,(1-z)\,\Gamma(n-1+z)\,}{\Gamma(z)} \, dz=\\[6mm]
& \displaystyle \label{d2x3dj2ed}\qquad\quad
=\,(-1)^{n-1} (n-1)! \!\int\limits_0^1 \!\frac{\,(1-z)\,n^{z-1}\,}{\Gamma(z)}
\left\{1+ \frac{\,(1-z)\,(2-z)\,}{2 n} +O(n^{-2})\right\} dz
\,,\qquad n\to\infty\,,
\end{array}
\ee
Then, using the MacLaurin series for $1/\Gamma(z)$ and proceeding analogously to \eqref{lc0239jnmc}, we find 
\begin{eqnarray}
\displaystyle\notag
G_n\,=\,\frac{1}{n!}\sum_{l=1}^{n} \frac{S_1(n,l)}{l+1} \;&&\displaystyle =\,  \frac{\,(-1)^{n-1}}{n}\cdot
\sum\limits_{l=1}^\infty \frac{(-1)^{l}}{\ln^{l+1}n}\cdot
\left[\frac{1-x}{\Gamma(x)}\right]^{(l)}_{x=1}
+  O\!\left(\frac{1}{\,n^2\ln n\,}\right) 		\\[3mm]
&&\displaystyle\label{ojc02m3d}
=\,  \frac{\,(-1)^{n-1}}{n}\cdot
\left\{\frac{1}{\,\ln^2 \! n\,} - \frac{2\gamma}{\,\ln^3\! n\,} - \frac{\,\pi^2-6\gamma^2}{\,2\ln^4\! n\,} 
+  O\!\left(\frac{1}{\,\ln^5\!n\,}\right) \!\right\}
\end{eqnarray}

\noindent\textbf{Historical remark}
The first--order approximation for Gregory's coefficients $G_n$ at $n\to\infty$, the first formul\ae~in \eqref{j38ndbr893r}, was found by
Ernst Schr\"oder in 1879 \cite[p.~115, Eq.~(25a)]{schroder_01}.
It was rediscovered by Johan Steffensen in 1924 \cite[pp.~2--4]{steffensen_01},  \cite[pp.~106--107]{steffensen_02},
and was slightly bettered in 1957 by Davis \cite[p.~14, Eq.~(14)]{davis_02}. Higher--order terms
of this asymptotics were obtained by S.~C.~Van Veen in 1950 \cite[p.~336]{van_veen_01}, \cite[p.~29]{norlund_01},
Gergő Nemes in 2011 \cite{nemes_01} and the author in 2014. S.~C.~Van Veen and Gergő Nemes used different methods
to derive their results and obtained different expressions; however, one can show that their formul\ae~are equivalent. 
The former employed an elegant contour integration method, while  
the latter used Watson's lemma. Our formula \eqref{ojc02m3d} differs from both Van Veen's result and 
Nemes' result, but it is also equivalent to them. 
Note also that many researchers
(e.g.~Nørlund \cite[p.~29]{norlund_01}, Davis \cite[p.~14, Eq.~(14)]{davis_02}, Nemes \cite{nemes_01})
incorrectly attribute this asymptotics to Steffensen, who only rediscovered it. The same asymptotics 
also appears in the well--known monograph \cite[p.~294]{comtet_01}, \label{lkjhco2encn}
but Comtet did not specify the source of the formula. As regards the bounds for $G_n$, in 1922 
Steffensen found that
\be\notag
\frac{\,1\,}{\,6n\,(n-1)\,}
<\,\big|G_n\big|\,
<\, \frac{\,1\,}{\,6n\,}\,,\qquad\quad n>2\,,
\ee
\cite[p.~198, Eq.~(27)]{steffensen_00}, \cite[p.~2, Eq.~(3)]{steffensen_01}, \cite[p.~106, Eq.~(9)]{steffensen_02}.
A stronger result (at large $n$) was stated by Kluyver in 1924
\be\notag
\big|G_n\big|\,<\,\frac{1}{\,n\ln n\,}\,,\qquad\quad n\geqslant2\,,
\ee
\cite[p.~144]{kluyver_02}. In 2010, Rubinstein \cite[p.~30, Theorem 1.1]{rubinstein_01} found a bound for more general numbers, 
from which it \emph{inter alia} follows that\footnote{Rubinstein's $\alpha_n(s)$ at $s=0$ are our $-\big|G_n\big|$.}
\be\notag
\big|G_n\big|\,\leqslant\frac{\,4\,(1+\ln(1+n))\,}{1+n}\,,\qquad\quad n\geqslant1\,.
\ee
This bound is nevertheless much weaker
than both preceeding bounds. In another recent paper \cite[p.~473]{coffey_02}, Coffey remarked that numerical simulations suggest that
$\big|G_n\big|$ should be lesser than $\,\frac{1}{n\ln^2 n}\,$ for all $n\geqslant2$, but that the proof of this 
result was missing.\footnote{Note that Coffey's $p_{n+1}$ are our $|G_n|$ (Coffey's notation are probably borrowed from 
Ser's paper \cite{ser_01}).} 
Inequalities \eqref{uih297hdr2} include the missing proof. 
By the way, as far as we know,
our bounds \eqref{jhged28gd}--\eqref{uih297hdr2} are currently the best bounds for Gregory's coefficients.\\

\noindent\textbf{Nota Bene} Analogously, one can show that more general results take place. 
For instance, for positive integer $k$
\be\notag
\sum_{l=1}^{n} \frac{\,\big|S_1(n,l)\big|\,}{l+k}\, =\int\limits_0^1 \! x^{k-1} (x)_n \, dx \,
\sim\,(n-1)!\!\int\limits_0^1 \!\frac{\,n^{x} x^{k-1}\,}{\Gamma(x)} \, dx\,,\qquad n\to\infty\,.
\ee
Whence, using  \eqref{poi2d293dm}--\eqref{jkhe9843hdhr} and \eqref{djoiejdnss}, we obtain
\begin{eqnarray}
\displaystyle\notag
\frac{1}{n!}\sum_{l=1}^{n} \frac{\big|S_1(n,l)\big|}{l+k} \;\,&&\displaystyle \sim\,
\frac{1}{\,\ln n\,}+ 
\sum\limits_{l=1}^\infty \frac{(-1)^{l}}{\ln^{l+1}n}\cdot
\left[\frac{x^{k-1}}{\Gamma(x)}\right]^{(l)}_{x=1}=\,  \frac{1}{\,\ln n\,} - \frac{\,k-1+\gamma\,}{\,\ln^2 \! n\,}  \\[3mm]
&&\displaystyle\quad
+ \, \frac{\,6\gamma^2 +12\gamma(k-1)-\pi^2+6(k^2-3k+2)\,}{\,6\ln^3 \! n\,} 
+  O\!\left(\frac{1}{\,\ln^4\!n\,}\right) \,,
\qquad 
\begin{array}{l}
k=1, 2, 3,\ldots\\[1mm]
n\to\infty\,.
\end{array}
\qquad
\end{eqnarray}
In particular, for $k=2$
\be\label{klxj2103xnd}
\frac{1}{n!}\sum_{l=1}^{n} \frac{\big|S_1(n,l)\big|}{l+2}\,
=\,  \frac{1}{\,\ln n\,} - \frac{1+\gamma}{\,\ln^2 \! n\,} - \frac{\,\pi^2-6\gamma^2-12\gamma\,}{\,6\ln^3 \! n\,} 
+  O\!\left(\frac{1}{\,\ln^4\!n\,}\right) 
\,,\qquad n\to\infty
\ee
Similarly,
\be\notag
\begin{array}{ll}
\displaystyle
\sum_{l=1}^{n} \frac{S_1(n,l)}{l+k}\;\,&\displaystyle =\, (-1)^n\!\int\limits_0^1 \! x^{k-1} (-x)_n \, dx \,
\sim\,\frac{\,(-1)^{n-1} (n-1)!\,}{n}\!\int\limits_0^1 \!\frac{\,n^{z} (1-z)^k\,}{\Gamma(z)} \, dz\\[6mm]
&\displaystyle\quad=\,\frac{\,(-1)^{n-1} (n-1)!\,}{n}\sum_{m=0}^k\binom{k}{m}(-1)^m
\!\int\limits_0^1 \!\frac{\,z^m n^{z} \,}{\Gamma(z)} \, dz\,,\qquad n\to\infty\,,
\end{array}
\ee
whence
\begin{eqnarray}
\displaystyle\notag
\frac{1}{n!}\sum_{l=1}^{n} \frac{S_1(n,l)}{l+k} \;\,&&\displaystyle \sim\,  \frac{\,(-1)^{n-1}}{n}
\sum_{m=0}^k\binom{k}{m}(-1)^m	\left\{ \frac{1}{\,\ln n\,}+ 
\sum\limits_{l=1}^\infty \frac{(-1)^{l}}{\ln^{l+1}n}\cdot
\left[\frac{x^m}{\Gamma(x)}\right]^{(l)}_{x=1} \right\}\\[3mm]
&&\displaystyle\label{kf3094fmf}
=\, \frac{(-1)^{n-1}}{n}\left\{\frac{k!}{\,\ln^{k+1} \! n\,} - \frac{\gamma\,(k+1)!}{\,\ln^{k+2} \! n\,} 
+  O\!\left(\frac{1}{\,\ln^{k+3} \! n\,}\right)\!\right\}\,,
\qquad 
\begin{array}{l}
k=1, 2, 3,\ldots\\[1mm]
n\to\infty\,,
\end{array}
\end{eqnarray}
In particular,
\be\label{klxj2103xnd}
\frac{1}{n!}\sum_{l=1}^{n} \frac{S_1(n,l)}{l+2}\,
=\,  \frac{\,(-1)^{n-1}}{n}
\left\{\frac{2}{\,\ln^3 \! n\,}- \frac{6\gamma}{\,\ln^4 \! n\,} - \frac{2(\pi^2-6\gamma^2)}{\,\ln^5 \! n\,} 
+  O\!\left(\frac{1}{\,\ln^6\!n\,}\right) \!\right\}
\,,\qquad n\to\infty
\ee

Finally, we remark that different asymptotical aspects, in which are involved Stirling numbers, are also discussed, at different extents, in works 
of Jordan \cite{jordan_00}, \cite[Chapt.~IV]{jordan_01}, Moser \& Wyman \cite{moser_01}, Wilf \cite{wilf_02}, Temme \cite{temme_02},
Hwang \cite{hwang_01}, Timashev \cite{timashev_01}, Grünberg \cite{grunberg_01} and Louchard \cite{louchard_01}.
Readers interested in a more deep
study of general asymptotical methods might also wish to consult the following literature: \cite[Chapt.~I, \S~4]{evgrafov_01_eng_1st}, 
\cite{evgrafov_03_eng}, \cite{olver_01}, \cite{dingle_01}.

\bibliographystyle{crelle}

\begin{thebibliography}{100}
\providecommand{\url}[1]{\texttt{#1}}
\providecommand{\urlprefix}{URL }
\expandafter\ifx\csname urlstyle\endcsname\relax
  \providecommand{\doi}[1]{doi:\discretionary{}{}{}#1}\else
  \providecommand{\doi}{doi:\discretionary{}{}{}\begingroup
  \urlstyle{rm}\Url}\fi

\bibitem{abramowitz_01}
\textit{M.~Abramowitz} and \textit{I.~A. Stegun}, Handbook of mathematical
  functions with formula, graphs and mathematical tables [Applied mathematics
  series \no 55], US Department of Commerce, National Bureau of Standards,
  1961.

\bibitem{adamchik_03}
\textit{V.~Adamchik}, On \text{Stirling} numbers and \text{Euler's} sums,
  Journal of Computational and Applied Mathematics, vol.~79, pp.~119--130
  (1997).

\bibitem{adelberg_01}
\textit{A.~Adelberg}, 2-\text{Adic} congruences of \text{Nörlund} numbers and
  of \text{Bernoulli} numbers of the second kind, Journal of Number Theory,
  vol.~73, no.~1, pp.~47--58  (1998).

\bibitem{alabdulmohsin_01}
\textit{I.~M. Alabdulmohsin}, Summability calculus, arXiv:1209.5739v1  (2012).

\bibitem{appel_01}
\textit{P.~Appel}, Développement en série entière de $(1+ax)^{\frac{1}{x}}$,
  Archiv der Mathematik und Physik, vol.~65, pp.~171--175  (1880).

\bibitem{arakawa_01}
\textit{T.~Arakawa}, \textit{T.~Ibukiyama} and \textit{M.~Kaneko}, Bernoulli
  Numbers and Zeta Functions, Springer Monographs in Mathematics, Japan, 2014.

\bibitem{artin_01}
\textit{E.~Artin}, Einführung in die \text{Theorie} der \text{Gammafunktion},
  B. G. Teubner, Leipzig, Germany, 1931.

\bibitem{bateman_01}
\textit{H.~Bateman} and \textit{A.~Erdélyi}, Higher Transcendental Functions
  [in 3 volumes], Mc Graw--Hill Book Company, 1955.

\bibitem{batir_01}
\textit{N.~Batir}, Very accurate approximations for the factorial function,
  Journal of Mathematical Inequalities, vol.~4, no.~3, pp.~335--344  (2010).

\bibitem{bellavista_01}
\textit{L.~V. Bellavista}, On the \text{Stirling} numbers of the first kind
  arising from probabilistic and statistical problems, Rendiconti del Circolo
  Matematico di Palermo, vol.~32, no.~1, pp.~19--26  (1983).

\bibitem{bender_01}
\textit{E.~A. Bender} and \textit{S.~G. Williamson}, Foundations of
  Combinatorics with Applications, Addison--Wesley, USA, 1991.

\bibitem{binet_01}
\textit{M.~J. Binet}, Mémoire sur les intégrales définies eulériennes et
  sur leur application à la théorie des suites, ainsi qu'à l'évaluation des
  fonctions des grands nombres, Journal de l'École Royale Polytechnique, tome
  XVI, cahier 27, pp.~123--343  (1839).

\bibitem{iaroslav_06}
\textit{\text{Ia.~V.} Blagouchine}, Rediscovery of \text{Malmsten's} integrals,
  their evaluation by contour integration methods and some related results, The
  Ramanujan Journal, vol.~35, no.~1, pp.~21--110 (erratum in press, erratum's
  DOI: 10.1007/s11139-015-9763-z)  (2014).

\bibitem{iaroslav_07}
\textit{\text{Ia.~V.} Blagouchine}, A theorem for the closed--form evaluation
  of the first generalized \text{Stieltjes} constant at rational arguments and
  some related summations, Journal of Number Theory (Elsevier), vol.~148,
  pp.~537--592 and vol.~151, pp.~276--277, arXiv:1401.3724, 2014  (2015).

\bibitem{iaroslav_09}
\textit{\text{Ia.~V.} Blagouchine}, Expansions of generalized \text{Euler's}
  constants into the series of polynomials in $\pi^{-2}$ and into the formal
  enveloping series with rational coefficients only, Journal of Number Theory
  (Elsevier), vol.~158, pp.~365--396, arXiv:1501.00740, 2015  (2016).

\bibitem{boole_01}
\textit{G.~Boole}, Calculus of finite differences (edited by J.~F.~Moulton, 4th
  ed.), Chelsea Publishing Company, New-York, USA, 1957.

\bibitem{bromwich_02}
\textit{T.~J. \text{I'a} Bromwich}, A note on \text{Stirling's} series and
  \text{Euler's} constant, The Messenger of Mathematics, vol.~36, pp.~81--85
  (1906).

\bibitem{bromwich_01}
\textit{T.~J. \text{I'a} Bromwich}, An introduction to the theory of infinite
  series, Macmillan and Co. Limited, St--Martin Street, London, 1908.

\bibitem{brychkov_02}
\textit{\text{Yu.~A.} Brychkov}, Power expansions of powers of trigonometric
  functions and series containing \text{Bernoulli} and \text{Euler}
  polynomials, Integral Transforms and Special Functions, vol.~20, no.~11,
  pp.~797--804  (2009).

\bibitem{brychkov_01}
\textit{\text{Yu.~A.} Brychkov}, On some properties of the generalized
  \text{Bernoulli} and \text{Euler} polynomials, Integral Transforms and
  Special Functions, vol.~23, no.~10, pp.~723--735  (2012).

\bibitem{burnside_01}
\textit{W.~Burnside}, A rapidly convergent series for \text{$\log N!$}, The
  Messenger of Mathematics, vol.~46, pp.~157--159  (1916--1917).

\bibitem{butzer_01}
\textit{P.~M. Butzer} and \textit{M.~Hauss}, \text{Stirling} functions of the
  first and second kinds; some new applications, In Israel Mathematical
  Conference Proceedings, Approximation, Interpolation, and Summability, in
  Honor of Amnon Jakimovski on his Sixty-Fifth Birthday (Ed. S. Baron and D.
  Leviatan), Bar-Ilan University, Tel Aviv, June 4--8 1990, pp.~89--108,
  Weizmann Press, Israel  (1991).

\bibitem{butzer_02}
\textit{P.~M. Butzer}, \textit{C.~Markett} and \textit{M.~Schmidt},
  \text{Stirling} numbers, central factorial numbers, and representations of
  the riemann zeta function, Results in mathematics, vol.~19, no.~3--4,
  pp.~257--274  (1991).

\bibitem{campbell_01}
\textit{R.~Campbell}, Les intégrales eulériennes et leurs applications,
  Dunod, Paris, 1966.

\bibitem{candelpergher_01}
\textit{B.~Candelpergher} and \textit{M.-A. Coppo}, A new class of identities
  involving \text{Cauchy} numbers, harmonic numbers and zeta values, The
  Ramanujan Journal, vol.~27, pp.~305--328  (2012).

\bibitem{carlitz_02}
\textit{L.~Carlitz}, Some theorems on \text{Bernoulli} numbers of higher order,
  Pacific Journal of Mathematics, vol.~2, no.~2, pp.~127--139  (1952).

\bibitem{carlitz_01}
\textit{L.~Carlitz}, A note on \text{Bernoulli} and \text{Euler} polynomials of
  the second kind, Scripta Mathematica, vol.~25, pp.~323--330  (1961).

\bibitem{carlitz_03}
\textit{L.~Carlitz} and \textit{F.~R. Olson}, Some theorems on \text{Bernoulli}
  and \text{Euler} numbers of higher order, Duke Mathematical Journal, vol.~21,
  no.~3, pp.~405--421  (1954).

\bibitem{cauchy_01}
\textit{A.-L. Cauchy}, Mémoire sur la théorie de la propagation des ondes à
  la surface d'un fluide pesant, d'une profondeur indéfinie, Mémoires
  présentés par divers savants à l'Académie Royale des Sciences de
  l'Institut de France, Sciences Mathématiques et Physique, T.~1, p.~130
  (1827).

\bibitem{cayley_00}
\textit{A.~Cayley}, On a theorem for the development of a factorial,
  Philosophical magazine, vol.~6, pp.~182--185  (1853).

\bibitem{cayley_01}
\textit{A.~Cayley}, On some numerical expansions, The Quarterly journal of pure
  and applied mathematics, vol.~3, pp.~366--369  (1860).

\bibitem{cayley_02}
\textit{A.~Cayley}, Note on a formula for \text{$\Delta^n 0^i/n^i$} when $n$,
  $i$ are very large numbers, Proceedings of the Royal Society of Edinburgh,
  vol.~14, pp.~149--153  (1887).

\bibitem{chabert_01}
\textit{J.-L. Chabert}, \textit{E.~Barbin}, \textit{M.~Guillemot},
  \textit{A.~Michel-Pajus}, \textit{J.~Borowczyk}, \textit{A.~Jebbar} and
  \textit{J.-C. Martzloff}, Histoire d'algorithmes: du caillou à la puce,
  Belin, France, 1994.

\bibitem{charalambides_01}
\textit{C.~A. Charalambides}, Enumerative Combinatorics, Chapman \& Hall/CRC,
  USA, 2002.

\bibitem{coffey_08}
\textit{M.~W. Coffey}, Addison-type series representation for the
  \text{Stieltjes} constants, Journal of Number Theory, vol.~130,
  pp.~2049--2064  (2010).

\bibitem{coffey_07}
\textit{M.~W. Coffey}, Certain logarithmic integrals, including solution of
  monthly problem tbd, zeta values, and expressions for the \text{Stieltjes}
  constants, arXiv:1201.3393v1  (2012).

\bibitem{coffey_02}
\textit{M.~W. Coffey}, Series representations for the \text{Stieltjes}
  constants, Rocky Mountain Journal of Mathematics, vol.~44, pp.~443--477
  (2014).

\bibitem{comtet_01}
\textit{L.~Comtet}, Advanced Combinatorics. The art of Finite and Infinite
  Expansions (revised and enlarged edition), D. Reidel Publishing Company,
  Dordrecht, Holland, 1974.

\bibitem{conway_01}
\textit{J.~H. Conway} and \textit{R.~K. Guy}, The Book of Numbers, Springer,
  New--York, 1996.

\bibitem{copson_01}
\textit{E.~T. Copson}, Asymptotic Expansions, Cambridge University Press, Great
  Britain, 1965.

\bibitem{davis_02}
\textit{H.~T. Davis}, The approximation of logarithmic numbers, American
  Mathematical Monthly, vol.~64, no.~8, part II, pp.~11--18  (1957).

\bibitem{davis_01}
\textit{P.~J. Davis}, \text{Leonhard Euler's} integral: A historical profile of
  the \text{Gamma} function, American Mathematical Monthly, vol.~66,
  pp.~849--869  (1959).

\bibitem{dingle_01}
\textit{R.~B. Dingle}, Asymptotic Expansions: their Derivation and
  Interpretation, Academic Press, USA, 1973.

\bibitem{espinosa_01}
\textit{O.~Espinosa} and \textit{V.~H. Moll}, On some integrals involving the
  \text{Hurwitz} zeta function: \text{Part I}, The Ramanujan Journal, vol.~6,
  pp.~150--188  (2002).

\bibitem{ettingshausen_01}
\textit{A.~von Ettingshausen}, Die combinatorische Analysis als
  Vorbereitungslehre zum Studium der theoretischen höhern Mathematik,
  J.~B.~Wallishausser, Vienna, 1826.

\bibitem{euler_02}
\textit{L.~Eulero}, Institutiones calculi differentialis cum eius usu in
  analysi finitorum ac doctrina serierum, Academi\ae~Imperialis Scientiarum
  Petropolitan\ae, 1755.

\bibitem{evgrafov_01_eng}
\textit{M.~A. Evgrafov}, \textit{K.~A. Bezhanov}, \textit{Y.~V. Sidorov},
  \textit{M.~V. Fedoriuk} and \textit{M.~I. Shabunin}, A Collection of Problems
  in the Theory of Analytic Functions (second edition) [in Russian], Nauka,
  Moscow, USSR, 1972.

\bibitem{evgrafov_03_eng}
\textit{M.~A. Evgrafov} and \textit{A.~Shields}, Asymptotic estimates and
  entire functions, Gordon and Breach, USA, 1961.

\bibitem{evgrafov_01_eng_1st}
\textit{M.~A. Evgrafov}, \textit{Y.~V. Sidorov}, \textit{M.~V. Fedoriuk},
  \textit{M.~I. Shabunin} and \textit{K.~A. Bezhanov}, A Collection of Problems
  in the Theory of Analytic Functions [in Russian], Nauka, Moscow, USSR, 1969.

\bibitem{forsyth_01}
\textit{A.~R. Forsyth}, On an approximate expression for $x!$, Report of the
  Fifty--Third Meeting of the British Association for the Advancement of
  Science held at Southport in September 1883, pp.~407--408  (1884).

\bibitem{gauss_02}
\textit{C.~F. Gauss}, Disquisitiones generales circa seriem infinitam
  $1+\dfrac{\alpha\beta}{1\cdot\gamma}x+
  \dfrac{\alpha(\alpha+1)\beta(\beta+1)}{1\cdot2\cdot\gamma(\gamma+1)}xx+
  \dfrac{\alpha(\alpha+1)(\alpha+2)\beta(\beta+1)(\beta+2)}{1\cdot2\cdot3\cdot%
\gamma(\gamma+1)(\gamma+2)}x^3+\mathrm{etc}$, Commentationes Societatis Regiae
  Scientiarum Gottingensis recentiores, Classis Mathematic\ae, vol. II,
  pp.~3--46 [republished later in ``\text{Carl} \text{Friedrich} \text{Gauss}
  \text{Werke}'', vol.~3, pp.~265--327, \text{Königliche} \text{Gesellschaft}
  der \text{Wissenschaften}, \text{Göttingen}, 1866]  (1813).

\bibitem{gautschi_01}
\textit{W.~Gautschi}, Some elementary inequalities relating to the gamma and
  incomplete gamma function, Journal of Mathematics and Physics, no.~38,
  pp.~77--81  (1959).

\bibitem{gelfond_01}
\textit{A.~O. Gelfond}, The calculus of finite differences (3rd revised
  edition) [in Russian], Nauka, Moscow, USSR, 1967.

\bibitem{gessel_01}
\textit{I.~Gessel} and \textit{R.~P. Stanley}, \text{Stirling} polynomials,
  Journal of Combinatorial Theory, vol.~A24, pp.~24--33  (1978).

\bibitem{glaisher_02}
\textit{G.~W.~L. Glaisher}, Congruences relating to the sums of products of the
  first $n$ numbers and to other sums of products, The Quarterly journal of
  pure and applied mathematics, vol.~31, pp.~1--35  (1900).

\bibitem{godefroy_01}
\textit{M.~Godefroy}, La fonction Gamma ; Théorie, Histoire, Bibliographie,
  Gauthier--Villars, Imprimeur Libraire du Bureau des Longitudes, de l'École
  Polytechnique, Quai des Grands--Augustins, 55, Paris, 1901.

\bibitem{goldstine_01}
\textit{H.~H. Goldstine}, A History of Numerical Analysis from the 16th through
  the 19th Century, Springer--Verlag, New--York, Heidelberg, Berlin, 1977.

\bibitem{gould_01}
\textit{H.~W. Gould}, \text{Stirling} number representation problems,
  Proceedings of the American Mathematical Society, vol.~11, no.~3, pp.~447-451
   (1960).

\bibitem{gould_03}
\textit{H.~W. Gould}, An identity involving \text{Stirling} numbers, Annals of
  the Institute of Statistical Mathematics, vol.~17, \no.~1, pp.265--269
  (1965).

\bibitem{gould_02}
\textit{H.~W. Gould}, Note on recurrence relations for \text{Stirling} numbers,
  Publications de l'Institut Mathématique, Nouvelle série, vol.~6 (20),
  pp.~115--119  (1966).

\bibitem{knuth_01}
\textit{R.~L. Graham}, \textit{D.~E. Knuth} and \textit{O.~Patashnik}, Concrete
  mathematics: \text{A} foundation for computer science (2nd), Addison--Wesley,
  USA, 1994.

\bibitem{grunberg_01}
\textit{D.~B. Grünberg}, On asymptotics, \text{Stirling} numbers, gamma
  function and polylogs, Results in Mathematics, vol.~ 49, no.~1--2,
  pp.~89--125  (2006).

\bibitem{hagen_01}
\textit{J.~G. Hagen}, Synopsis der höheren Analysis. Vol.~1. Arithmetische und
  algebraische Analyse, von Felix L.~Dames, Taubenstraße 47, Berlin, Germany,
  1891.

\bibitem{hansen_01}
\textit{E.~R. Hansen}, A Table of Series and Products, Prentice--Hall, 1975.

\bibitem{hauss_01}
\textit{M.~Hauss}, Verallgemeinerte Stirling, Bernoulli und Euler Zahlen, deren
  Anwendungen und schnell konvergente Reihe für Zeta Funktionen
  (Ph.D.~dissertation), Aachen, Germany, 1995.

\bibitem{hayman_01}
\textit{W.~K. Hayman}, A generalisation of \text{Stirling's} formula, Journal
  für die reine und angewandte Mathematik, vol.~196, pp.~67--95  (1956).

\bibitem{hermite_01}
\textit{C.~Hermite}, Extrait de quelques lettres de \text{M.~Ch.~Hermite} à
  \text{M.~S.~Pincherle}, Annali di matematica pura ed applicata, serie III,
  tomo V, pp.~57--72  (1901).

\bibitem{hindenburg_01}
\textit{C.~F. Hindenburg}, Der polynomische Lehrsatz das wichtigste Theorem der
  ganzen Analysis nebst einigen Verwandten und andern Sätzen : Neu bearbeitet
  von \text{Tetens}, \text{Klügel}, \text{Kramp}, \text{Pfaff} und
  \text{Hindenburg}, bei Gerhard Fleischer, Leipzig, 1796.

\bibitem{howard_01}
\textit{F.~T. Howard}, Extensions of congruences of \text{Glaisher} and
  \text{Nielsen} concerning \text{Stirling} numbers, The Fibonacci Quarterly,
  vol.~28, no.~4, pp.~355--362  (1990).

\bibitem{howard_03}
\textit{F.~T. Howard}, Nörlund number \text{$B_n^{(n)}$}, in ``Applications of
  Fibonacci Numbers,'' vol.~5, pp.~355--366, Kluwer Academic, Dordrecht
  (1993).

\bibitem{howard_02}
\textit{F.~T. Howard}, Congruences and recurrences for \text{Bernoulli} numbers
  of higher orders, The Fibonacci Quarterly, vol.~32, no.~4, pp.~316--328
  (1994).

\bibitem{hwang_01}
\textit{H.~K. Hwang}, Asymptotic expansions for the \text{Stirling} numbers of
  the first kind, Journal of Combinatorial Theory, ser.~A 71, pp.~343--351
  (1995).

\bibitem{jeffreys_02}
\textit{H.~Jeffreys} and \textit{B.~S. Jeffreys}, Methods of mathematical
  physics (second edition), University Press, Cambridge, Great Britain, 1950.

\bibitem{jordan_02}
\textit{C.~Jordan}, Sur des polynômes analogues aux polynômes de
  \text{Bernoulli}, et sur des formules de sommation analogues à celle de
  \text{MacLaurin--Euler}, Acta Scientiarum Mathematicarum (Szeged), vol.~4,
  no.~3-3, pp.~130--150  (1928--1929).

\bibitem{jordan_00}
\textit{C.~Jordan}, On \text{Stirling's Numbers}, Tohoku Mathematical Journal,
  First Series, vol.~37, pp.~254--278  (1933).

\bibitem{jordan_01}
\textit{C.~Jordan}, The calculus of finite differences, Chelsea Publishing
  Company, USA, 1947.

\bibitem{kenter_01}
\textit{F.~K. Kenter}, A matrix representation for \text{Euler's} constant
  $\gamma$, The American Mathematical Monthly, vol.~106, pp.~452--454  (1999).

\bibitem{kluyver_02}
\textit{J.~C. Kluyver}, \text{Euler's} constant and natural numbers,
  Proc.~K.~Ned.~Akad. Wet., vol.~27, no.~1--2, pp.~142--144  (1924).

\bibitem{knopp_01}
\textit{K.~Knopp}, Theory and applications of infinite series (2nd edition),
  Blackie \& Son Limited, London and Glasgow, UK, 1951.

\bibitem{knuth_02}
\textit{D.~E. Knuth}, Two notes on notation, American Mathematical Monthly,
  vol.~99, no.~5, pp.~403--422  (1992).

\bibitem{korn_01}
\textit{G.~A. Korn} and \textit{T.~M. Korn}, Mathematical Handbook for
  Scientists and Engineers. Definitions, Theorems, and Formulas for Reference
  and Review (second, enlarged and revised edition), McGraw--Hill Book Company,
  New--York, 1968.

\bibitem{kowalenko_02}
\textit{V.~Kowalenko}, Generalizing the reciprocal logarithm numbers by
  adapting the partition method for a power series expansion, Acta
  Applicand\ae~Mathematic\ae, vol.~106, pp.~369--420  (2009).

\bibitem{kowalenko_01}
\textit{V.~Kowalenko}, Properties and applications of the reciprocal logarithm
  numbers, Acta Applicand\ae~Mathematic\ae, vol.~109, pp.~413--437  (2010).

\bibitem{kramp_01}
\textit{C.~Kramp}, Élements d'arithmétique universelle, L'imprimerie de
  Th.~F.~Thiriart, Cologne, 1808.

\bibitem{kratzer_01}
\textit{A.~Kratzer} and \textit{W.~Franz}, Transzendente Funktionen,
  Akademische Verlagsgesellschaft, Leipzig, Germany, 1960.

\bibitem{kruchinin_01}
\textit{V.~Kruchinin}, Composition of ordinary generating functions,
  arXiv:1009.2565  (2010).

\bibitem{kruchinin_03}
\textit{V.~V. Kruchinin} and \textit{D.~V. Kruchinin}, Powers of generating
  functions and their applications [in Russian], Tomsk University Press,
  Russia, 2013.

\bibitem{kruchinin_02}
\textit{V.~V. Kruchinin} and \textit{D.~V. Kruchinin}, Composita and its
  properties, Journal of Analysis \& Number Theory, vol.~2, no.~2, pp.~37--44
  (2014, enlarged preprint published in March 2011, arXiv:1103.2582).

\bibitem{krylov_01}
\textit{V.~I. Krylov}, Approximate calculation of integrals, The Macmillan
  Company, New-York, USA, 1962.

\bibitem{skramer_01}
\textit{S.~Krämer}, Die Eulersche Konstante $\gamma$ und verwandte Zahlen
  (unpublished Ph.D.~manuscript, pers.~comm.), Göttingen, Germany, 2014.

\bibitem{kuznecov_01}
\textit{D.~S. Kuznetsov}, Special functions (2nd edition) [in Russian],
  Vysshaya Shkola, Moscow, 1965.

\bibitem{lanczos_01}
\textit{C.~Lanczos}, A precision approximation of the gamma function, SIAM
  Journal on Numerical Analysis, vol.~1, pp.~86--96  (1964).

\bibitem{laplace_01}
\textit{P.-S. Laplace}, Traité de Mécanique Céleste. Tomes I--V, Courcier,
  Imprimeur--Libraire pour les Mathématiques, quai des Augustins, \no 71,
  Paris, France, 1800--1805.

\bibitem{laplace_02}
\textit{P.-S. Laplace}, Traité analytique des probabiblités, M\up{me} V\up{e}
  Courcier, Imprimeur--Libraire pour les Mathématiques, quai des Augustins,
  \no 57, Paris, France, 1812.

\bibitem{lienard_01}
\textit{R.~Li\'enard}, Nombres de \text{Cauchy}, Intermédiaire des Recherches
  Mathématiques, vol.~2, no.~5, p.~38  (1946).

\bibitem{louchard_01}
\textit{G.~Louchard}, Asymptotics of the \text{Stirling} numbers of the first
  kind revisited: \text{A} saddle point approach, Discrete Mathematics and
  Theoretical Computer Science, vol.~12, no.~2, pp.~167--184  (2010).

\bibitem{malmsten_01}
\textit{C.~J. Malmstén}, De integralibus quibusdam definitis seriebusque
  infinitis, Journal für die reine und angewandte Mathematik, vol.~38,
  pp.~1--39  (1849, work dated at May 1, 1846).

\bibitem{mascheroni_01}
\textit{L.~Mascheronio}, Adnotationes ad calculum integralem \text{Euleri} in
  quibus nonnulla problemata ab \text{Eulero} proposita resolvuntur, Ex
  Typographia Petri Galeatii, Ticini, 1790.

\bibitem{merlini_01}
\textit{D.~Merlini}, \textit{R.~Sprugnoli} and \textit{\text{M.~Cecilia}
  Verri}, The \text{Cauchy} numbers, Discrete Mathematics (Elsevier), vol.~306,
  pp.~1906--1920  (2006).

\bibitem{mezo_01}
\textit{I.~Mez\H{o}}, Gompertz constant, \text{Gregory} coefficients and a
  series of the logarithm function, Journal of Analysis \& Number Theory,
  vol.~2, no.~2, pp.~33--36  (2014).

\bibitem{mitrinovic_01}
\textit{D.~S. Mitrinovi\'c} and \textit{R.~S. Mitrinovi\'c}, Sur les nombres de
  \text{Stirling} et les nombres de \text{Bernoulli} d'ordre supérieur,
  Publications de la faculté d'électrotechnique de l'Université à
  Bélgrade, Série Mathématique et Physique, no. 43, pp.~1--63  (1960).

\bibitem{mortici_01}
\textit{C.~Mortici}, An ultimate extremely accurate formula for approximation
  of the factorial function, Archiv der Mathematik, vol.~93, pp.~37--45
  (2009).

\bibitem{moser_01}
\textit{L.~Moser} and \textit{M.~Wyman}, Asymptotic development of the
  \text{Stirling} numbers of the first kind, Journal of the London Mathematical
  Society, vol.~s1--33, no.~2, pp.~133--146  (1958).

\bibitem{murray_01}
\textit{F.~J. Murray}, Formulas for factorial \text{$N$}, Mathematics of
  Computation, vol.~39, vol.~160, , pp.~655--662  (1982).

\bibitem{nemes_02}
\textit{G.~Nemes}, On the coefficients of the asymptotic expansion of
  \text{$n!$}, Journal of Integer Sequences, vol.~13, article 10.6.6  (2010).

\bibitem{nemes_01}
\textit{G.~Nemes}, An asymptotic expansion for the \text{Bernoulli} numbers of
  the second kind, Journal of Integer Sequences, vol.~14, article 11.4.8
  (2011).

\bibitem{nemes_03}
\textit{G.~Nemes}, Generalization of \text{Binet's} gamma function formulas,
  Integral Transforms and Special Functions, vol.~24, no.~8, pp.~597--606
  (2013).

\bibitem{netto_01}
\textit{E.~Netto}, Lehrbuch der Combinatorik (2nd Edn.), Teubner, Leipzig,
  Germany, 1927.

\bibitem{nielsen_04}
\textit{N.~Nielsen}, Recherches sur les polynômes et les nombres de
  \text{Stirling}, Annali di Matematica Pura ed Applicata, vol.~10, no.~1,
  pp.~287--318  (1904).

\bibitem{nielsen_01}
\textit{N.~Nielsen}, Handbuch der \text{Theorie} der \text{Gammafunktion}, B.
  G. Teubner, Leipzig, Germany, 1906.

\bibitem{nielsen_03}
\textit{N.~Nielsen}, Recherches sur les polynômes de \text{Stirling},
  Hovedkommissionaer: Andr. Fred. Høst \& Søn, Kgl. Hof-Boghandel, Bianco
  Lunos Bogtrykkeri, København, Denmark, 1920.

\bibitem{norlund_02}
\textit{N.~E. Nörlund}, Vorlesungen über \text{Differenzenrechnung},
  Springer, Berlin, 1924.

\bibitem{norlund_01}
\textit{N.~E. Nörlund}, Sur les valeurs asymptotiques des nombres et des
  polynômes de \text{Bernoulli}, Rendiconti del Circolo Matematico di Palermo,
  vol.~10, no.~1, pp.~27--44  (1961).

\bibitem{olson_01}
\textit{F.~R. Olson}, Arithmetic properties of \text{Bernoulli} numbers of
  higher order, Duke Mathematical Journal, vol.~22, no.~4, pp.~641--653.
  (1955).

\bibitem{olver_01}
\textit{F.~W.~J. Olver}, Asymptotics and Special Functions, Academic Press,
  USA, 1974.

\bibitem{paris_01}
\textit{R.~B. Paris}, Asymptotic approximations for \text{$n!$}, Applied
  Mathematical Sciences, Vol.~5, no.~37, pp.~1801--1807  (2011).

\bibitem{polya_01_eng}
\textit{G.~P\'olya} and \textit{G.~Szeg\H{o}}, Problems and Theorems in
  Analysis I: Series, Integral calculus, Theory of functions, Springer--Verlag,
  Berlin, Germany, 1978.

\bibitem{proskuriyakov_01_eng}
\textit{I.~V. Proskuriyakov}, A Collection of Problems in Linear Algebra
  (fourth edition) [in Russian], Nauka, Moscow, USSR, 1970.

\bibitem{prudnikov_en}
\textit{A.~P. Prudnikov}, \textit{Y.~A. Brychkov} and \textit{O.~I. Marichev},
  Integrals and Series. Vol. I--IV, Gordon and Breach Science Publishers, 1992.

\bibitem{qi_01}
\textit{F.~Qi}, An integral representation, complete monotonicity, and
  inequalities of \text{Cauchy} numbers of the second kind, Journal of Number
  Theory, vol.~144, pp.~244--255  (2014).

\bibitem{rigaud_01}
\textit{S.~J. Rigaud}, Correspondence of scientific men of the seventeenth
  century, including letters of Barrow, Flamsteed, Wallis, and Newton, printed
  from the Originals [in 2 vols.], Oxford at the University Press, 1841.

\bibitem{riordan_01}
\textit{J.~Riordan}, An Introduction to Combinatorial Analysis, John Wiley \&
  Sons, Inc., USA, 1958.

\bibitem{rubinstein_01}
\textit{M.~O. Rubinstein}, Identities for the \text{Riemann} zeta function, The
  Ramanujan Journal, vol.~27, pp.~29--42  (2012).

\bibitem{rubinstein_02}
\textit{M.~O. Rubinstein}, Identities for the \text{Hurwitz} zeta function,
  \text{Gamma} function, and \text{$L$}-functions, The Ramanujan Journal,
  vol.~32, pp.~421--464  (2013).

\bibitem{rzadkowski_01}
\textit{G.~Rz\c{a}dkowski}, Two formulas for successive derivatives and their
  applications, Journal of Integer Sequences [electronic only], vol.~12,
  article 09.8.2  (2009).

\bibitem{salmieri_01}
\textit{A.~Salmeri}, Introduzione alla teoria dei coefficienti fattoriali, from
  ``Giornale di Matematiche di Battaglini'', vol.~90 (no.~10, serie 5),
  pp.~44--54  (1962).

\bibitem{sato_01}
\textit{H.~Sato}, On a relation between the \text{Riemann} zeta function and
  the \text{Stirling} numbers, Integers: Electronic Journal of Combinatorial
  Number Theory, vol.~8, no.~1  (2008).

\bibitem{schlaffli_01}
\textit{L.~Schläffli}, Sur les coëfficients du développement du produit
  $(1+x)(1+2x)\cdots \big(1+(n-1)\big)$ suivant les puissances ascendantes de
  $x$, Journal für die reine und angewandte Mathematik, vol.~43, pp.~1--22
  (1852).

\bibitem{schlaffli_02}
\textit{L.~Schläffli}, Ergänzung der abhandlung über die entwickelung des
  products $(1+x)(1+2x)\cdots \big(1+(n-1)\big)$ in band \text{XLIII} dieses
  journals, Journal für die reine und angewandte Mathematik, vol.~67,
  pp.~179--182  (1867).

\bibitem{schlomilch_04}
\textit{O.~Schlömilch}, Recherches sur les coefficients des facultés
  analytiques, Journal für die reine und angewandte Mathematik, vol.~44,
  pp.~344--355  (1852).

\bibitem{schlomilch_05}
\textit{O.~Schlömilch}, Compendium der höheren Analysis, Druck und Verlag von
  Friedrich Vieweg und Sohn, Braunschweig, Germany, 1853.

\bibitem{schlomilch_06}
\textit{O.~Schlömilch}, Compendium der höheren Analysis (2nd edn., in two
  volumes), Druck und Verlag von Friedrich Vieweg und Sohn, Braunschweig,
  Germany, 1861, 1866.

\bibitem{schlomilch_03}
\textit{O.~Schlömilch}, Nachschrift hierzu, Zeitschrift für angewandte
  Mathematik und Physik, vol.~25, pp.~117--119  (1880).

\bibitem{schroder_01}
\textit{E.~Schröder}, Bestimmung des infinitären \text{Werthes} des
  \text{Integrals} $\int\limits_0^1 (u)_n\, du$, Zeitschrift für angewandte
  Mathematik und Physik, vol.~25, pp.~106--117  (1880).

\bibitem{ser_01}
\textit{J.~Ser}, Sur une expression de la fonction $\zeta(s)$ de
  \text{Riemann}, Comptes-rendus hebdomadaires des séances de l'Académie des
  Sciences, Série 2, vol.~182, pp.~1075--1077  (1926).

\bibitem{shen_01}
\textit{L.-C. Shen}, Remarks on some integrals and series involving the
  \text{Stirling} numbers and $\zeta(n)$, The Transactions of the American
  Mathematical Society, vol.~347, no.~4, pp.~1391--1399  (1995).

\bibitem{shirai_01}
\textit{S.~Shirai} and \textit{K.~ichi Sato}, Some identities involving
  \text{Bernoulli} and \text{Stirling} numbers, Journal of Number Theory,
  vol.~90, pp.~130--142  (2001).

\bibitem{sondow_02}
\textit{J.~Sondow}, Double integrals for \text{Euler's} constant and
  $\ln(4/\pi)$ and an analog of \text{Hadjicostas's} formula, American
  Mathematical Monthly, vol.~112, pp.~61--65  (2005).

\bibitem{spouge_01}
\textit{J.~L. Spouge}, Computation of the gamma, digamma, and trigamma
  functions, SIAM Journal on Numerical Analysis, vol.~31, no.~3, pp.~931--944
  (1994).

\bibitem{srivastava_03}
\textit{H.~M. Srivastava} and \textit{J.~Choi}, Series Associated with the Zeta
  and Related Functions, Kluwer Academic Publishers, the Netherlands, 2001.

\bibitem{stamper_01}
\textit{P.~C. Stamper}, Table of \text{Gregory} coefficients, Mathematics of
  Computation, vol.~20, p.~465  (1966).

\bibitem{stanley_01}
\textit{R.~P. Stanley}, Enumerative Combinatorics (1st Edn., 2nd printing),
  Cambridge University Press, 1997.

\bibitem{steffensen_00}
\textit{J.~F. Steffensen}, On certain formulas of approximate summation and
  integration, Journal of the Institute of Actuaries, vol.~53, no.~2,
  pp.~192--201  (1922).

\bibitem{steffensen_01}
\textit{J.~F. Steffensen}, On \text{Laplace's} and \text{Gauss'}
  summation--formulas, Skandinavisk Aktuarietidskrift (Scandinavian Actuarial
  Journal), no.~1, pp.~1--15  (1924).

\bibitem{steffensen_02}
\textit{J.~F. Steffensen}, Interpolation (2nd Edn.), Chelsea Publishing
  Company, New--York, USA, 1950.

\bibitem{stirling_01}
\textit{J.~Stirling}, Methodus differentialis, sive Tractatus de summatione et
  interpolatione serierum infinitarum, Gul.~Bowyer, Londini, 1730.

\bibitem{temme_02}
\textit{N.~M. Temme}, Asymptotic estimates of \text{Stirling} numbers, Studies
  in Applied Mathematics, vol.~89, pp.~233--243  (1993).

\bibitem{demoivre_01}
\textit{A.~\text{De Moivre}}, Miscellanea analytica de seriebus et quadraturis
  (with a supplement of 21 pages), J.~Thonson \& J.~Watts, Londini, 1730.

\bibitem{gunter_02_eng}
\textit{\text{N.~M.~Gunther (Günter)}} and \textit{\text{R.~O.~Kuzmin
  (Kusmin)}}, A Collection of Problems on Higher Mathematics. Vol.~2 (12th
  edition) [in Russian], Gosudarstvennoe izdatel'stvo tehniko--teoreticheskoj
  literatury, Leningrad, USSR, 1949.

\bibitem{gunter_03_eng}
\textit{\text{N.~M.~Gunther (Günter)}} and \textit{\text{R.~O.~Kuzmin
  (Kusmin)}}, A Collection of Problems on Higher Mathematics. Vol.~3 (4th
  edition) [in Russian], Gosudarstvennoe izdatel'stvo tehniko--teoreticheskoj
  literatury, Leningrad, USSR, 1951.

\bibitem{van_veen_01}
\textit{S.~C. \text{Van Veen}}, Asymptotic expansion of the generalized
  \text{Bernoulli} numbers \text{$B_n^{(n-1)}$} for large values of $n$ ($n$
  integer), Indagationes Mathematic\ae~(Proceedings of the Koninklijke
  Nederlandse Akademie van Wetenschappen. Series A, Mathematical sciences),
  vol.~13, pp.~335--341  (1951).

\bibitem{timashev_01}
\textit{A.~N. Timashev}, On asymptotic expansions of \text{Stirling} numbers of
  the first and second kinds, Discrete Mathematics and Applications, vol.~8,
  no.~5, pp.~533--544  (1998).

\bibitem{tricomi_01}
\textit{F.~G. Tricomi} and \textit{A.~Erdélyi}, The asymptotic expansion of a
  ratio of gamma functions, Pacific Journal of Mathematics, vol.~1, no.~1,
  pp.~133--142  (1951).

\bibitem{newton_01}
\textit{H.~W. Turnbull}, The correspondence of Isaac Newton [vols.~1--7], Royal
  Society at the University Press, Cambridge, 1959--1977.

\bibitem{tweedie_01}
\textit{C.~Tweedie}, The \text{Stirling} numbers and polynomials, Proceedings
  of the Edinburgh Mathematical Society, vol.~37, pp.~2--25  (1918).

\bibitem{vorobiev_01}
\textit{N.~N. Vorobiev}, Theory of series (4th edition, enlarged and revised)
  [in Russian], Nauka, Moscow, USSR, 1979.

\bibitem{wachs_01}
\textit{S.~Wachs}, Sur une propriété arithmétique des nombres de
  \text{Cauchy}, Bulletin des Sciences Mathématiques, deuxi\`eme s\'erie.
  vol.~71, pp.~219--232  (1947).

\bibitem{watson_01}
\textit{G.~N. Watson}, An expansion related to \text{Stirling's} formula,
  derived by the method of steepest descents, The Quarterly journal of pure and
  applied mathematics, vol.~48, pp.~1--18  (1920).

\bibitem{weisstein_04}
\textit{E.~W. Weisstein}, CRC Concise Encyclopedia of Mathematics (2nd Edn.),
  Chapman \& Hall/CRC, USA, 2003.

\bibitem{whittaker_01}
\textit{E.~Whittaker} and \textit{G.~N. Watson}, A course of modern analysis.
  An introduction to the general theory of infinite processes and of analytic
  functions, with an account of the principal transcendental functions (third
  edition), Cambridge at the University Press, Great Britain, 1920.

\bibitem{wilf_02}
\textit{H.~S. Wilf}, The asymptotic behavior of the \text{Stirling} numbers of
  the first kind, Journal of Combinatorial Theory, ser.~A 64, pp.~344--349
  (1993).

\bibitem{wilf_01}
\textit{H.~S. Wilf}, Generatingfunctionology (2nd), Academic Press, Inc., USA,
  1994.

\bibitem{wilton_02}
\textit{J.~R. Wilton}, A proof of \text{Burnside's} formula for
  \text{$\log\Gamma(x+1)$} and certain allied properties of \text{Riemann's}
  $\zeta$-function, The Messenger of Mathematics, vol.~52, pp.~90--93
  (1922--1923).

\bibitem{wrench_01}
\textit{J.~W. Wrench}, Concerning two series for the gamma function,
  Mathematics of Computation, vol.~22, pp.~617--626  (1968).

\bibitem{young_02}
\textit{P.~T. Young}, Congruences for \text{Bernoulli}, \text{Euler} and
  \text{Stirling} numbers, Journal of Number Theory, vol.~78, pp.~204--227
  (1999).

\bibitem{young_01}
\textit{P.~T. Young}, A 2-adic formula for \text{Bernoulli} numbers of the
  second kind and for the \text{Nörlund} numbers, Journal of Number Theory,
  vol.~128, pp.~2951--2962  (2008).

\bibitem{young_03}
\textit{P.~T. Young}, Rational series for multiple zeta and log gamma
  functions, Journal of Number Theory, vol.~133, pp.~3995--4009  (2013).

\bibitem{zhao_01}
\textit{F.-Z. Zhao}, Sums of products of \text{Cauchy} numbers, Discrete
  Mathematics, vol.~309, pp.~3830--3842  (2009).

\end{thebibliography}

\newpage
\begin{figure}[!t]   
\hspace{-6mm}
\includegraphics[width=1.05\textwidth]{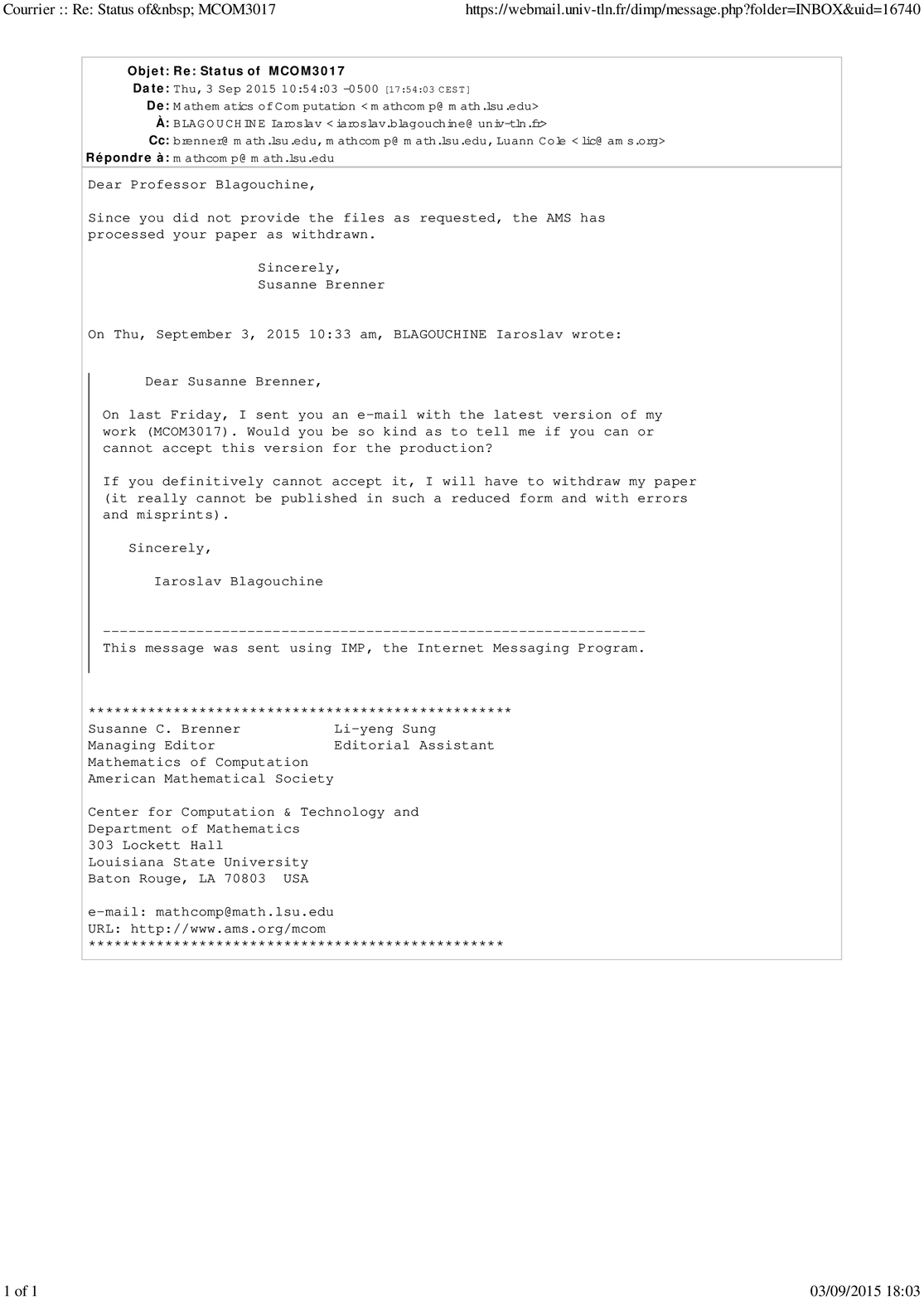}
\end{figure}

\newpage
\begin{figure}[!t]   
\hspace{-6mm}
\includegraphics[width=1.05\textwidth]{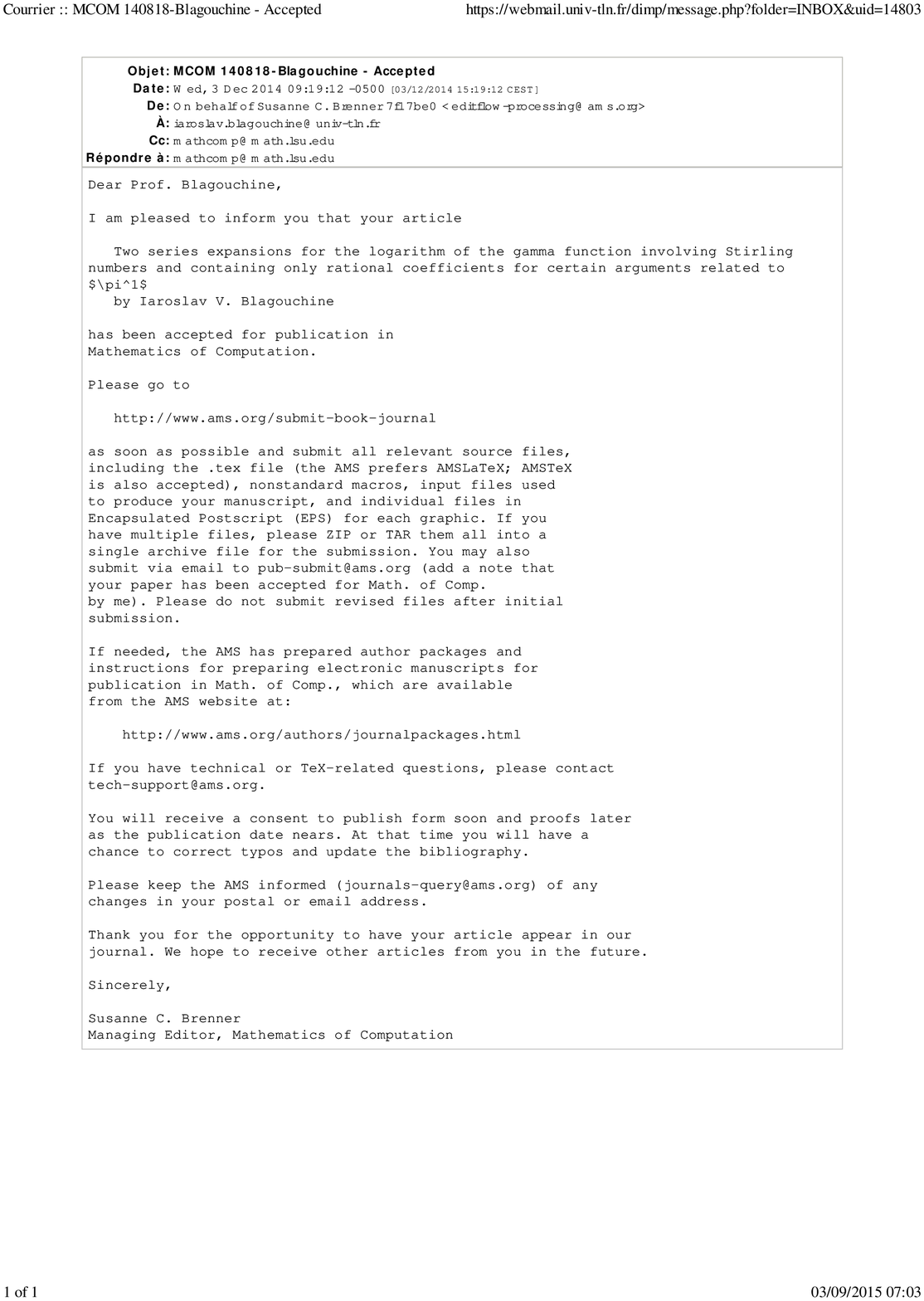}
\end{figure}

\newpage
\begin{figure}[!t]   
\hspace{-6mm}
\includegraphics[width=1.05\textwidth]{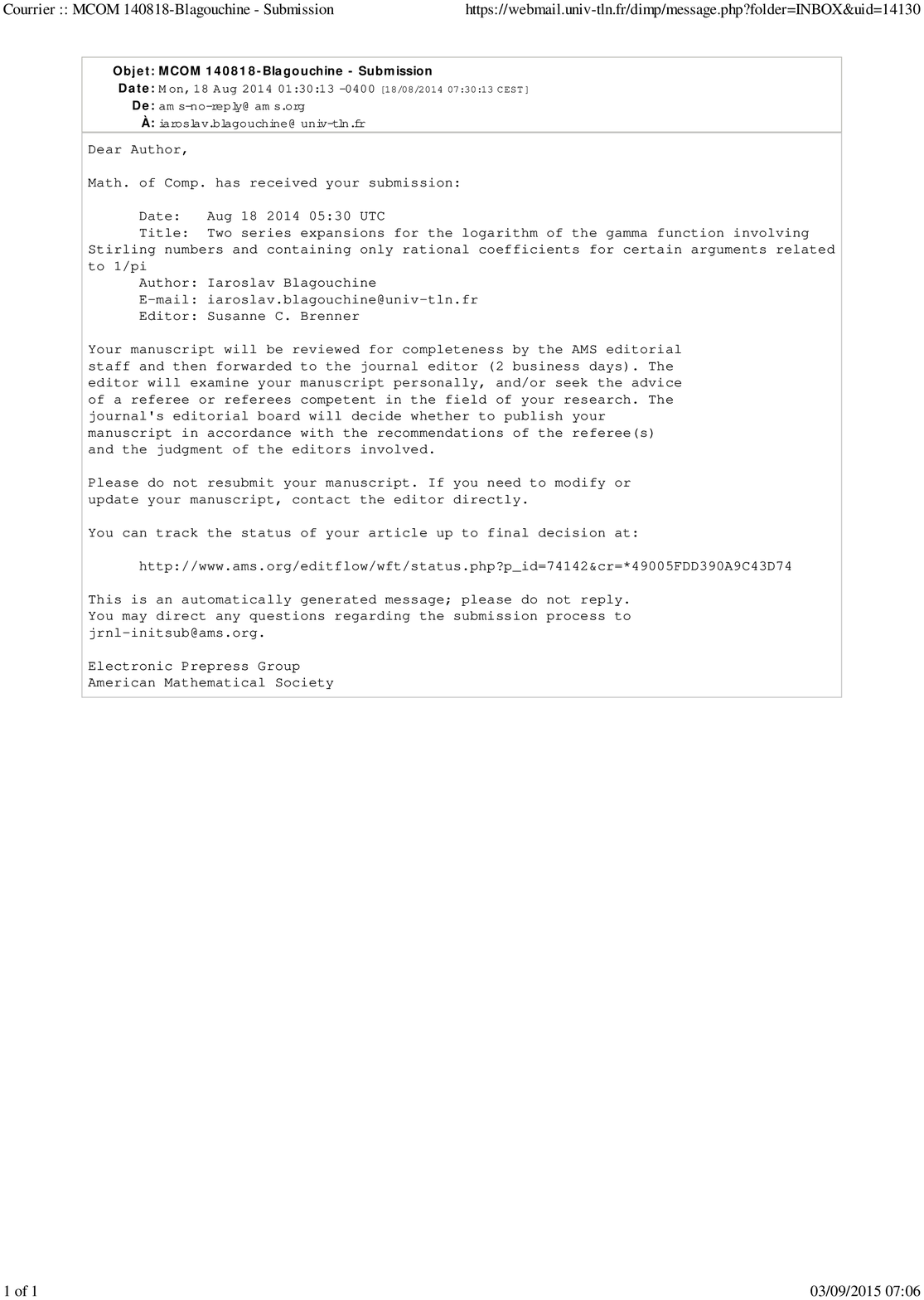}
\end{figure}

\end{document}